\documentclass[11pt]{article}
\usepackage[letterpaper,margin=1in]{geometry}
\usepackage[english]{babel}
\usepackage{amsmath}
\usepackage{amssymb}
\usepackage{amsthm}
\usepackage{bbm,bm}
\usepackage{url}
\usepackage{pdfpages}
\usepackage{float}
\usepackage{bookmark}
\usepackage{subfigure}
\usepackage{mathtools}
\usepackage{colonequals}
\usepackage{siunitx}

\usepackage{xcolor}

\definecolor{cblack}{rgb}{0,0,0}
\definecolor{cblue}{rgb}{0.121569,0.466667,0.705882}
\definecolor{corange}{rgb}{1.000000,0.498039,0.054902}
\definecolor{cgreen}{rgb}{0.172549,0.627451,0.172549}
\definecolor{cred}{rgb}{0.839216,0.152941,0.156863}
\definecolor{cpurple}{rgb}{0.580392,0.403922,0.741176}
\definecolor{cbrown}{rgb}{0.549020,0.337255,0.294118}
\definecolor{cpink}{rgb}{0.890196,0.466667,0.760784}
\definecolor{cgray}{rgb}{0.498039,0.498039,0.498039}
\definecolor{cgreen2}{rgb}{0.7372549019607844, 0.7411764705882353, 0.13333333333333333}

\definecolor{clightgray}{rgb}{0.6,0.6,0.6}
\definecolor{cllightgray}{rgb}{0.8,0.8,0.8}

\usepackage{hyperref}
\hypersetup{
  pdfencoding=auto,
  linkcolor  = cblue,
  citecolor  = cgreen,
  urlcolor   = corange,
  colorlinks = true,
}

\newtheorem{thm}{Theorem}[section]

\newtheorem{definition}[thm]{Definition}
\newtheorem{lem}[thm]{Lemma}
\newtheorem{prop}[thm]{Proposition}

\newtheorem{cor}[thm]{Corollary}

\newtheorem{rmk}[thm]{Remark}
\newtheorem{assumption}[thm]{Assumption}

\numberwithin{equation}{section}

\newcommand{\bg}{\bm{g}}

\newcommand{\bw}{\bm{w}}
\newcommand{\bx}{\bm{x}}
\newcommand{\by}{\bm{y}}
\newcommand{\bz}{\bm{z}}

\newcommand{\bkappa}{\bm{\kappa}}
\newcommand{\bGamma}{\bm{\Gamma}}
\newcommand{\bsigma}{\bm{\sigma}}

\newcommand{\sE}{\mathcal E}
\newcommand{\sO}{\mathcal O}
\newcommand{\sU}{\mathcal U}

\newcommand{\CC}{\mathbb{C}}
\newcommand{\EE}{\mathbb{E}}
\newcommand{\FF}{\mathbb{F}}

\newcommand{\ZZ}{\mathbb{Z}}
\newcommand{\PP}{\mathbb{P}}
\newcommand{\RR}{\mathbb{R}}
\renewcommand{\SS}{\mathbb{S}}

\newcommand{\bA}{\mathbf{A}}
\newcommand{\bB}{\mathbf{B}}

\newcommand{\bE}{\mathbf{E}}

\newcommand{\bQ}{\mathbf{Q}}

\newcommand{\bW}{\mathbf{W}}
\newcommand{\bX}{\mathbf{X}}
\newcommand{\bY}{\mathbf{Y}}

\newcommand{\sA}{\mathcal{A}}
\newcommand{\sB}{\mathcal{B}}
\newcommand{\sC}{\mathcal{C}}

\newcommand{\sF}{\mathcal{F}}

\newcommand{\sN}{\mathcal{N}}

\newcommand{\sS}{\mathcal{S}}

\newcommand{\sX}{\mathcal{X}}

\newcommand{\sym}{\mathrm{sym}}
\newcommand{\GOE}{\mathrm{GOE}}
\newcommand{\GUE}{\mathrm{GUE}}
\newcommand{\GFE}{\mathrm{G}\FF\mathrm{E}}

\newcommand{\Res}{\mathrm{Res}}
\newcommand{\herm}{\mathrm{herm}}

\newcommand{\spec}{\mathrm{spec}}
\newcommand{\Round}{\mathrm{Round}}

\newcommand{\Haar}{\mathrm{Haar}}

\renewcommand{\sc}{\mathrm{sc}}

\newcommand{\sgn}{\mathrm{sgn}}

\newcommand{\diag}{\mathrm{diag}}
\newcommand{\offdiag}{\mathrm{offdiag}}

\newcommand{\F}{\mathrm{F}}

\newcommand{\Diag}{\mathrm{Diag}}
\newcommand{\id}{\mathrm{id}}
\newcommand{\conj}{\mathrm{conj}}
\newcommand{\cpxi}{\bm i}

\newcommand{\reg}{\mathrm{reg}}

\renewcommand{\emptyset}{\varnothing}

\newcommand{\eqlaw}{\stackrel{\text{(law)}}{=}}

\newcommand{\boldone}{\mathbf{1}}

\renewcommand{\Im}{\mathrm{Im}}
\renewcommand{\Re}{\mathrm{Re}}

\newcommand{\what}{\widehat}
\newcommand{\wtilde}{\widetilde}

\newcommand\numberthis{\addtocounter{equation}{1}\tag{\theequation}}

\renewcommand{\epsilon}{\varepsilon}

\newcommand{\Ex}{\mathop{\mathbb{E}}}

\DeclareMathOperator*{\argmin}{arg\,min}

\usepackage{authblk}

\title{Universal entrywise eigenvector fluctuations in delocalized spiked matrix models and asymptotics of rounded spectral algorithms}
\date{August 20, 2026}
\author{Shujing Chen\thanks{Email: \texttt{schen344@jhu.edu}}\,\,}
\author{Dmitriy Kunisky\thanks{Email: \texttt{kunisky@jhu.edu}.}}
\affil{Department of Applied Mathematics \& Statistics, Johns Hopkins University}

\begin{document}

\pagenumbering{roman}

\maketitle
\thispagestyle{empty}

\begin{abstract}
    We consider the distribution of the top eigenvector $\widehat{v}$ of a spiked matrix model of the form $H = \theta vv^* + W$, in the supercritical regime where $H$ has an outlier eigenvalue of comparable magnitude to $\|W\|$.
    We show that, if $v$ is sufficiently delocalized, then the distribution of the individual entries of the projector $\widehat{v}\widehat{v}^*$ (not, we emphasize, merely the inner product $|\langle \widehat{v}, v\rangle|^2$) is universal over a large class of generalized Wigner matrices $W$ having independent entries, depending only on the first two moments of the distributions of the entries of $W$.
    This complements the observation of Capitaine and Donati-Martin (2021) that these distributions are \emph{not} universal when $v$ is instead sufficiently localized.
    Further, for $W$ having entrywise variances close to constant and thus resembling a Wigner matrix, we show by comparing to $W$ drawn from the Gaussian orthogonal or unitary ensembles that averages of entrywise functions of $\widehat{v}\widehat{v}^*$ behave as they would if $\widehat{v}$ had Gaussian fluctuations around a suitable multiple of $v$.
    We also establish such results for several possibly dependent spiked matrices, showing that, if such matrices are entrywise uncorrelated, then their leading eigenvectors behave as they would with independent Gaussian fluctuations.
    We apply these results to spectral algorithms with rounding procedures for synchronization problems over the cyclic and circle groups, obtaining the first precise asymptotic error rates for such algorithms.
    Using our analysis of multiple spiked matrices, we also show that \emph{multi-frequency} spectral algorithms using estimates from several matrices often have asymptotic error rate superior to that of naive spectral algorithms using just one matrix.
\end{abstract}

\clearpage

\pagestyle{empty}
\tableofcontents

\clearpage

\pagestyle{plain}
\pagenumbering{arabic}

\section{Introduction}

We will study a universality phenomenon in \emph{spiked matrix models}, of which we take the simple rank-one case
\[ H = \theta vv^{*} + W, \]
where $\theta \in \RR$, $v \in \SS^{n - 1}(\CC) \subset \CC^n$ (the unit sphere of $\CC^n$), and $W \in \CC^{n \times n}_{\herm}$ is a Hermitian ``noise'' matrix normalized such that $\|W\|$ is of constant order.
From a statistical point of view, $H$ can be seen as an ``observation'' of the hidden ``signal'' $v$, from which one seeks to estimate $v$.
A natural choice of estimator is $\what{v} = v_1(H)$, the unit eigenvector of $H$ corresponding to the largest eigenvalue of $H$, which we denote $\what{\lambda} = \lambda_1(H)$.
Our general goal will be to understand in a fine-grained way the quality of estimate of $v$ that we obtain from $\what{v}$.

Let us explain some of what is known about such models under the following assumption on $W$, which we will also use in the first part of our main results.
\begin{definition}[Generalized Wigner matrix]
    \label{def:gen-wigner}
    A \emph{generalized Wigner matrix} with $W \in \CC^{n \times n}_{\herm}$ with parameters $(\gamma_W, \xi_W)$ is a random matrix having the following properties:
\begin{enumerate}
    \item The entries on and above the diagonal $(W_{ij})_{1 \leq i \leq j \leq n}$ are independent.
    \item $\EE W_{ij} = 0$ for all $i, j \in [n]$.
    \item The magnitudes $|W_{ij}|$ admit a uniform tail bound of the form
    \begin{equation}
        \label{eq:entries-subexp}
        \PP[\sqrt{n} \cdot |W_{ij}| \ge t] \leq \xi_W^{-1} \exp(-t^{\xi_W})
    \end{equation}
    We will call this condition the rescaled entries $\sqrt{n} \cdot |W_{ij}|$ being \emph{$\xi_W$-sub-Weibull} (Definition~\ref{def:power-subexp}).
    \item The numbers $\sigma_{ij}^2 \colonequals \EE |W_{ij}|^2$ satisfy
    \begin{equation}
    \label{eq:variance-bounds}
    \gamma_W^{-1} \leq n \cdot \sigma_{ij}^2 \leq \gamma_W \text{ for } i\neq j,\ i, j \in [n],
    \end{equation}
    as well as
    \begin{equation}
    \label{eq:gen-wigner-var-sums}
    \sum_{j = 1}^n \sigma_{ij}^2 = 1 \text{ for all } i \in [n].
    \end{equation}
\end{enumerate}
\end{definition}
\noindent
This class of matrices is a well-studied \emph{universality class}, sharing many basic behaviors.
For instance, as a useful point of reference for the results we discuss next, under these assumptions the empirical distribution of $W$ almost surely converges weakly to the semicircle distribution supported on $[-2, 2]$, and $\lambda_1(W)$ and $\|W\|$ both converge in probability to 2 (see, e.g., \cite{Anderson_Guionnet_Zeitouni_2009} for an exposition of these classical results, or our Theorem~\ref{thm:semicircle}).

We mostly take $v$ to be deterministic, and later in our results we will mention the case of $v$ random, in which case it is random independently of $W$ (in this situation results for $v$ deterministic can be applied by conditioning on the value of $v$).

The following describes the main features of the \emph{phase transition} that the behavior of the spectrum of $H$ undergoes as $\theta$ varies, a version for this setting of the \emph{Baik--Ben Arous--P\'{e}ch\'{e}} transition.
We sketch how the proof of this version follows from the \emph{isotropic local law} proved by \cite{bloemendal2015isotropiclocallawssample} in Theorem~\ref{thm:spiked-non-asymp} below; this technique has been used by several works such as \cite{knowles2012isotropicsemicirclelawdeformation,Knowles_2014} in the past and we merely adapt the particular local law on which it is based.
\begin{thm}
    \label{thm:bbp}
Let $W = W^{(n)} \in \CC^{n \times n}_{\herm}$ be a sequence of generalized Wigner matrices with parameters $(\gamma_W, \xi_W)$ not depending on $n$ and let $v = v^{(n)} \in \SS^{n - 1}(\CC)$.
Write $H = H^{(n)} = \theta vv^{*} + W$, $\what{\lambda} = \what{\lambda}^{(n)} = \lambda_1(H^{(n)})$ for the largest eigenvalues of these matrices, and $\what{v} = \what{v}^{(n)}$ for the associated eigenvectors of $H$.
Then, the following hold, with all convergences in probability as $n \to \infty$:
\begin{align}
    \what{\lambda} &\to \lambda(\theta) \colonequals \left\{\begin{array}{ll} 2 & \text{if } \theta \leq 1 \\ \theta + \theta^{-1} & \text{if } \theta > 1\end{array}\right\}, \label{eq:def-lambda} \\
    |\langle \what{v}, v \rangle|^2 &\to \rho(\theta)^2 \colonequals \left\{\begin{array}{ll} 0 & \text{if } \theta \leq 1 \\ 1 - \theta^{-2} & \text{if } \theta > 1\end{array}\right\}. \label{eq:def-rho}
\end{align}
\end{thm}
\noindent
It will also be useful for us to have a specific notation for
\begin{equation}
\label{eq:def-tau}
\tau(\theta)^2 \colonequals 1 - \rho(\theta)^2 = \left\{\begin{array}{ll} 1 & \text{if } \theta \leq 1 \\ \theta^{-2} & \text{if } \theta > 1\end{array}\right\},
\end{equation}
the typical value of $1 - |\langle \what{v}, v \rangle|^2$, the magnitude of the component of $\what{v}$ orthogonal to $v$.

This gives a precise understanding of the top eigenvalue $\what{\lambda}$ and eigenvector $\what{v}$ of $H$ to leading order, though for the eigenvector only in terms of its correlation with its noiseless counterpart $v$.
For both quantities, it is reasonable to ask about the next-order fluctuations: how do a rescaling of $\what{\lambda} - \lambda(\theta)$ and $\what{v} - \rho(\theta)v$ behave?
(The latter is not quite well-defined for the reason described in Remark~\ref{rmk:outer-prod} below, but let us ignore this issue for the moment.)

While it would be natural to conjecture that both of these quantities should also behave universally over a generalized Wigner matrices or some other such large class, in fact this is not the case.
In our setting, it is known that \emph{non-universality} of the fluctuations of both eigenvalues and eigenvectors---that is, a sensitivity to the specific distribution of entries of $W$---arises when $v$ is \emph{localized}, having some large entries.
See, for instance, \cite{Capitaine_2009,capitaine2011centrallimittheoremseigenvalues,renfrew2012finiterankdeformationswigner,benaychgeorges2011fluctuationsextremeeigenvaluesfinite,knowles2012isotropicsemicirclelawdeformation,Knowles_2014} for such results on eigenvalues, \cite{capitaine2020nonuniversalityfluctuationsoutlier,bao2020singularvectorsingularsubspace}, and \cite{marcinek2020highdimensionalnormalitynoisy,fan2020asymptotictheoryeigenvectorsrandom} for similar results about eigenvectors but in different settings (where the signal is of macroscopic rank and the noise is Gaussian in the first case, and where $\theta = \theta(n) \to \infty$ in the second case).
In the former works, it is also shown that the fluctuations of eigenvalues are universal when $v$ is sufficiently delocalized.

\begin{figure}
    \centering
    \includegraphics[width=1\linewidth]{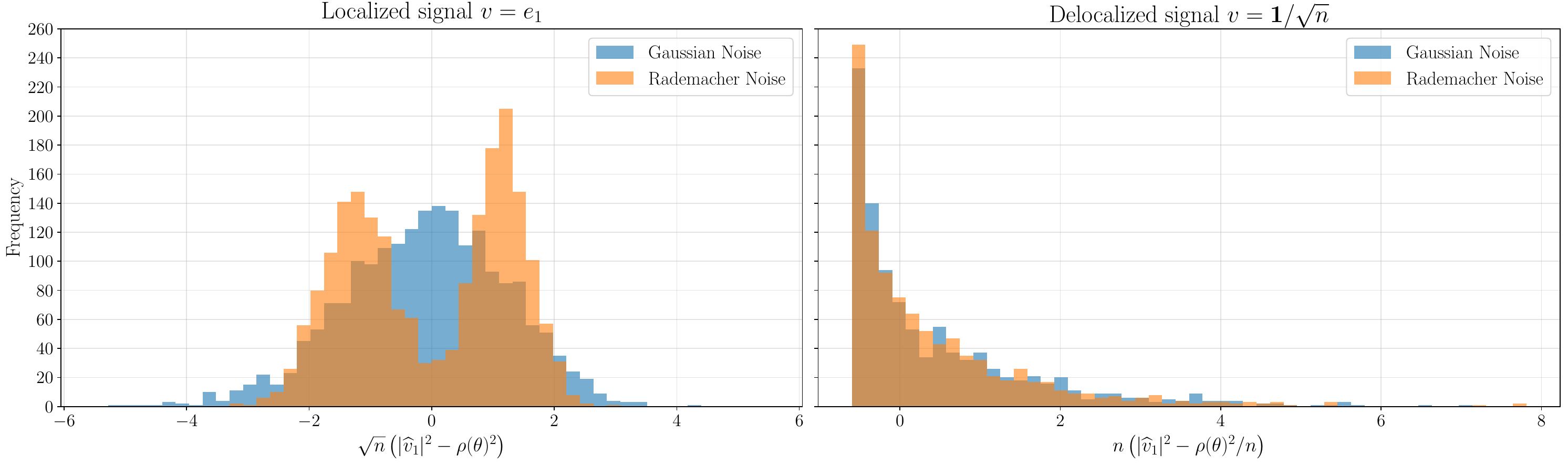}
    \caption{\textbf{Universality and non-universality in spiked matrix models.}
    We illustrate a simple example of the contrast between universality for delocalized signals and non-universality for localized signals.
    Let $P = \what{v}\what{v}^*$ be the projection onto the leading eigenvector of a spiked matrix model of dimension $n = \num{2000}$ with supercritical signal strength $\theta = \num{1.5}$.
    We plot the fluctuations of the diagonal entry $P_{11} = |\what{v}_1|^2$ under its natural scaling in each signal regime.
    We compare noise matrices whose upper-triangular entries are drawn independently from Gaussian and Rademacher distributions.
    In the left panel, the signal $v = e_1$ is localized and $P_{11} = |\langle \what{v}, v\rangle|^2$; the visibly different laws of $\sqrt{n}(P_{11} - \rho(\theta)^2)$ illustrate the non-universality phenomenon proved in \cite{capitaine2020nonuniversalityfluctuationsoutlier}.
    In the right panel, the signal $v = \boldone / \sqrt{n}$ is delocalized, and the agreement of the laws of $n(P_{11} - \rho(\theta)^2/n)$ illustrates the universality of entry statistics proved in Theorem~\ref{thm:main-univ}.}
    \label{fig:non-univ}
\end{figure}

\subsection{Summary of contributions}
\label{sec:summary-contrib}

To the best of our knowledge, until now a gap has remained in the above line of work concerning a situation that is particularly important in applications: it is not known that the \emph{eigenvector} fluctuations---those of $\what{v}$---\emph{are} universal provided that $v$ is \emph{delocalized}, in any stronger sense beyond just the behavior of the specific inner product $\langle \what{v}, v\rangle$.
We illustrate a basic version of this phenomenon in Figure~\ref{fig:non-univ}.

This is the task that we take up in this paper: we will show such a universality result in a rather strong entrywise sense, in the style of the results for eigenvectors of purely random matrices with no low-rank perturbation of \cite{knowles2011eigenvectordistributionwignermatrices}.
We treat both the behavior of the individual entries $n \cdot \widehat{v}_i\overline{\widehat{v}_j}$ and averages of the form
\begin{equation}
\frac{1}{n^2}\sum_{i, j = 1}^n \psi(n \cdot v_i\overline{v_j}, n \cdot \what{v}_i\overline{\what{v}_j}). \label{eq:psi-avg}
\end{equation}
In both cases, we scale by $n$ since the typical scale of each of these numbers is $\Theta(1/n)$.
If we view $\psi$ as a measurement of proximity of two complex numbers, the latter gives us access to various ways of viewing the amount of error made by a spectral algorithm by measuring various notions of entrywise distance between $vv^*$ and $\what{v}\what{v}^*$.

In particular, when we have prior information that $v$ belongs to some special subset of $\CC^n$, we may include \emph{rounding} functions in the definition of $\psi$, that encode a procedure estimating $vv^*$ by first computing $\what{v}\what{v}^*$ and then ``projecting'' its value in some suitable entrywise sense to the set in which the entries $v_i\overline{v_j}$ of $vv^*$ can possibly lie.
As a simple example, if we know in advance that $v \in \{\pm 1 / \sqrt{n}\}^n$ and $W$ is real-valued, then it is sensible to take the entrywise sign $\sgn(\what{v}_i \overline{\what{v}_j}) / n$ as an estimate of $v_i \overline{v_j}$.
We can then measure the natural notion of error of what fraction of these rounded signs disagree with those of $vv^*$, i.e., the quantity:
\[ \frac{1}{n^2}\sum_{i, j = 1}^n \boldone\{\sgn(\what{v}_i \overline{\what{v}_j}) \neq \sgn(v_i \overline{v_j})\}, \]
arguably a more natural and meaningful notion of error for this discrete setting than continuous quantities like $\|\what{v}\what{v}^* - vv^*\|_{\F}^2$.
To the best of our knowledge, ours are the first precise asymptotics for the amount of entrywise error made by such algorithms under general non-Gaussian models of the noise $W$.

We further extend these results to the case of several spiked matrices $H^{(a)} = \theta v^{(a)} v^{(a)^*} + W^{(a)}$ for $a \in [L]$.
In this case, the main structural assumption on the noise matrices $W^{(a)}$ is that the tuples $(W^{(1)}_{ij}, \dots, W^{(L)}_{ij})$ are independent over different pairs $(i, j)$ (up to Hermitian symmetry); however, we note that \emph{within} these tuples the same entry of different $W^{(a)}$ can be dependent or even deterministic functions of one another.
This has an important consequence for the special case where these tuples are uncorrelated, because in this situation, one can compare to the case of $W^{(a)}$ having independent Gaussian entries.
Thus, the case of $H^{(a)}$ being entrywise \emph{uncorrelated} in fact belongs to the same universality class as a case of $H^{(a)}$ being \emph{independent} (as for Gaussian random variables being uncorrelated implies being independent).

We call this phenomenon \emph{pseudoindependence} of spiked matrix models, and we propose that it has interesting ramifications for the design of spectral algorithms.
In particular, if the $v^{(a)}$ are deterministic functions of one another, then clearly in the case of $W^{(a)}$ independent one expects to gain more information about $v^{(a)}$ from more independent matrix observations, and by our analysis the same holds under pseudoindependence.
This suggests that \emph{multi-spectral algorithms}, ones using not just a given spiked matrix $H$ but several deterministic entrywise transformations thereof, may be superior to direct spectral algorithms in some settings.
We illustrate this proposal by showing that multi-spectral algorithms are indeed superior for the problem of synchronization over the cyclic and circle groups.

\subsection{Theoretical results}

We always use the parameter $L$ to designate the number of matrices we work with.
For the sake of brevity we state our results for the case of general $L \geq 1$, but the reader may find it instructive to read the results with the simplifying choice $L = 1$, to recover the case of a single spiked matrix discussed above.

\subsubsection{Universality}

Our first main result says that, provided that the $v^{(a)}$ are sufficiently delocalized and that the $W^{(a)}$ are sufficiently unstructured, statistics of the noisy top eigenvectors $\what{v}^{(a)}$ do not depend on the entry distribution of the $W^{(a)}$.
The specific delocalization assumption we make on the $v^{(a)}$ is the following:
\begin{assumption}
    \label{asm:v}
    For $a\in[L]$, we assume that $v^{(a)} \in \CC^n$ has $\|v^{(a)}\| = 1$, and that
    \[ \max_{a\in[L]}\|v^{(a)}\|_{\infty} \leq C_v n^{-1/2 + \epsilon_v} \]
    for some $\epsilon_v \in (0, 1/20)$ and $C_v > 0$.
    We call $(\epsilon_v, C_v)$ the parameters of this assumption on $v^{(a)}$.
\end{assumption}

\noindent
Our first main result is then as follows.

\begin{thm}[Universality of entry statistics]
    \label{thm:main-univ}
    Let $\theta > 1$ and $v^{(1)}, \dots, v^{(L)} \in \SS^{n - 1}(\CC)$ satisfy Assumption~\ref{asm:v}.
    Let $W^{(1)}, \dots, W^{(L)}, X^{(1)}, \dots, X^{(L)}$ be generalized Wigner matrices (Definition~\ref{def:gen-wigner}) with uniform parameters $(\gamma_W, \xi_W)$ and $(\gamma_X, \xi_X)$ such that the tuples $(W^{(1)}_{ij}, \dots, W^{(L)}_{ij})$ are independent over $1 \leq i \leq j \leq n$, and likewise for the $(X^{(1)}_{ij}, \dots, X^{(L)}_{ij})$.
    Letting $w^{(ij)} \in \RR^{2L}$ contain the real and imaginary parts of $W^{(a)}_{ij}$ for $a \in [L]$ and likewise $x^{(ij)} \in \RR^{2L}$ for $X^{(a)}_{ij}$, suppose that the second moments of these vectors match (recalling that their first moments are fixed to be zero and so automatically match by the assumption of the $W^{(a)}$ and $X^{(a)}$ being generalized Wigner matrices): for all $i, j \in [n]$, we have
    \[ \EE w^{(ij)} w^{(ij)^{\top}} = \EE x^{(ij)} x^{(ij)^{\top}}. \]
    Let $\what{v}^{(W, a)}$ be the eigenvector associated to the largest eigenvalue of $\theta v^{(a)}v^{(a)^{*}} + W^{(a)}$ and $\what{v}^{(X, a)}$ that associated to the largest eigenvalue of $\theta v^{(a)} v^{(a)^{*}} + X^{(a)}$.
    Let $\phi: \CC^L \to \RR$ be a $\sC^5$ function such that $\phi$ and its first five derivatives are bounded.
    Then, for any $\epsilon > 0$,
    \begin{align*}
    &\max_{1 \leq i, j \leq n} \bigg| \EE \phi\left(n \cdot \what{v}^{(W, 1)}_i\overline{\what{v}^{(W, 1)}_j}, \dots, n \cdot \what{v}^{(W, L)}_i\overline{\what{v}^{(W, L)}_j}\right) \\
    &\hspace{2cm} - \EE \phi\left(n \cdot \what{v}^{(X, 1)}_i\overline{\what{v}^{(X, 1)}_j}, \dots, n \cdot \what{v}^{(X, L)}_i\overline{\what{v}^{(X, L)}_j}\right) \bigg| \leq C n^{-1/2 + 10\epsilon_v + \epsilon},
    \end{align*}
    for a constant $C$ depending only on the parameters $L, \epsilon, \theta, \phi$ in the statement of this result, the parameters $(\epsilon_v, C_v)$ of Assumption~\ref{asm:v} on $v$, and the generalized Wigner matrix parameters $(\gamma_W, \xi_W)$ and $(\gamma_X, \xi_X)$.
\end{thm}

\begin{rmk}
    \label{rmk:outer-prod}
    The reason for taking products of pairs of entries of the various $\what{v}$ above is that $\what{v}$ itself is only defined up to a global phase (or a global sign flip in the real-valued case), while $\what{v}\what{v}^{*}$, the orthogonal projection to the $\what{v}$ direction, is a well-defined geometric object; above we formulate our result in terms of its entries.
\end{rmk}

\subsubsection{Gaussian approximation and pseudoindependence}

The above is a conceptually interesting result, but one that does not allow concrete calculations of the actual values of entrywise quantities like the above.
In the special case where the $W_{ij}^{(a)}$ are uncorrelated for a fixed $(i, j)$ over different $a$, however, by comparing to special Gaussian distributions of noise matrices we can obtain tractable approximations.
We consider which $W$ can be compared by Theorem~\ref{thm:main-univ} to the following classical distributions of random matrices:
\begin{definition}[Gaussian orthogonal and unitary ensembles]\label{def:GOE-GUE}
Let $n \geq 1$. We then define the following distributions of random matrices in $\CC^{n \times n}_{\herm}$.
\begin{itemize}
\item The \emph{Gaussian orthogonal ensemble}, denoted $\GOE(n)$ for given dimension $n$, is the law of $W \in \RR^{n \times n}_{\sym}$ with $W_{ij}$ independent for $1 \leq i \leq j \leq n$ drawn as
\begin{align*}
W_{ij} &\sim \sN\left(0, \frac{1}{n}\right) \text{ for } 1 \leq i < j \leq n, \\
W_{ii} &\sim \sN\left(0, \frac{2}{n}\right) \text{ for } 1 \leq i \leq n.
\end{align*}
Equivalently, $\{\sqrt{n} \cdot W_{ij}:1 \leq i<j \leq n\}\cup\{\sqrt{n / 2} \cdot W_{ii}:1\le i\le n\}$ are i.i.d.\ with law $\mathcal N(0,1)$.

\item The \emph{Gaussian unitary ensemble}, denoted $\GUE(n)$ for given dimension $n$, is the law of $W \in \CC^{n \times n}_{\herm}$ with $W_{ij}$ independent for $1 \leq i \leq j \leq n$ drawn as
\begin{align*}
W_{ij} &\sim \sN_{\CC}\left(0, \frac{1}{n}\right) \text{ for } 1 \leq i < j \leq n, \\
W_{ii} &\sim \sN_{\RR}\left(0, \frac{1}{n}\right) \text{ for } 1 \leq i \leq n.
\end{align*}
Here $\sN_{\RR}(0, \sigma^2) = \sN(0, \sigma^2)$ is the ordinary Gaussian measure, while $\sN_{\CC}(0, \sigma^2)$ is the complex Gaussian measure, the law of $Z = X + \cpxi Y$ with $X, Y \sim \sN_{\RR}(0, \sigma^2 / 2)$ independently, so that $\EE|Z|^2 = \sigma^2$.
Equivalently, a GUE matrix $W$ has $\{\sqrt{2n} \cdot \Re W_{ij}: 1 \leq i < j \leq n\} \cup \{\sqrt{2n} \cdot \Im W_{ij}: 1 \leq i < j \leq n\} \cup \{\sqrt{n} \cdot W_{ii}: 1 \leq i \leq n\}$ are i.i.d.\ with law $\sN(0, 1)$.
\end{itemize}
\end{definition}
\noindent
Here and throughout we denote the imaginary unit by $\cpxi$ to avoid conflict with indices named $i$.
It is easy to check that GOE and GUE matrices are both generalized Wigner matrices, up to a negligible renormalization in the GOE case.

When our matrices are drawn from these distributions, the behavior of $\what{v}$ can be analyzed quite precisely.
That is because $X \sim \GOE(n)$ and $X \sim \GUE(n)$ are respectively \emph{orthogonally} and \emph{unitarily invariant}, satisfying the property that $QXQ^{*}$ has the same law as $X$ for any orthogonal or any unitary $Q$, respectively.
Because of this, the component of $\what{v}$ that is orthogonal to $v$ is itself a likewise invariant random vector in this subspace, whose norm by Theorem~\ref{thm:bbp} is close to $\tau(\theta) = \theta^{-1}$.
Thus, setting aside the sign ambiguity in $\what{v}$ mentioned in Remark~\ref{rmk:outer-prod}, approximating this component by a Gaussian random vector (which is also orthogonally invariant) of comparable norm, we expect
\begin{equation}
\what{v}(X) \stackrel{\text{(d)}}{\approx} \sN_{\FF}\left(\rho(\theta) v, \frac{\tau(\theta)^2}{n}(I_n - vv^{*}) + \frac{\sigma(\theta)^2}{n}vv^* \right), \label{eq:gauss-approx}
\end{equation}
for some $\sigma(\theta) = \Theta(1)$ (the correct value of this constant may also be determined explicitly, but in our case of $v$ delocalized this second term will not play a role in the relevant asymptotics).
Here, $\FF \in \{\RR, \CC\}$ is the field we are working over depending on whether $X$ is a real or complex Wigner matrix.
Further, since the mean is a constant multiple of $v$, we should be able to neglect the contribution of $\frac{1}{n}(\tau(\theta)^2 + \sigma(\theta)^2)vv^{*}$ in the covariance.
Similar ideas have appeared before in the literature in \cite{couillet2012fluctuationsspikedrandommatrix,lebeau2024asymptoticgaussianfluctuationseigenvectors}, but we give a full derivation of a precise version of such an approximation specific to our setting in Section~\ref{sec:gauss-models}.

Through Theorem~\ref{thm:main-univ}, together with such analysis, we expect to obtain more concrete predictions for generalized Wigner matrices whose first two moments match those of either a GOE or GUE matrix.
That leads to the following related definition of a more restrictive class of random matrices, in which we also allow some slack in the assumption of exact moment matching that we will see we may deal with in the proof of our next result, and which makes these results more useful for applications.

\begin{definition}[Weakly Wigner tuples]\label{def:weakly-Wigner-tuples}
Let $L\geq 1$. For $1\leq p \leq q \leq n$, define
\[w^{(pq)}\colonequals \left(\Re W^{(1)}_{pq},\Im W^{(1)}_{pq},\ldots,\Re W^{(L)}_{pq},\Im W^{(L)}_{pq}\right)\in\RR^{2L}.\]
A \emph{weakly Wigner tuple} $\bW = (W^{(1)},\ldots W^{(L)})$ with parameters $(\xi_W, \epsilon_W, C_W)$, where each $W^{(a)}\in\FF_a^{n\times n}$ with $\FF_a \in \{\RR, \CC\}$ is Hermitian, is a tuple of matrices satisfying the following properties:
\begin{enumerate}
    \item The vectors $w^{(pq)}$ over $1\leq p \leq q \leq n$ are independent.
    \item $\EE W^{(a)}_{pq} = 0$ and thus $\EE w^{(pq)} = \bm{0}$ for all $p,q\in[n]$, $a\in[L]$.
    \item The $\sqrt{n} \cdot |W^{(a)}_{pq}|$ are $\xi_W$-sub-Weibull for all $p,q\in[n]$, $a\in[L]$.
\end{enumerate}
Define the covariance matrices
\[\Sigma^{(pq)}_{\bW} \colonequals \EE w^{(pq)} w^{(pq) \top} \in \RR^{2L\times 2L}_{\sym} \]
which consist of the blocks
\[\Sigma^{(pq)}_{\bW,ab} \colonequals \begin{pmatrix} \EE\Re W^{(a)}_{pq}\Re W^{(b)}_{pq} & \EE\Re W^{(a)}_{pq}\Im W^{(b)}_{pq} \\ \EE\Im W^{(a)}_{pq}\Re W^{(b)}_{pq} & \EE\Im W^{(a)}_{pq} \Im W^{(b)}_{pq}\end{pmatrix}.\]
We further require the following properties from these covariances:
\begin{enumerate}
    \item[4.] $\|\Sigma_{\bW}^{(pp)}\| \leq C_W / n$ for all $1 \leq p \leq n$.
    \item[5.] Let $\Sigma_{\RR} = \begin{pmatrix}
        1 & 0 \\ 0 & 0
    \end{pmatrix}$, $\Sigma_{\CC} \colonequals \begin{pmatrix}
        1/2 & 0 \\ 0 & 1/2
    \end{pmatrix}$, and $\Sigma_{\GFE} \colonequals \frac{1}{n}\Diag(\Sigma_{\FF_1},\ldots,\Sigma_{\FF_L})$.
    Then, for each $\FF_a \in \{\RR,\CC\}$, $a\in[L]$,
    \[ \|\Sigma^{(pq)}_{\bW} - \Sigma_{\GFE}\| \leq \frac{C_W}{n^{1+\epsilon_W}}, \text{ for all } 1\leq p < q \leq n.\]
\end{enumerate}
\end{definition}
\noindent
We will repeatedly use the variable $\FF_a \in \{\RR, \CC\}$ specifying whether the matrices we are working with are real or complex weakly Wigner matrices.

We note that it is actually not quite true that all matrices in weakly Wigner tuples are generalized Wigner matrices, because the exact condition \eqref{eq:gen-wigner-var-sums} need not hold.
However, the weakly Wigner tuples assumptions imply that it holds up to an $O(n^{-\epsilon_W})$ error, which we will see in our arguments suffices for our purposes.

Our second main result formalizes the above sketch of an argument combined with Theorem~\ref{thm:main-univ}.
We obtain a precise prediction for averages of entrywise functions of $\what{v}\what{v}^{*}$ for $\what{v} = \what{v}(W)$ defined with weakly $\FF$-Wigner matrices $W$ for $\FF \in \{\RR, \CC\}$ matching what we predicted for GOE and GUE matrices, respectively, based on their special invariance properties.
As will be useful in applications, we may also allow these functions to depend on the corresponding entries of $vv^{*}$, allowing us to analyze various notions of the quality of entrywise approximation that $\what{v}\what{v}^*$ gives to $vv^*$, as we suggested above in \eqref{eq:psi-avg}.

\begin{thm}
    \label{thm:main-gauss}
    Let $L\geq 1$ and $\theta > 1$, suppose that $\bW$ is a weakly Wigner tuple (Definition \ref{def:weakly-Wigner-tuples}) with parameters $(\xi_W, \epsilon_W, C_W)$. Let $v^{(a)} \in \SS^{n-1}(\FF_a)$ satisfy Assumption~\ref{asm:v} uniformly for all $a\in[L]$. For each $a\in[L]$, let $\what{v}^{(W,a)}$ be the eigenvector associated to the largest eigenvalue of $\theta v^{(a)}v^{(a)*} + W^{(a)}$.
    Let $\psi: \CC^L \times \CC^L \to \RR$ be a $\sC^5$ function such that the values of $\psi$ and its first five derivatives are bounded.
    Define the associated function $\Psi: (\CC^{n \times n}_{\herm})^L \times (\CC^{n \times n}_{\herm})^L \to \RR$ by
    \[ \Psi(\bA, \bB) \colonequals \frac{1}{n^2} \sum_{i, j = 1}^n \psi\left(\left(A^{(a)}_{ij}\right)_{a=1}^L,\left(B^{(a)}_{ij}\right)_{a=1}^L\right). \]
    Let $g^{(a)} \sim \sN_{\FF_a}(0, I_n)$, $a\in[L]$, be independent standard Gaussian random vectors over the corresponding fields, i.e., having entries $g^{(a)}_i \sim \sN_{\FF_a}(0, 1)$ drawn i.i.d.
    Then, for any $\epsilon > 0$,
    \begin{align*}
    &\Bigg|\EE\Psi\left(\left(nv^{(a)}v^{(a)*}\right)_{a=1}^L,\left(n\what{v}^{(W,a)}\what{v}^{(W,a)*}\right)_{a=1}^L \right) \\ &\hspace{2cm} -\EE\Psi\left(\left(nv^{(a)}v^{(a)*}\right)_{a=1}^L, \left((\rho(\theta)\sqrt{n} v^{(a)}+\tau(\theta)g^{(a)})(\rho(\theta)\sqrt{n}v^{(a)}+\tau(\theta)g^{(a)})^* \right)_{a=1}^L \right) \Bigg| \\ &\hspace{4cm}\leq C\left(n^{-1/2+10\epsilon_v+\epsilon}+n^{-\epsilon_W/2+\epsilon} \right),
    \end{align*}
    where $C$ is a constant depending only on $L,(\FF_a)_{a=1}^L, \epsilon, \theta, \psi$ in the statement of this result, the parameters $(\epsilon_v, C_v)$ of Assumption~\ref{asm:v} on $v^{(a)}$, and the weakly Wigner tuple parameters $(\xi_W,\epsilon_W,C_W)$.

    Suppose further that $v^{(1)},\ldots, v^{(L)}$ are also drawn at random and independently of $\bW$ such that $\sqrt{n} ( v^{(1)}_i,\ldots,v^{(L)}_i)$ for $i\in[n]$ are i.i.d.\ according to some compactly supported probability measure $\mu$ on $\SS^0(\FF_1) \times \dots \times \SS^0(\FF_L)$.\footnote{Concretely, these entries are allowed to lie in $\SS^0(\RR) = \{-1, +1\}$ or the complex unit circle $\SS^0(\CC) = U(1)$.}
    Let $(u_1, \dots, u_L), (v_1, \dots, v_L) \sim \mu$ and $g_a, h_a \sim\sN_{\FF_a}(0,1)$ for $a \in [L]$ all be independent.
    Then, for any $\epsilon > 0$,
    \begin{align*}
    &\Bigg|\EE\Psi\left(\left(nv^{(a)}v^{(a)*}\right)_{a=1}^L,\left(n\what{v}^{(W,a)}\what{v}^{(W,a)*}\right)_{a=1}^L\right) \\ &\hspace{2cm} -\EE\psi\left(\left(u_a\overline{v_a}\right)_{a=1}^L,\left((\rho(\theta)u_a+\tau(\theta)g_a)\overline{(\rho(\theta)v_a+\tau(\theta)h_a)}\right)_{a=1}^L\right)\Bigg| \\ &\hspace{4cm}\leq C\left(n^{-1/2+\epsilon} + n^{-\epsilon_W/2+\epsilon}\right),
    \numberthis
    \label{eq:single-letter}
    \end{align*}
    where $C$ is a constant depending only on $L,(\FF_a)_{a=1}^L, \epsilon, \theta, \psi$ in the statement of this result, the entrywise distribution of the signal $\mu$, and the weakly Wigner tuple parameters $(\xi_W, \epsilon_W, C_W)$.
\end{thm}

\begin{rmk}
    We emphasize the important coordinatewise assumption that, for each $a\in[L]$, $v^{(a)}\in\SS^{n-1}(\FF_a)$, where $\FF_a=\RR$ if $W^{(a)}$ is real symmetric and $\FF_a=\CC$ if $W^{(a)}$ is complex Hermitian; in particular, we exclude the case of $v^{(a)} \in \CC^n\setminus\RR^n$ while $W^{(a)}$ is real symmetric.
\end{rmk}

We think of $\psi$ as an \emph{entrywise loss} function for the estimation of $(v^{(a)}v^{(a)*})_{a=1}^L$ by $(\what{v}^{(a)}\what{v}^{(a)*})_{a=1}^L$, measuring some notion of distance between vectors in $\CC^L$.
$\Psi$ then averages this loss over the entries of a pair of $L$-tuples of matrices.
The final bound \eqref{eq:single-letter} then describes a ``$L$-letter formula'' for the expectation of any such loss: though it is a complicated expectation, it only involves $O(L)$ many random variables while giving precise estimates for large $n$.
In particular, if we have an asymptotic sequence of $v^{(a)} = v^{(n,a)} \in \SS^{n - 1}(\FF_a)$, $a\in[L]$, drawn as above for some fixed joint entrywise distribution $\mu$ on $\FF_1\times\dots\times\FF_L$ not depending on $n$, a sequence of weakly Wigner tuples $\bW = \bW^{(n)}=(W^{(n,1)},\dots,W^{(n,L)})$ satisfying the definition with any parameters $(\xi_W, \epsilon_W, C_W)$ not depending on $n$, and $\what{v}^{(a)} = \what{v}^{(n,a)}$ the top eigenvectors of $\theta v^{(n,a)} v^{(n,a)^*} + W^{(n,a)}$, then the above implies the exact asymptotic result
\begin{equation}
    \label{eq:single-letter-exact}
    \begin{aligned}
    &\lim_{n \to \infty} \EE \Psi\left(\left(n \cdot v^{(a)}v^{(a)*}\right)_{a=1}^L, \left(n \cdot \what{v}^{(a)}\what{v}^{(a)*}\right)_{a=1}^L\right)\\
    &\hspace{2.5cm}= \EE \psi\left(\left(u_a\overline{v_a}\right)_{a=1}^L, \left(\left(\rho(\theta)u_a + \tau(\theta)g_a\right)\overline{\left(\rho(\theta)v_a + \tau(\theta)h_a\right)}\right)_{a=1}^L
    \right).
    \end{aligned}
\end{equation}
Thus, we characterize the asymptotic expectation of any entrywise measurement of the loss achieved by such a spectral algorithm as $n \to \infty$ by an expectation of fixed dimension on the right-hand side, over just the $L$-tuples $(v_a)_{a=1}^L$, $(u_a)_{a=1}^L$, $(g_a)_{a=1}^L$, and $(h_a)_{a=1}^L$, which is complicated to evaluate in closed form for most $\psi$ but can easily be estimated numerically by straightforward Monte Carlo integration methods.

\begin{rmk}[Pseudoindependence]
\label{rmk:pseudoindependence}
We emphasize again that dependence across $a$ is allowed within both tuples: the signals $v^{(1)},\ldots,v^{(L)}$ may be dependent or deterministic functions of one another, and, for each $(i,j)$, the noise entries $W^{(1)}_{ij},\ldots,W^{(L)}_{ij}$ may likewise be dependent or deterministic functions of one another.
In this case, the result may be viewed as describing the \emph{pseudoindependent} asymptotic behavior of such matrices even though in reality they are highly dependent.
\end{rmk}

\subsection{Applications to spectral and multi-spectral algorithms}
\label{sec:intro-appl}

We give applications to the following setting where spectral algorithms have been used with entrywise rounding in the literature \cite{singer2009angularsynchronizationeigenvectorssemidefinite,Cucuringu_2012,Romanov_2020,gao2019multifrequencyphasesynchronization,cucuringu2021extensionangularsynchronizationproblem}.
Let $G$ be a compact Lie group, and write $\Haar(G)$ for its Haar measure.
We suppose that we draw $x_1, \dots, x_n \stackrel{\mathrm{i.i.d.}}{\sim} \Haar(G)$, from which we can form the matrix of pairwise differences:
\[ M_{ij} \colonequals x_i x_j^{-1}, \]
a matrix $M \in G^{n \times n}$ that is \emph{$G$-Hermitian} in the sense that $M_{ji} = M_{ij}^{-1}$.
We observe a noisy version of this matrix,
\begin{equation}
\label{eq:sync-Y}
Y_{ij} \sim \left\{\begin{array}{ll} M_{ij} & \text{with probability } p, \\ \Haar(G) & \text{with probability } 1 - p\end{array}\right\}
\end{equation}
for each $1 \leq i < j \leq n$, and set $Y_{ii} = e$ the group identity\footnote{The diagonal entries will not contribute to the statistics we consider, so this choice is inconsequential.} and $Y_{ji} = Y_{ij}^{-1}$ to preserve the $G$-Hermitian property.
Here $p = p(n) \in (0, 1)$ is some parameter governing the amount of information available in the observation $Y$.
Our goal is to produce an estimator $\what{M} = \what{M}(Y)$ of $M$.
(Note that, for the same reason as in Remark~\ref{rmk:outer-prod}, it is not sensible to try to recover the $x_i$ themselves.)
Such problems in general are referred to as \emph{group synchronization}, and this particular noise model is called the \emph{truth-or-Haar model}, so named by \cite{perry2016optimalitysuboptimalitypcaspiked}.

To assess the quality of the estimator $\what{M}$, we will introduce a \emph{loss function} $\ell: G \times G \to \RR$.
In principle in our results this can be nearly arbitrary, but to interpret our results naturally it should be viewed as a distance-like function of two elements of $G$.
We would then like to assess, given an estimator, the average loss
\[ \EE\left[\frac{1}{n^2}\sum_{i, j = 1}^n \ell(M_{ij}, \what{M}(Y)_{ij})\right]. \]
We expect that the average loss should also concentrate around its expectation under very mild assumptions, and we always observe strong concentration in our experiments, but we do not pursue establishing such results theoretically here.
As a simple concrete example, for $G = \ZZ / K$ we may take simply
\[ \ell(x, y) \colonequals \boldone\{x \neq y \}, \]
in which case the average loss is
\[ \frac{1}{n^2} \sum_{i, j = 1}^n \ell(M_{ij}, \what{M}_{ij}) = \frac{\#\{(i, j) \in [n]^2: M_{ij} \neq \what{M}_{ij}\}}{n^2}, \]
the fraction of incorrectly estimated entries of $M$.

While quite general $G$ have been considered in the literature, we consider two particular examples that are well-suited to applying our theoretical results: the finite cyclic groups $G = \ZZ / K$ for some $K \geq 2$, and the infinite circle group $G = U(1) \cong \mathrm{SO}(2)$, which we can identify with $\RR / 2\pi$.
Note that the cyclic groups are subgroups of $U(1)$; all of these cases are sometimes referred to as \emph{angular synchronization} problems.

The case $K = 2$ has an important additional interpretation: it is equivalent to a dense version of the \emph{stochastic block model}, a much-studied model of community detection in random graphs.
In this case, the $x_i$ are Boolean and may be identified with a label given to each vertex $i$ in a graph, and $M_{ij} \in \{-1, +1\}$ describes whether vertices $i$ and $j$ have the same label or different labels.
In this case, our results describe the asymptotic rate of mislabeling for essentially arbitrary rounded spectral estimators of these pairwise community membership relations.
Some of our numerical experiments will concern this case, but we do not focus on it otherwise since we will see that, because $\ZZ / 2$ has only one non-trivial group character, multi-spectral algorithms (at least in our framework for them) are not applicable to this example.

Let us focus on $G = U(1)$ for the purposes of exposition; our results below apply equally well to the finite cyclic groups as well.
In this case, Singer in \cite{singer2009angularsynchronizationeigenvectorssemidefinite} observed that a reasonable spectral algorithm is as follows.
(See Section~\ref{sec:related} below for further references.)
Given $Y \in G^{n \times n}$, we may build $H \in \CC^{n \times n}_{\herm}$ by taking $H_{ij} = \chi(Y_{ij})$ entrywise for a non-trivial character (or one-dimensional representation) $\chi$ of $G$.
A natural choice when $G = U(1)$ identified with $\RR / 2\pi$ is to simply take
\[ H_{ij} = \frac{1}{\sqrt{n}} \exp(\bm i Y_{ij}). \]
This is close to a spiked matrix model, and thus the top eigenvector $\what{v}$ of $H$ gives a good estimate of the vector $v$ of $v_i = \exp(\bm i g_i)$.
We may then round the $\what{v}_i \overline{\what{v}_j}$ to the unit circle of $\CC$ and take its argument as an estimate $M_{ij} = g_ig_j^{-1}$.

Subsequent literature has explored \emph{multi-frequency} algorithms, which instead of applying a single $\chi$ to $G$ apply several different characters.
Yet it is less clear how to leverage these different frequencies in spectral algorithms; the authors of \cite{Perry_2018}, who studied analogous approximate message passing algorithms, wrote:
\begin{quote}
``In sharp contrast to spectral methods, which offer no reasonable way to couple the frequencies together, AMP produces an estimate that is orders of magnitude more accurate than what is possible with a single frequency.''
\end{quote}

An algorithm somewhat in the spirit of capturing multi-frequency information with spectral algorithms was later proposed in \cite{gao2019multifrequencyphasesynchronization}, but here we explore and will obtain much more precise predictions concerning a direct implementation of this idea.
We propose the following natural class of algorithm.
First, from $Y \in G^{n \times n}$, we build a tuple of Hermitian matrices $H^{(1)},\ldots, H^{(L)}$ by applying several characters of $G$ to the entries of $Y$.
To be concrete, for fixed integers $k_1,\ldots,k_L$ such that the characters below are nontrivial, we take:
\begin{equation}
    \label{eq:sync-chi}
    \chi^{(a)}(x) \colonequals \left\{\begin{array}{ll} \exp(2\pi\cpxi k_a x / K) & \text{if } G = \ZZ / K, \\ \exp(\cpxi k_a x) & \text{if } G = U(1)\end{array}\right\},
\end{equation}
and set, for each $a \in [L]$,
\[ H^{(a)}_{ij} \colonequals \frac{1}{\sqrt{n}}\chi^{(a)}(Y_{ij}). \]
Note that each $H^{(a)}$ is Hermitian since we have chosen $Y$ to be $G$-Hermitian and $\chi^{(a)}(x^{-1})=\overline{\chi^{(a)}(x)}$.
We will see that the covariance condition of Theorem~\ref{thm:main-gauss} is satisfied provided we choose characters that are not only distinct but also not conjugate; this is sensible since applying conjugate characters to $Y$ would yield conjugate matrices $H^{(a)}$ containing essentially the same information.
Thus we assume that, for $a\neq b$, $\chi^{(a)}\notin \{\chi^{(b)},\overline{\chi^{(b)}} \}$. We will see in our proofs that this setup indeed makes $(H^{(1)}, \dots, H^{(L)})$ a weakly Wigner tuple in the sense of Definition~\ref{def:weakly-Wigner-tuples}.

Next, let $\what{v}^{(a)} = v_1(H^{(a)})$ for each $a\in[L]$ be the top eigenvector of each of these matrices.
Finally, let $\Round: \CC^L \to G$ be a function that \emph{rounds} tuples of complex numbers to $G$.
We think of $\Round$ as an inverse of $\bm{\chi} = (\chi^{(1)}, \dots, \chi^{(L)})$, extended to an ``approximate inverse'' on all of $\CC^L$ rather than just $\bm{\chi}(G) \subset \CC^L$.
Using this rounding function, we define the estimator
\[ \what{M}_{ij} \colonequals \Round\left( \left(n \cdot \what{v}^{(a)}_i \overline{\what{v}^{(a)}_j} \right)_{a=1}^L \right). \]
One natural class that we will focus on in our experiments is
\[\Round_{\varrho}(z_1, \dots, z_L) \colonequals \argmin_{x\in G} \left(\frac{1}{L}\sum_{a=1}^L\big|\chi^{(a)}(x) - z_a\big|^{\varrho}\right)^{1/\varrho},\]
the minimum of the generalized mean with parameter $\varrho$ of the distance of $z_a$ to a character value $\chi^{(a)}(x)$ over the $L$ coordinates.
We may take any $\varrho \in [-\infty, \infty]$, the generalized mean becoming a minimum and maximum at the respective extremes.

We note that, comparing again to Singer's spectral algorithm, the case $L=1$ and $k_1=1$ reduces, after identifying $G=U(1)$ with the unit circle, to $\Round_{\varrho}(z)\colonequals z/|z|$ regardless of the choice of $\varrho$, and likewise for $G=\ZZ/K$ to the function that rounds a complex number to the nearest $K$th root of unity.
In both cases, these are the natural rounding functions for producing an estimator $\what{M}$ from a spectral estimate; our $\Round_{\varrho}$ extends these to natural ways to round multi-frequency estimates.

The following uses our earlier results to describe, in considerable generality, the asymptotic average loss of rounded spectral estimators for spectral ($L = 1$) and multi-spectral or multi-frequency ($L \geq 2$) algorithms for group synchronization.

\begin{thm}
    \label{thm:sync}
    Suppose in the above setting that
    $\chi^{(a)} \notin \{\chi^{(b)},\overline{\chi^{(b)}}\}$ for $a\neq b$, that the map $x \mapsto \bm\chi(x) = (\chi^{(1)}(x), \dots, \chi^{(L)}(x))$ is injective on $G$,
    and that
    \[ p = p(n) = \frac{\theta}{\sqrt{n}} \] for some $\theta > 1$.
    Let $\FF_a = \RR$ if $\chi^{(a)}(G) \subseteq \RR$ and $\FF_a = \CC$ otherwise, and view $\FF_1 \times \cdots \times \FF_L$ as a real Euclidean space.
    Let $\Round: \CC^L \to G$ be a Borel function whose restriction to $\FF_1 \times \cdots \times \FF_L$ is continuous Lebesgue-almost everywhere.
    Let $\ell: G \times G \to \RR$ be arbitrary if $G = \ZZ / K$, or a smooth function if $G = U(1)$.
    Let $x, y \sim \Haar(G)$ and $g_a, h_a \sim \sN_{\FF_a}(0, 1)$ for $a \in [L]$ all be drawn independently.
    Then, we have
    \begin{align*}
    &\lim_{n \to \infty} \EE\left[\frac{1}{n^2} \sum_{i, j = 1}^n \ell(M_{ij}, \what{M}_{ij})\right] \\ &\hspace{1cm}= \EE \left[\ell\bigg(xy^{-1}, \Round\left(\left( (\rho(\theta)\chi^{(a)}(x) + \tau(\theta)g_a) \overline{(\rho(\theta)\chi^{(a)}(y) + \tau(\theta)h_a)}\right)_{a=1}^L\right)\bigg)\right]. \numberthis \label{eq:sync-asymp}
    \end{align*}
\end{thm}
\noindent
As in the abstract case earlier, this gives a ``nearly closed form'' for the error rate of such estimators under various notions of error, provided we can estimate the integral of fixed dimension on the right-hand side.
We note that our condition $\theta > 1$ is natural here, since the results of \cite{Kunisky-2025-LCDA2} give evidence that, at least for finite $G$, this problem is \emph{computationally hard} when $\theta < 1$, requiring super-polynomial time to even distinguish an observation $Y \in G^{n \times n}$ from a matrix of i.i.d.\ (up to the $G$-Hermitian relations) Haar-distributed group elements.
Similar results with the same $\theta = 1$ threshold under additive Gaussian noise have also been obtained by \cite{kireeva2024computationallowerboundsmultifrequency, Li-2026-LowerBoundMultiFrequencySynchronization}.

\begin{figure}
    \centering
    \includegraphics[width=1\linewidth]{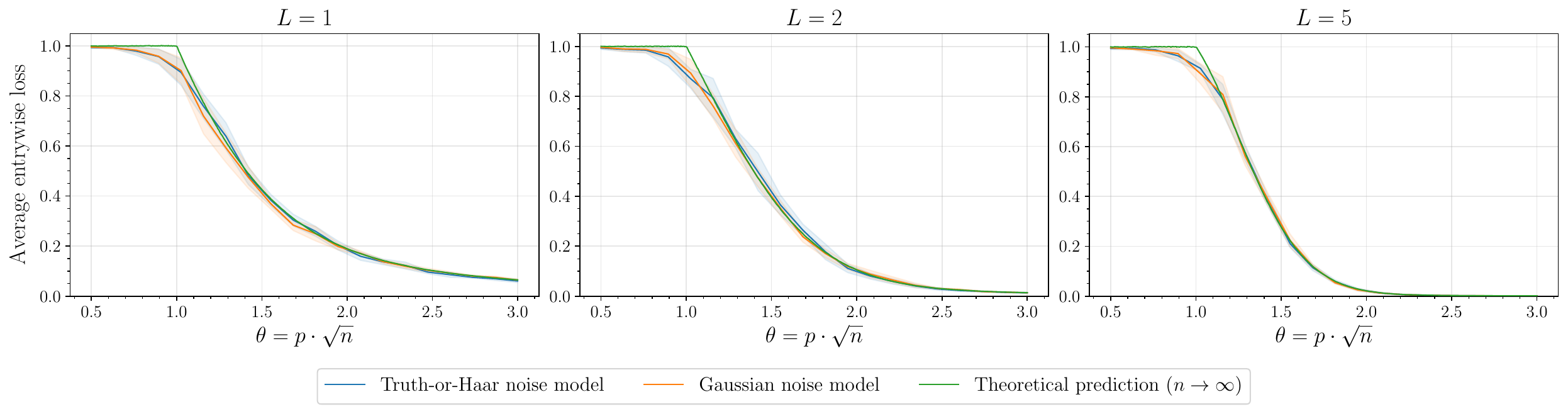}
    \caption{\textbf{Agreement of empirical and predicted performance of synchronization algorithms.} We plot the empirical average loss of the rounded multi-frequency spectral algorithm for group synchronization for $G = U(1)$ (with other parameters as detailed at the end of Section~\ref{sec:intro-appl}) versus the theoretical prediction of Theorem~\ref{thm:sync} for various numbers of frequencies $L$.}
    \label{fig:emp-pred-U}
\end{figure}

\begin{figure}
    \centering
    \includegraphics[width=0.9\linewidth]{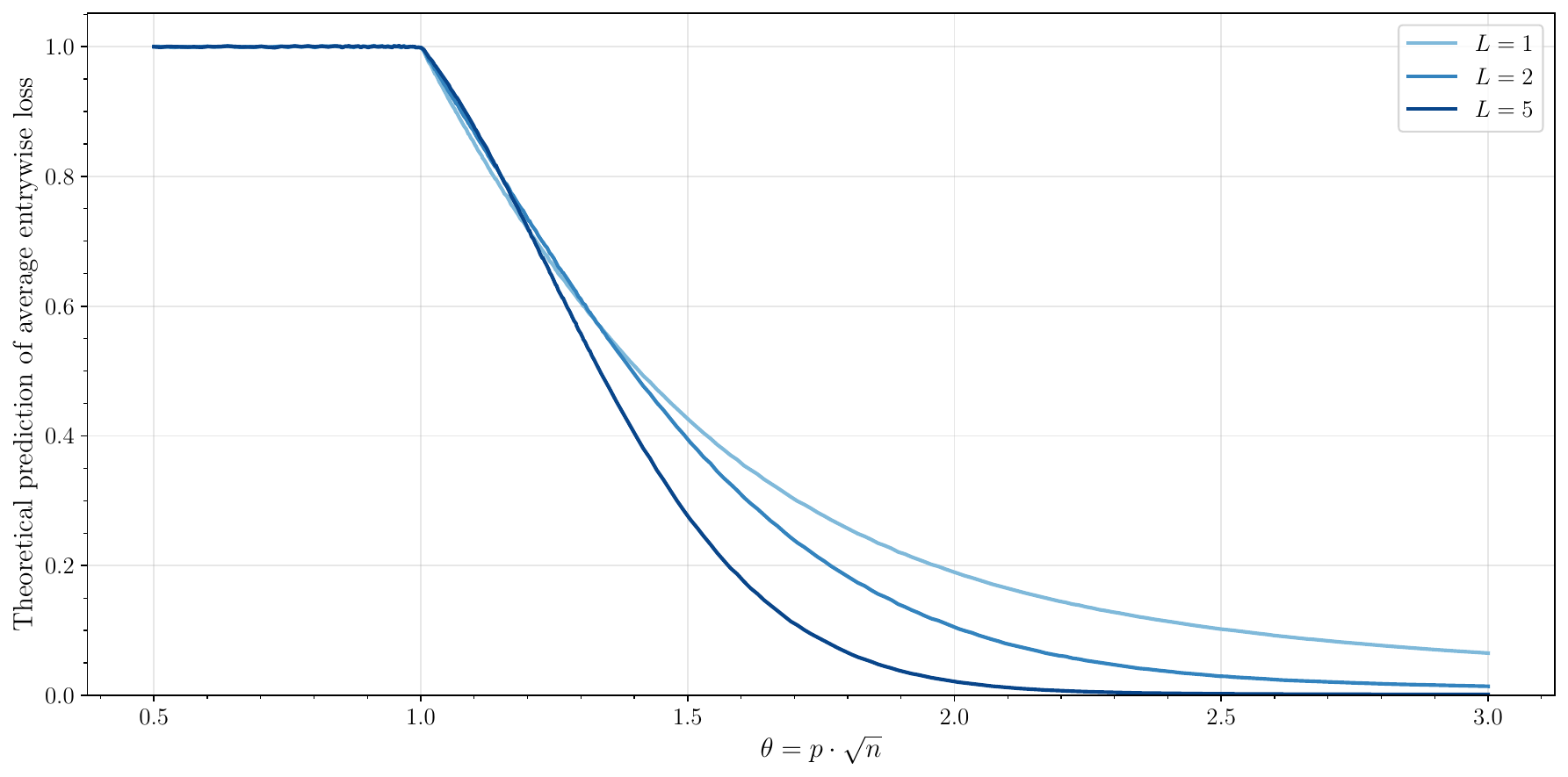}
\caption{\textbf{Predicted superiority of multi-frequency spectral algorithms.} In the same setting as Figure~\ref{fig:emp-pred-U}, for group synchronization over $G = U(1)$, we plot only the theoretical predictions of average entrywise loss per Theorem~\ref{thm:sync} for various numbers of frequencies $L$.}
    \label{fig:multifreq-U}
\end{figure}

We give some illustrative numerical experiments here for the case $G = U(1)$ rounding function $\Round_{\varrho}$ with $\varrho = 2$, and loss function $\ell(x, y) = 1 - \cos(x - y)$; in Section~\ref{sec:appl} we give additional supporting results for the claim of Gaussianity, for the quality of our predictions for different $\varrho$, and for discrete groups $G = \ZZ / K$.
First, in Figure~\ref{fig:emp-pred-U}, for matrices of dimension $n = 500$ and numbers of frequencies $L \in \{1, 2, 5\}$ we plot empirical average losses achieved under the truth-or-Haar noise model as well as the analogous Gaussian additive noise model, and the theoretical prediction of Theorem~\ref{thm:sync}, observing excellent agreement in all cases (especially away from the transition point $\theta = 1$).
Second, in Figure~\ref{fig:multifreq-U}, we plot on a single pair of axes only the theoretical predictions for each number of frequencies (the right-hand side of \eqref{eq:sync-asymp}).
We observe that more frequencies lead to lower estimation error for a large signal strength $\theta$, confirming the superiority of multi-frequency algorithms over the direct spectral algorithm of \cite{singer2009angularsynchronizationeigenvectorssemidefinite}.
Further, as the number of frequencies $L$ increases, the value of $\theta$ where the corresponding algorithm becomes superior to ones for smaller $L$ decreases; thus, in a sense the limit $L \to \infty$ appears to give the ``best'' of these algorithms (though their computational cost also increases with $L$).
It remains unclear to us how to explain why, when comparing two finite $L_1 < L_2$, there is a regime of $\theta$ close to 1 for which the algorithm using only $L_1$ frequencies achieves lower error.
We leave this as an interesting question for future work.

We remark that, in both cases, the theoretical predictions given by evaluating the prediction in Theorem~\ref{thm:sync} above are approximated by Monte Carlo estimation of the expectation, and thus have an inherent error, but we take a sufficient number of trials that this is vanishingly small.

\subsection{Related work}
\label{sec:related}

\paragraph{Spiked matrix models}
The spiked matrix model originates in statistics in the work of \cite{10.1214/aos/1009210544}, who proposed a version where the signal is applied to a covariance matrix of vector observations, though a special case of a model similar to ours was discussed much earlier in the foundational work \cite{FurediKomlos1981Eigenvalues}.
Those predictions were proved in the seminal work of \cite{10.1214/009117905000000233}.
Subsequently, many variants of such results appeared; the first results to the effect of our Theorem~\ref{thm:bbp} for models with additive noise appear to have been shown soon after by \cite{Peche2006LargestEigenvalue,F_ral_2007}.
Other more general versions are shown, for instance, in \cite{Capitaine_2009,benaychgeorges2010eigenvalueseigenvectorsfinitelow}, and an overview of this line of work from the mathematical point of view may be found in \cite{Capitaine-2017-HDRRandomMatrices}, and from the statistical point of view in \cite{PAUL20141,JohnstonePaul2018PCAHighDim}.
The general approach we take to analyzing spiked matrix models using the resolvent appears in \cite{benaychgeorges2010eigenvalueseigenvectorsfinitelow,benaychgeorges2011fluctuationsextremeeigenvaluesfinite} and was used more systematically together with \emph{isotropic local laws} by \cite{knowles2012isotropicsemicirclelawdeformation,Knowles_2014,knowles2016anisotropiclocallawsrandom}, whose techniques we will follow closely.
Another relevant work for studying eigenvectors by this method is the earlier one of \cite{knowles2011eigenvectordistributionwignermatrices}, though this concerns matrices like our $W$ with no low-rank perturbation.

\paragraph{Universality and non-universality}
As we discussed briefly above, the general phenomenon of non-universality of fluctuations in spiked matrix models depending on the structure of the signal was gradually uncovered by works including \cite{Capitaine_2009,capitaine2011centrallimittheoremseigenvalues,renfrew2012finiterankdeformationswigner,benaychgeorges2011fluctuationsextremeeigenvaluesfinite,knowles2012isotropicsemicirclelawdeformation,Knowles_2014}.
For eigenvectors, non-universality for localized signals was shown by \cite{capitaine2020nonuniversalityfluctuationsoutlier}, with regard to a quantity similar to $\langle \what{v}, v\rangle$.
We note that this is a special projection of $\what{v}$ to study because, for $\theta > 1$, per Theorem~\ref{thm:bbp} it is of magnitude $\Theta(1)$; in contrast, for $v$ delocalized the entries $\what{v}_i$ are of magnitude $o(1)$ and therefore considerably more delicate to understand.
Non-asymptotic entrywise eigenvector bounds, including an application to $\ZZ / 2$ synchronization, were obtained by \cite{abbe2019entrywiseeigenvectoranalysisrandom}.
Distributions of similar quantities in asymmetric and rectangular matrix denoising models and in supercritical spiked covariance models are considered by \cite{bao2020singularvectorsingularsubspace,bao2020statisticalinferenceprincipalcomponents}.
The results perhaps most similar to ours are those of \cite{fan2020asymptotictheoryeigenvectorsrandom,bhattacharya2026outliereigenvalueseigenvectors}, but both consider, in our normalization, the case $\theta = \theta(n) \to \infty$, in which case $\what{v}$ is very close to $v$ and its Gaussian fluctuations have vanishing magnitude.

\paragraph{Approximate message passing and state evolution}
One important other branch of the literature where similar statements of Gaussian fluctuations appear is in the characterization of \emph{approximate message passing (AMP)} algorithms.
See, e.g., \cite{feng2021unifyingtutorialapproximatemessage} for a general survey of this area.
These algorithms perform a general kind of \emph{nonlinear power method}, alternating multiplying a vector by a matrix like our $H$ and applying an entrywise nonlinearity.
Various expressions of the Gaussian distribution of the result of such algorithms are known as \emph{state evolution} descriptions.
In principle, the power method itself for computing $\what{v}$ is a special case of AMP where the nonlinearities are omitted, and thus various averages of nonlinear functions of $\what{v}$ can be described as measurements of the state of a suitable AMP algorithm.
However, there is an important caveat that much of the work on AMP concerns only asymptotic results as $n \to \infty$ for a finite number of iterations $t = O(1)$ of such an algorithm, from a random initialization.
Such an algorithm essentially would compute $H^tx$ for a random $x$ independent of $H$, which does not faithfully approximate $\what{v}$, and thus such analysis cannot be used in our setting directly.

In the AMP literature there have been essentially two workarounds considered to this kind of issue (variants of which are also relevant to other, more sophisticated instances of AMP).
First, one may perform non-asymptotic analysis of a number of iterations $t = t(n)$ depending on $n$, as pioneered by \cite{Rush_2018} and developed further by \cite{li2023nonasymptoticframeworkapproximatemessage,li2023approximatemessagepassingrandom,cademartori2023nonasymptoticanalysisgeneralizedapproximate}.
However, such analysis requires strong assumptions on $W$ such as Gaussianity or orthogonal invariance.
The recent work \cite{han2025entrywisedynamicsuniversalitygeneral} is in a similar spirit to ours, but only can treat $t(n) \leq (\log n)^{1/3}$ iterations of an AMP algorithm, while $t(n) = \Theta(\log n)$ are required to faithfully estimate $\what{v}$ in spiked matrix models.

Second, a line of work of \cite{montanari2019estimationlowrankmatricesapproximate,mondelli2021approximatemessagepassingspectral,mondelli2021optimalcombinationlinearspectral}, also discussed in Section~3.2 of \cite{feng2021unifyingtutorialapproximatemessage}, has sought to analyze AMP initialized not from a random vector $x$ independent of $H$, but from precisely our $\what{v}$, a so-called \emph{spectral initialization} from the top eigenvector.
Of course, if this were possible, then we could study $\what{v}$ very directly using such analysis.
However, these works involve various tricks that rely deeply on the Gaussian structure of $W$, and in their current form do not seem to imply any analysis of $\what{v}$ comparable to ours for non-Gaussian $W$, or indeed for any $W$ besides very symmetric Gaussian models like the GOE and GUE---the information about $\what{v}$ that these approaches use stems from precisely the same invariance properties of the GOE and GUE that we take advantage of in our proof of Theorem~\ref{thm:main-gauss}.

\paragraph{Group synchronization}
The group synchronization problem asks to recover group-valued labels from noisy pairwise relative measurements.
Algorithmic approaches include Singer's eigenvector and semidefinite programming methods for angular synchronization \cite{singer2009angularsynchronizationeigenvectorssemidefinite}, message passing over compact groups \cite{Perry_2018}, and multi-frequency phase synchronization \cite{gao2019multifrequencyphasesynchronization}.
Yang, Wee, and Fan \cite{YWF-2024-MutualInformationQuadraticCompactGroups} characterize the
  information-theoretic limits of inference in several quadratic models over compact groups, including
  angular and phase synchronization.
Other work proves lower bounds for low-degree polynomials and low-coordinate-degree algorithms in
  Gaussian multi-frequency and truth-or-Haar synchronization models
  \cite{kireeva2024computationallowerboundsmultifrequency,Kunisky-2025-LCDA2,Li-2026-LowerBoundMultiFrequencySynchronization}, giving evidence for computational hardness below the threshold $\theta=1$.
Our results complement these by giving precise formulas for the asymptotic performance of specific entrywise rounded spectral estimators when $\theta > 1$.

\subsection{Notation}

We use the standard notation $[n] = \{1, \dots, n\}$.
We occasionally use the special function
\[ \varphi(n) \colonequals (\log n)^{\log\log n}. \]
We write $\cpxi = \sqrt{-1}$ for the imaginary unit, not to be confused with the indices $i \in [n]$.
For $\FF \in \{\RR, \CC\}$, we denote the associated unit sphere in dimension $n$ by $\SS^{n - 1}(\FF) = \{x \in \FF^n: \|x\| = 1\}$.
For a proposition $P$ about various variables defined in our arguments, we write $\boldone\{P\}$ to equal 1 if $P$ is true and 0 if $P$ is false.
Reusing this notation, in probabilistic arguments we also use $\boldone_{\sE}$ for the indicator random variable of an event $\sE$.

We write $\CC^{n \times n}_{\herm}$ for the set of Hermitian matrices and $\RR^{n \times n}_{\sym}$ for the set of real symmetric matrices.
$I = I_n$ denotes the $n \times n$ identity matrix; we omit the subscript if it is clear from context.
For $X \in \CC^{n \times n}_{\herm}$, we write $\lambda_1(X) \geq \lambda_2(X) \geq \cdots \geq \lambda_n(X)$ for its ordered real eigenvalues and $v_1(X), \dots, v_n(X)$ for the associated eigenvectors, provided the corresponding eigenvalues are simple.

The asymptotic notations $\lesssim, \gtrsim, O(\cdot), \Theta(\cdot), o(\cdot), \omega(\cdot), \Omega(\cdot)$ have their usual meanings, always with respect to the limit $n \to \infty$.
We sometimes write subscripts on these notations for parameters that the implicit constants depend on, but when proving a result whose statement gives this dependence explicitly, we omit those parameters for the sake of brevity.
In general, all of the quantities we have called ``parameters'' in the statements, with subscripts of ``$v$'' or ``$W$'', are viewed as constants for the purposes of our proofs.

\section{Preliminaries}

\subsection{Probability}

Let us give a name to the tail bound assumption \eqref{eq:entries-subexp} that we make on the magnitudes of the entries of our generalized Wigner matrices.
\begin{definition}
    \label{def:power-subexp}
    We say that a random variable $X$ is \emph{$\xi$-sub-Weibull} if, for all $t > 0$,
    \[ \PP[|X| \geq t] \leq \xi^{-1} \exp(-t^{\xi}). \]
\end{definition}

\begin{prop}
    \label{prop:subexp-moments}
    If $X$ is $\xi$-sub-Weibull and $p \geq 1$, then there exists $C(\xi, p)$ such that $\EE|X|^p \leq C(\xi, p)$.
\end{prop}
\begin{proof}
    We have
    \begin{align*}
        \EE|X|^p
        &= \int_{t = 0}^{\infty} \PP[|X|^p \geq t]\,dt \\
        &= \int_{t = 0}^{\infty} \PP[|X| \geq t^{1/p}]\,dt \\
        &\leq \xi^{-1} \int_{t = 0}^{\infty} \exp(-t^{\xi/p})\,dt,
    \end{align*}
    where the remaining integral is finite and depends only on $\xi$ and $p$.
\end{proof}

The following notation for tail bounds up to a polynomial amount of ``slack'' will be useful throughout.
See Appendix~A of \cite{AEK-2017-UniversalityWignerTypeMatrices} for further generalities about this definition.
\begin{definition}[Polynomial stochastic domination]\label{def:stochastic-domination}
For sequences of non-negative random variables $X = X_n$ and $Y = Y_n$, we say $X$ is \emph{polynomially stochastically dominated} by $Y$, written as $X \prec Y$ or $(X_n)\prec (Y_n)$, if for every $\varepsilon>0$ and for all $D>0$, there exists $n_0 = n_0(\varepsilon,D)$ such that, for all $n\ge n_0$,
\[
\PP[X_n > n^{\varepsilon}Y_n] \le n^{-D}.
\]
Also, if $\sE_n$ is a sequence of events, we say that $\sE_n$ occurs with \emph{polynomially high probability} if, for every $D > 0$, there exists $n_0 = n_0(D)$ such that, for all $n \geq n_0$,
\[ \PP[\sE_n^c] \leq n^{-D}. \]
\end{definition}
\noindent
The following property of polynomial stochastic domination is elementary to verify by the union bound.
\begin{prop}
    \label{prop:stoch-dom-sum-prod}
    Suppose that $X_{n, 1}, \dots, X_{n, k}, Y_{n, 1}, \dots, Y_{n, k}$ are non-negative random variables such that $(X_{n, i}) \prec (Y_{n, i})$ for each fixed $i \in [k]$.
    Then,
    \begin{align*}
        \sum_{i = 1}^k X_{n, i} &\prec \sum_{i = 1}^k Y_{n, i}, \\
        \prod_{i = 1}^k X_{n, i} &\prec \prod_{i = 1}^k Y_{n, i}.
    \end{align*}
\end{prop}
\noindent
The following gives a convenient condition under which polynomial stochastic domination results can be converted to control of expectations.
\begin{prop}\label{prop:stochastic-domination-to-moment-bound}
Let $(X_n)_{n\ge1}$ be a sequence of non-negative random variables and $(Y_n)_{n\ge1}$ be a sequence of deterministic non-negative numbers. Suppose that the following hold:
\begin{enumerate}
\item $X_n\prec Y_n$.
\item $\sup_{n \geq 1} \EE X_n^2 < \infty$.
\item There exists $C>0$ such that $Y_n \ge n^{-C}$ for all $n$.
\end{enumerate}
Then, for any $\epsilon > 0$,
\[ \EE X_n \lesssim_{\epsilon} n^{\epsilon} Y_n. \]
\end{prop}
\begin{proof}
Let $K = \sup_{n \geq 1} \EE X_n^2$.
We have
\begin{align*}
    \EE X_n
    &= \EE X_n \boldone_{X_n \leq n^{\epsilon}Y_n} + \EE X_n \boldone_{X_n > n^{\epsilon} Y_n}
    \intertext{and using the bound from the indicator on the first term and the Cauchy-Schwarz inequality on the second,}
    &\leq n^{\epsilon} Y_n + (\EE X_n^2)^{1/2}(\PP[X_n > n^{\epsilon} Y_n])^{1/2}
    \intertext{Now, choosing some $D \geq 2C$, if $n \geq n_0(\epsilon, D)$, then we have}
    &\leq n^{\epsilon}Y_n + K^{1/2} n^{-D/2} \\
    &\leq n^{\epsilon}Y_n + K^{1/2} Y_n \\
    &\leq (1 + K^{1/2})n^{\epsilon}Y_n,
\end{align*}
completing the proof.
\end{proof}

\subsection{Matrix resolvents and low-rank perturbations}

We will work extensively with the resolvents of random matrices, whose basic properties we enumerate below.

\begin{definition}
    The \emph{resolvent} of a matrix $X$ at the point $z \in \CC$ is
    \[ \Res_z(X) = (zI - X)^{-1}, \]
    defined at any $z$ that is not an eigenvalue of $X$.
    In particular, if $X \in \CC^{n \times n}_{\herm}$, then $\Res_z(X)$ is defined away from a compact subset of $\RR$.
\end{definition}

\begin{prop}[Riesz projection formula]
\label{prop:proj-int}
Let $\Gamma$ be a counterclockwise simple contour in $\CC$ that does not pass through any eigenvalues of a matrix $X \in \CC^{n \times n}_{\herm}$ and encircles a single eigenvalue $\lambda$ of $X$.
Let $P$ be the orthogonal projection to the eigenspace of $X$ associated to $\lambda$.
Then,
\[ P = \frac{1}{2\pi i}\oint_{\Gamma} \Res_z(X)\,dz. \]
\end{prop}

\begin{prop}[Resolvent identities]
    \label{prop:res-id}
    Let $X, Y \in \CC^{n \times n}_{\herm}$ and $w, z \in \CC \setminus (\spec(X) \cup \spec(Y))$.
    Then, the following hold:
    \begin{enumerate}
        \item (First resolvent identity) $R_X(z) - R_X(w) = (w - z) R_X(z) R_X(w)$.
        \item (Second resolvent identity) $R_X(z) - R_Y(z) = R_X(z) (X - Y) R_Y(z)$.
        \item (Resolvent expansion) $R_X(z) = \sum_{s = 0}^t R_Y(z) ((X - Y) R_Y(z))^s + R_X(z)((X - Y)R_Y(z))^{t + 1}$ for any $t \geq 0$.
    \end{enumerate}
\end{prop}
\noindent
The two resolvent identities are standard, while the resolvent expansion follows from repeatedly applying the second resolvent identity.

One important use of resolvents is in characterizing the eigenvalues and eigenvectors of a low-rank perturbation of a matrix.
We will use these results on random spiked matrix models, but they are true deterministically as well.
We state the results for the general case
\[ H = \sum_{i = 1}^k \theta_i v_iv_i^* + W = V\Diag(\theta)V^* + W. \]

To study the eigenvalues, we use the following device, used in derivations of the main phase transition phenomena of spiked matrix models by, e.g., \cite{knowles2012isotropicsemicirclelawdeformation,Knowles_2014}; see also Section 5 of \cite{benaychgeorges2010eigenvalueseigenvectorsfinitelow}.
This characterization of the eigenvalues of a low-rank perturbation is sometimes called the \emph{secular equation}.
\begin{prop}
    \label{prop:res-eig}
    Suppose that $\theta_i \neq 0$ for every $i \in [k]$ and that $\lambda \notin \spec(W)$.
    Then, $\lambda \in \spec(H)$ if and only if
    \[ \det(\Diag(\theta)^{-1} - V^* R_W(\lambda) V) = 0. \]
    If $k = 1$, $\theta = \theta_1$, and $v = v_1$, then this reduces to
    \[ v^* R_W(\lambda) v = \frac{1}{\theta}. \]
\end{prop}
\begin{proof}
    Since $\lambda \notin \spec(W)$, the matrix $R_W(\lambda)^{-1} = \lambda I - W$ is well-defined and non-singular.
    We therefore have
    \begin{align*}
        \det(\lambda I - H)
        &= \det(\lambda I - W - V\Diag(\theta)V^*) \\
        &= \det(R_W(\lambda)^{-1})
        \det(I - R_W(\lambda)V\Diag(\theta)V^*) \\
        &= \det(R_W(\lambda)^{-1})
        \det(I - \Diag(\theta)V^*R_W(\lambda)V) \\
        &= \det(R_W(\lambda)^{-1})\det(\Diag(\theta))
        \det(\Diag(\theta)^{-1} - V^*R_W(\lambda)V),
    \end{align*}
    where the third line follows from Sylvester's determinant identity.
    The first two determinant factors in the last line are non-zero, which gives the result.
\end{proof}

We also use the following calculation of the resolvent of a spiked matrix model as a low-rank perturbation of the resolvent of the underlying noise matrix.
\begin{prop}
    \label{prop:res-pert}
    Suppose that $\theta_i \neq 0$ for every $i \in [k]$ and that $z \notin \spec(W) \cup \spec(H)$.
    Then,
    \[ R_H(z) = R_W(z) + R_W(z) V(\Diag(\theta)^{-1} - V^*R_W(z)V)^{-1}V^* R_W(z). \]
\end{prop}
\begin{proof}
    We calculate directly:
    \begin{align*}
        R_{W + V\Diag(\theta)V^*}(z)
        &= (zI - W - V\Diag(\theta)V^*)^{-1} \\
        &= (I - R_W(z) V\Diag(\theta)V^*)^{-1} R_W(z) \\
        &= \left(I + R_W(z)V(\Diag(\theta)^{-1} - V^*R_W(z)V)^{-1}V^*\right)R_W(z) \\
        &= R_W(z) + R_W(z) V(\Diag(\theta)^{-1} - V^*R_W(z)V)^{-1}V^* R_W(z),
    \end{align*}
    where we have used the Woodbury matrix inverse identity in the middle.
\end{proof}
\noindent
Using this, we may compute the eigenvectors of $H$.
\begin{prop}
    \label{prop:res-eigvec}
    Suppose that $\lambda \in \spec(H) \setminus \spec(W)$ is an eigenvalue with some multiplicity $\ell$, and let $P$ be the orthogonal projection to the associated ($\ell$-dimensional) eigenspace.
    Then, the dimension of $\ker(\Diag(\theta)^{-1} - V^*R_W(\lambda)V)$ is $\ell$.
    Letting $U \in \CC^{k \times \ell}$ have as its columns an orthonormal basis for this kernel, we have
    \[ P =  R_W(\lambda) V U(U^* V^* R_W(\lambda)^2 V U)^{-1}U^* V^* R_W(\lambda). \]
    Equivalently, the eigenspace of $H$ associated to $\lambda$ is the column space of $R_W(\lambda)V U$, spanned by $R_W(\lambda)Vu$ over all $u$ such that $(\Diag(\theta)^{-1} - V^*R_W(\lambda)V)u = 0$.

    If $k = 1$, $\theta = \theta_1$, $v = v_1$, and $\lambda$ is a simple eigenvalue with eigenvector $\what{v}$ (matching our earlier notation in the case of $\lambda$ the largest eigenvalue of $H$), then this reduces to
    \[ \what{v}\what{v}^* = \frac{1}{v^* R_W(\lambda)^2 v} R_W(\lambda) vv^* R_W(\lambda), \]
    and $\what{v}$ up to rescaling is $R_W(\lambda)v$.
\end{prop}
\begin{proof}
    Let $\Gamma$ be a closed contour encircling $\lambda$ but no other eigenvalues of $W$ or $H$.
    We have, by Propositions~\ref{prop:proj-int} and~\ref{prop:res-pert},
    \begin{align*}
    P
    &= \frac{1}{2\pi i} \oint_{\Gamma} R_H(z)\,dz \\
    &= \frac{1}{2\pi i} \oint_{\Gamma} \big(R_W(z) + R_W(z)V(\Diag(\theta)^{-1} - V^*R_W(z) V)^{-1}V^*R_W(z)\big)\,dz
    \intertext{and, since $\Gamma$ avoids the eigenvalues of $W$, $R_W(z)$ is analytic on an open neighborhood of the interior of $\Gamma$, and thus by Cauchy's integral theorem the first term integrates to zero and we have}
    &= \frac{1}{2\pi i} \oint_{\Gamma} R_W(z)V(\Diag(\theta)^{-1} - V^*R_W(z) V)^{-1}V^*R_W(z)\,dz.
    \intertext{Now, by Proposition~\ref{prop:res-eig}, the integrand has poles precisely $z$ equaling the eigenvalues of $H$.
    Since the only one of these that $\Gamma$ encircles is $\lambda$, the above equals the (matrix-valued) residue at this pole.
    Again by our assumption, $R_W(z)$ is analytic in an open neighborhood of the interior of $\Gamma$, and so we may take the (matrix-valued) residue of the inner term $(\Diag(\theta)^{-1} - V^*R_W(z)V)^{-1}$ at $z = \lambda$.
    Let $U \in \CC^{k \times \ell}$ have as its columns an orthonormal basis of $\ker(\Diag(\theta)^{-1} - V^*R_W(\lambda)V)$.
    Since $\frac{d}{dz}R_W(z) = -R_W(z)^2$, this is given by $U(U^* V^* R_W(\lambda)^2 V U)^{-1}U^*$.
    Thus we may evaluate and we find}
    &= R_W(\lambda) V U(U^* V^* R_W(\lambda)^2 V U)^{-1}U^* V^* R_W(\lambda).
    \end{align*}
    Note that, since $\lambda \in \RR$, $R_W(\lambda)$ is Hermitian, and so this is just the formula for the orthogonal projection to the span of the vectors $R_W(\lambda) VU$.
\end{proof}

\subsection{Random matrix theory}

\subsubsection{Semicircle limit theorems}

Let us review some of the properties of generalized Wigner random matrices we have mentioned above.
\begin{definition}
    The \emph{semicircle distribution}, denoted $\mu_{\sc}$, is the probability measure with density
    \[ \frac{1}{2\pi} \sqrt{4 - x^2} \,\boldone\{x \in [-2, 2]\} \,dx \]
    with respect to Lebesgue measure.
\end{definition}
\noindent
The following general limit theorem is a direct consequence of the much more precise results of \cite{erdos2011rigidityeigenvaluesgeneralizedwigner} on generalized Wigner matrices.
\begin{thm}
    \label{thm:semicircle}
    Let $W^{(n)} \in \CC^{n \times n}_{\herm}$ be a sequence of generalized Wigner matrices (Definition~\ref{def:gen-wigner}) with parameters $(\gamma_W, \xi_W)$ not depending on $n$.
    Then, the following convergences hold in probability:
    \begin{align*}
        \frac{1}{n}\sum_{i = 1}^n f(\lambda_i(W^{(n)})) &\to \int fd\mu_{\sc}, \\
        \lambda_1(W^{(n)}) &\to 2, \\
        \lambda_n(W^{(n)}) &\to -2, \\
        \|W^{(n)}\| &\to 2.
    \end{align*}
    In the first claim, $f: \RR \to \RR$ is any bounded continuous function or any polynomial.
\end{thm}

Many of the calculations to follow will also involve the following object evaluated with the semicircle measure, a few of whose properties we establish now.
\begin{definition}
    For a probability measure $\mu$ on $\RR$ supported on some $K \subseteq \RR$, its \emph{Cauchy transform} is the function $G_{\mu}: \CC \setminus K \to \CC$ defined by
    \[ G_{\mu}(z) = \int \frac{1}{z - x} d\mu(x). \]
\end{definition}

\begin{prop}
    \label{prop:Gsc}
    The Cauchy transform $G_{\mu_{\sc}}(z)$ of the semicircle distribution is defined on $\CC \setminus [-2, 2]$ and satisfies:
    \begin{align*}
        G_{\mu_{\sc}}(z) &= \frac{1}{2}(z - \sqrt{z^2 - 4}), \\
        G_{\mu_{\sc}}^{\prime}(z) &= \frac{1}{2}\left(1 - \frac{z}{\sqrt{z^2 - 4}}\right).
    \end{align*}
    In particular, on $(2, \infty) \subset \RR$, $G_{\mu_{\sc}}^{\prime}$ increases monotonically from $-\infty$ to 0.
\end{prop}

\subsubsection{Tail bounds}

We follow works like \cite{knowles2012isotropicsemicirclelawdeformation} in adopting the notation
\[ \varphi(n) \colonequals (\log n)^{\log\log n} \]
for a nuisance factor that appears in the following results.
For intuition, it is useful to keep in mind that $\varphi(n)$ grows slower than any polynomial of $n$ but faster than any polynomial of $\log n$: for any $C > 0$,
\[ (\log n)^C \ll \varphi(n) \ll n^{1/C}. \]
The following concrete norm bounds will be useful.
\begin{prop}[Wigner matrix norm deviation bounds \cite{erdos2011rigidityeigenvaluesgeneralizedwigner}]
    \label{prop:wig-norm}
    Let $W \in \CC^{n \times n}_{\herm}$ be a generalized Wigner matrix.
    Then,
    \[ \PP[\|W\| \geq 2 + \varphi(n) n^{-2/3}] \leq C\exp(-\varphi(n)) \]
    for a constant $C$ depending only on the parameters $(\gamma_W, \xi_W)$ of the generalized Wigner matrix.
\end{prop}

\begin{prop}[Resolvent norm bounds]
    \label{prop:res-moments}
    Let $W \in \CC^{n \times n}_{\herm}$ be a generalized Wigner matrix, and let $z \in \CC$ have $\Im(z) > n^{-C}$ and $\Re(z) > 2 + c$ for some $c, C > 0$.
    Then,
    \[ \EE\|R_W(z)\|^k \leq K \]
    for some constant $K$ depending only on $c, C, k$, and the parameters $(\gamma_W, \xi_W)$ of the generalized Wigner matrix.
    The same also holds for $W$ replaced by the matrix $Q$ formed by setting entries $(a, b)$ and $(b, a)$ to zero in $W$, for any $1 \leq a \leq b \leq n$.
\end{prop}
\begin{proof}
    We have
    \begin{align*}
        \|R_W(z)\|
        &= \|(zI - W)^{-1}\| \\
        &\leq \frac{1}{d(z, \spec(W))} \\
        &\leq \frac{1}{\Im(z)}.
    \end{align*}
    Fix some $\delta < c$.
    Define the event
    \[ \sE = \{\|W\| \leq 2 + \delta\} = \{\spec(W) \subseteq [-2 -\delta, 2 + \delta]\}. \]
    Then, we have using Proposition~\ref{prop:wig-norm} that
    \begin{align*}
        \EE \|R_W(z)\|^k
        &= \EE \boldone_{\sE} \|R_W(z)\|^k + \EE \boldone_{\sE^c} \|R_W(z)\|^k \\
        &\leq \left(\frac{1}{c - \delta}\right)^k + \Im(z)^{-k} \PP[\sE^c] \\
        &\leq O_{c, k}(1) + n^{Ck} \exp(-\varphi(n)),
    \end{align*}
    and the result for $W$ follows since $\varphi(n)$ grows faster than any power of $\log n$.
    The result for $Q$ follows by the same argument since $\|W - Q\| = |W_{ab}|$, which satisfies a sub-Weibull tail bound.
\end{proof}

\subsubsection{Local laws}

We will use the following powerful result about deterministic quadratic forms with the resolvent of a generalized Wigner random matrix.
\begin{thm}[Isotropic local law outside the spectrum, Theorem 2.15 and Remark 2.6 of \cite{bloemendal2015isotropiclocallawssample}]
    \label{thm:ill}
    Fix $c, C > 0$.
    Let $W \in \CC^{n \times n}_{\herm}$ be a generalized Wigner matrix, let $x,y \in \SS^{n-1}(\CC)$ be deterministic, and define
    \[ \Omega = \Omega(c, C) = \{z \in \CC: \Re(z) > 2 + c, \Im(z) \geq 0, |z| \leq C\}. \]
    Then,
    \begin{equation}
        \label{eq:ill-dom}
    \sup_{z \in \Omega} \left|x^*R_W(z)y - G_{\mu_{\sc}}(z)\cdot x^*y\right| \prec n^{-1/2},
    \end{equation}
    with the constants $n_0(\epsilon, D)$ in the polynomial stochastic domination depending only on $c, C$, and the parameters $(\gamma_W, \xi_W)$ of the generalized Wigner matrix.
\end{thm}
\begin{proof}
    For $\Im(z)>0$, the result follows from the cited theorem: the error away from the spectrum in the upper half-plane is of order $n^{-1/2}$, and Remark 2.6 of the reference makes this estimate uniform in $z$.
    Next, by Proposition~\ref{prop:wig-norm}, with polynomially high probability the spectrum of $W$ is disjoint from $\Omega\cap\RR$.
    On this event, both $R_W(z)$ and $G_{\mu_{\sc}}(z)$ are continuous as $\Im(z) \downarrow 0$, so the same estimate also holds on the real axis by the limiting argument of Remark 2.7 of \cite{bloemendal2015isotropiclocallawssample}.
\end{proof}

\begin{cor}
    \label{cor:ill-Q}
    Let $W \in \CC^{n \times n}_{\herm}$ be a generalized Wigner matrix.
    For some $i, j \in [n]$, let $Q$ be formed by setting entries $(i, j)$ and $(j, i)$ in $W$ to zero.
    Then, the conclusion \eqref{eq:ill-dom} of Theorem~\ref{thm:ill} holds with $W$ replaced by $Q$.
\end{cor}
\begin{proof}
    It suffices to show that
    \[ \sup_{z \in \Omega} \|R_W(z) - R_Q(z)\| \prec n^{-1/2}. \]
    By the second resolvent identity, we have
    \begin{align*}
    \|R_W(z) - R_Q(z)\|
    &= \|R_W(z)(W - Q)R_Q(z)\| \\
    &\leq |W_{ij}| \cdot \|R_W(z)\| \cdot \|R_Q(z)\| \\
    &\leq |W_{ij}| \cdot \|R_W(z)\| \cdot (\|R_W(z)\| + \|R_W(z) - R_Q(z)\|),
    \intertext{and rearranging this gives}
    &\leq \frac{|W_{ij}| \cdot \|R_W(z)\|^2}{(1 - |W_{ij}| \cdot \|R_W(z)\|)}
    \end{align*}
    provided that $|W_{ij}| \cdot \|R_W(z)\| < 1$.
    We have $|W_{ij}| \prec n^{-1/2}$ since by assumption it is sub-Weibull, while it follows from Proposition~\ref{prop:wig-norm} that $\sup_z \|R_W(z)\| \prec 1$ over the range of $z$ in the statement, and the result follows by a suitable union bound.
\end{proof}

\noindent
We will also use the following corollary, which essentially says that one may differentiate with respect to $z$ inside the expression appearing in the local law and still obtain the same quality of guarantee.
\begin{cor}[Derivative of isotropic local law]\label{cor:dill}
    In the setting of Theorem~\ref{thm:ill}, we also have
    \[ \sup_{z \in \Omega} \left|x^*R_W(z)^2y + G_{\mu_{\sc}}^{\prime}(z) x^*y\right| \prec n^{-1/2}, \]
    with the same dependence of constants.
\end{cor}
\begin{proof}
    For deterministic unit vectors $a,b \in \CC^n$, define
    \[ f_{a,b}(z) \colonequals a^*R_W(z)b-G_{\mu_{\sc}}(z)a^*b. \]
    Let
    \[ \Omega^{\prime} \colonequals \{z \in \CC: \Re(z)>2+c/2, |z|\leq C+c/2\}. \]
    By Proposition~\ref{prop:wig-norm}, with polynomially high probability $\spec(W)$ is disjoint from $\Omega^{\prime}$, so $f_{a,b}$ is analytic on an open neighborhood of this domain.
    Theorem~\ref{thm:ill}, applied with parameters $c/2$ and $C+c/2$, controls $f_{x,y}$ and $f_{y,x}$ on the part of $\Omega^{\prime}$ in the closed upper half-plane.
    For $z$ in the lower half-plane, Hermitian symmetry gives
    \[ f_{x,y}(z)=\overline{f_{y,x}(\overline{z})}. \]
    Consequently,
    \[ \sup_{z\in\Omega^{\prime}}|f_{x,y}(z)|\prec n^{-1/2}. \]
    For any $w\in\Omega$, let $C_w$ be the circle of radius $c/4$ centered at $w$.
    This circle is contained in $\Omega^{\prime}$, including when $w$ lies on or near the real axis.
    Cauchy's integral formula therefore gives
    \[ \sup_{w\in\Omega}|f_{x,y}^{\prime}(w)|
    =\sup_{w\in\Omega}\left|\frac{1}{2\pi i}\oint_{C_w}\frac{f_{x,y}(z)}{(z-w)^2}\,dz\right|
    \prec n^{-1/2}. \]
    Since $\frac{d}{dz}R_W(z)=-R_W(z)^2$, we have
    \[ f_{x,y}^{\prime}(z)=-x^*R_W(z)^2y-G_{\mu_{\sc}}^{\prime}(z)x^*y, \]
    which proves the result.
\end{proof}

\subsubsection{Spiked matrix model estimates}

Combining the above results, we may prove estimates on the largest eigenvalue and associated eigenvector of spiked matrix models with generalized Wigner matrix noise that will be useful later.
These arguments for deducing such bounds from isotropic local laws are standard; we outline them omitting some details for the sake of completeness.

\begin{thm}
    \label{thm:spiked-non-asymp}
    Let $W$ be a generalized Wigner matrix, $\theta > 1$, and $v \in \SS^{n - 1}(\CC)$ be deterministic.
    Write $\what{\lambda}$ for the largest eigenvalue of $H = \theta vv^* + W$ and $\what{v}$ for the associated eigenvector.
    Then, we have
    \begin{align*}
        |\what{\lambda} - \lambda(\theta)| &\prec n^{-1/2}, \\
        \big| |\langle \what{v}, v \rangle|^2 - \rho(\theta)^2 \big| &\prec n^{-1/2},
    \end{align*}
    with the constants in the polynomial stochastic domination depending only on $\theta$ and the parameters $(\gamma_W, \xi_W)$ of the generalized Wigner matrix.
    Further, for any $\epsilon > 0$, with polynomially high probability $H$ has at most one eigenvalue in the interval $[2 + \epsilon, \infty)$.
\end{thm}
\begin{proof}
    Write $\lambda=\lambda(\theta)=\theta+\theta^{-1}$ and choose $\delta>0$ such that
    \[ I\colonequals[\lambda-\delta,\lambda+\delta]\subset(2,\infty). \]
    By Proposition~\ref{prop:wig-norm}, with polynomially high probability $\lambda_1(W)<\lambda-\delta$.
    On this event, the function
    \[ m_W(E)\colonequals v^*R_W(E)v \]
    is defined and strictly decreasing on $I$, since $m_W^{\prime}(E)=-v^*R_W(E)^2v<0$.
    Theorem~\ref{thm:ill} gives
    \[ \sup_{E\in I}|m_W(E)-G_{\mu_{\sc}}(E)|\prec n^{-1/2}. \]
    Proposition~\ref{prop:Gsc} shows that $G_{\mu_{\sc}}$ is strictly decreasing on $I$ and
    \[ G_{\mu_{\sc}}(\lambda)=\frac{1}{\theta}. \]
    The values of $G_{\mu_{\sc}}$ at the two endpoints of $I$ lie on opposite sides of $1/\theta$, with fixed gaps.
    The uniform local law therefore gives the same inequalities for $m_W$ with polynomially high probability, so continuity and strict monotonicity show that the equation $m_W(E)=1/\theta$ has a unique solution in $I$.
    By Proposition~\ref{prop:res-eig}, this solution is an eigenvalue of $H$.
    The interlacing inequality gives $\lambda_2(H)\leq\lambda_1(W)$, so this solution is the simple largest eigenvalue $\what{\lambda}$ of $H$.
    Since $G_{\mu_{\sc}}^{\prime}$ is bounded away from zero on $I$, the mean value theorem and the uniform local law give
    \[ |\what{\lambda}-\lambda|\prec n^{-1/2}. \]

    Proposition~\ref{prop:res-eigvec} and the secular equation give
    \[ |\langle\what{v},v\rangle|^2
    =\frac{(v^*R_W(\what{\lambda})v)^2}{v^*R_W(\what{\lambda})^2v}
    =\frac{\theta^{-2}}{v^*R_W(\what{\lambda})^2v}. \]
    Corollary~\ref{cor:dill}, evaluated uniformly at $\what{\lambda}\in I$, yields
    \[ \left|v^*R_W(\what{\lambda})^2v-\frac{1}{\theta^2-1}\right|\prec n^{-1/2}, \]
    where we also used $|\what{\lambda}-\lambda|\prec n^{-1/2}$ and the smoothness of $G_{\mu_{\sc}}^{\prime}$ on $I$.
    Therefore,
    \[ \big||\langle\what{v},v\rangle|^2-(1-\theta^{-2})\big|\prec n^{-1/2}, \]
    which is the claimed estimate because $\rho(\theta)^2=1-\theta^{-2}$.

    Finally, for any fixed $\epsilon>0$, Proposition~\ref{prop:wig-norm} gives $\lambda_1(W)<2+\epsilon$ with polynomially high probability.
    The interlacing inequality $\lambda_2(H)\leq\lambda_1(W)$ then shows that $H$ has at most one eigenvalue in $[2+\epsilon,\infty)$.
\end{proof}

\section{Universality of entrywise statistics: Proof of Theorem~\ref{thm:main-univ}}

Recall the setting: we have two tuples of generalized Wigner matrices,
\begin{align*}
  \bW &= (W^{(1)},\ldots, W^{(L)}), \\
  \bX &= (X^{(1)},\ldots, X^{(L)}),
\end{align*}
where each $W^{(a)}, X^{(a)} \in \CC^{n \times n}_{\herm}$ for $a\in[L]$.
For $1\le p\le q\le n$, we write
\[w^{(pq)} = (\Re W^{(1)}_{pq}, \Im W^{(1)}_{pq}, \ldots, \Re W^{(L)}_{pq}, \Im W^{(L)}_{pq}) \in \RR^{2L}.\]
Thus the whole tuple $\bW$ is equivalently encoded by
$\bw=\left(w^{(pq)}\right)_{1\le p\le q\le n}$.

For $p<q$, the $(p,q)$-th entry of $a$-th matrix  is recovered as  \[W^{(a)}_{pq} = \Re W^{(a)}_{pq} + \cpxi \Im W^{(a)}_{pq} = w^{(pq)}_{2a-1} + \cpxi w^{(pq)}_{2a}\]
The lower-triangular entries are determined by Hermitian symmetry, and on the diagonal we have
\begin{align*}
  W^{(a)}_{pp} &= w^{(pp)}_{2a-1}, \\
  w^{(pp)}_{2a} &= 0
\end{align*}
For each $a\in[L]$, we are given a deterministic vector $v^{(a)} \in \SS^{n - 1}(\CC)$.
On $v^{(a)}$ we have only assumed that $\|v^{(a)}\| = 1$ and that, for each $a \in [L]$,
\[ \|v^{(a)}\|_{\infty} \leq C_v n^{-1/2+\epsilon_v}, \]
for some $\epsilon_v \in (0, 1/20)$.
For $a\in[L]$, we define
\begin{align*}
    H^{(W,a)} &= H(W^{(a)}) = \theta v^{(a)}v^{(a)*} + W^{(a)}, \\
    \what{\lambda}^{(W,a)} &= \what{\lambda}(W^{(a)}) = \lambda_1(H(W^{(a)})), \\
    \what{v}^{(W,a)} &= \what{v}(W^{(a)}) = v_1(H(W^{(a)})), \\
    F_{ij}(\bw) &= \phi\left(n \cdot \what{v}^{(W,1)}_i \cdot \overline{\what{v}^{(W,1)}_j}, \ldots,n \cdot \what{v}^{(W,L)}_i \cdot \overline{\what{v}^{(W,L)}_j}\right)
\end{align*}
for some fixed $i,j\in [n]$ and $\phi: \CC^{L} \to \RR$ five times continuously differentiable and with bounded values and derivatives. We view $\phi$ as a function on $\RR^{2L}$.
The analogous statements hold for the $\bX$-tuple, with $W,w$ replaced by $X,x$.

We will use the notation defined earlier
\[ \lambda = \lambda(\theta) \colonequals \theta + \theta^{-1} > 2. \]

Our goal is to show that, uniformly in $i,j$, \[\EE F_{ij}(\bw) \approx \EE F_{ij}(\bx).\]
We proceed in two steps.
First, we make some initial simplifications using our calculations related to the eigenspaces of $H$ from Proposition~\ref{prop:res-eigvec}.
In particular, we reduce $F_{ij}(\bw)$ to a simpler resolvent statistic $\wtilde{F}_{ij}(\bw)$.
Then, we complete the remaining proof by the Lindeberg method.

\subsection{Initial simplifications}
\label{sec:initial-simp}

Let $\epsilon \in (0, \epsilon_v)$ be a constant for the purposes of the proof to be chosen later, and fix $0 < \delta_1 < \delta_2$ such that $\lambda \in (2 + \delta_1, 2 + \delta_2)$.
We view $\epsilon, \delta_1, \delta_2, L$ as well as the function $\phi$ and the generalized Wigner matrix parameters $(\gamma_W, \xi_W, \gamma_X, \xi_X)$, all as constants for the purposes of this proof, and the dependence of various asymptotic notations on these parameters is not mentioned from now on.

For each $a\in[L]$, let $\mathcal{E}^{(a)}$ be the event that the following conditions hold for $W^{(a)}$:
\begin{enumerate}
    \item There is a single eigenvalue of $H^{(W,a)}$ in $[2 + \delta_1, 2 + \delta_2)$, of multiplicity 1, which is also the top eigenvalue $\what{\lambda}^{(W,a)} = \lambda_1(H(W^{(a)}))$.
    \item $|\what{\lambda}^{(W,a)} - \lambda| \leq n^{-1/2 + \epsilon}$.
    \item $\|W^{(a)}\| < 2 + \delta_1$, and in particular $\what{\lambda}^{(W,a)} \notin \spec(W^{(a)})$.
    \item The following bounds hold:
    \begin{align*}
        |e_i^* R_{W^{(a)}}(\what{\lambda}^{(W,a)}) v^{(a)} - G_{\mu_{\sc}}(\what{\lambda}^{(W,a)}) v^{(a)}_i| &\leq n^{-1/2 + \epsilon}, \\
        |v^{(a)*} R_{W^{(a)}}(\what{\lambda}^{(W,a)}) e_j - G_{\mu_{\sc}}(\what{\lambda}^{(W,a)}) \overline{v^{(a)}_j}| &\leq n^{-1/2 + \epsilon}, \\
        |v^{(a)*} R_{W^{(a)}}(\what{\lambda}^{(W,a)})^2 v^{(a)} + G_{\mu_{\sc}}^{\prime}(\what{\lambda}^{(W,a)})| &\leq n^{-1/2 + \epsilon}.
    \end{align*}
\end{enumerate}

\begin{prop}
    \label{prop:E-prob}
    $\sE \colonequals \bigcap_{a=1}^{L}\sE^{(a)}$ holds with polynomially high probability.
\end{prop}
\begin{proof}
    The result follows by combining the isotropic local law (Theorem~\ref{thm:ill}), its derivative (Corollary~\ref{cor:dill}), and Theorem~\ref{thm:spiked-non-asymp} on the eigenvalues of spiked matrix models. For each fixed $a\in[L]$, this follows from the corresponding scalar estimates. Since $L$ is fixed, taking the intersection over $a\in[L]$ only changes the implicit constants.
\end{proof}

Define
\begin{align*}
    \Lambda_1 &\colonequals -G_{\mu_{\sc}}^{\prime}(2 + \delta_2), \\
    \Lambda_2 &\colonequals -G_{\mu_{\sc}}^{\prime}(2 + \delta_1).
\end{align*}
By Proposition~\ref{prop:Gsc}, we have $0 < \Lambda_1 < \Lambda_2$ and $-G_{\mu_{\sc}}^{\prime}(\lambda) \in [\Lambda_1, \Lambda_2]$ since $-G_{\mu_{\sc}}^{\prime}(z)$ is monotonically decreasing on the interval $(2, \infty)$.

Since $\phi$ is bounded, we have by Proposition~\ref{prop:E-prob} the bound:
\begin{equation}
    \label{eq:f-trunc}
    |\EE F_{ij}(\bw) - \EE \boldone_{\mathcal{E}} F_{ij}(\bw)| \leq \|\phi\|_{L^{\infty}} \PP[\mathcal{E}^c] \lesssim n^{-1}.
\end{equation}
Applying Proposition~\ref{prop:res-eigvec} to each coordinate $a\in[L]$, we obtain:
\begin{align*}
    \EE \boldone_{\mathcal{E}} F_{ij}(\bw)
    &= \EE \boldone_{\mathcal{E}} \phi\left( \left(n \cdot \what{v}^{(W,a)}_i \cdot \overline{\what{v}^{(W,a)}_j}\right)_{a=1}^{L}\right) \\
    &= \EE \boldone_{\mathcal{E}} \phi\left( \left(n \cdot \frac{\left(e_i^* R_{W^{(a)}}(\what{\lambda}^{(W,a)}) v^{(a)}\right)\left(v^{(a)*} R_{W^{(a)}}(\what{\lambda}^{(W,a)}) e_j\right)}{v^{(a)*} R_{W^{(a)}}(\what{\lambda}^{(W,a)})^2 v^{(a)}} \right)_{a=1}^{L}\right).
\end{align*}

We will now show that $\EE F_{ij}(\bw)$ is close to what we get when we perform two operations on the above expression: (1) we replace the random eigenvalues $\what{\lambda}^{(W,a)}$ by their deterministic typical location $\lambda$ (which is the same for each $a \in [L]$), and (2) we apply the isotropic local law to the denominator.
There is a small but important caveat that we must attend to:

\begin{rmk}[Resolvent for discrete models]
    \label{rem:res-cpx}
    When the entries of $W^{(a)}$ have discrete distributions, which our Definition~\ref{def:gen-wigner} does not rule out, it is possible that $\lambda \in \spec(W^{(a)})$ with small but positive probability.
    On this event, $R_{W^{(a)}}(\lambda)$ is undefined, and in particular expectations involving $R_{W^{(a)}}(\lambda)$ are undefined.
\end{rmk}
\noindent
To deal with this, consider the slight further perturbation
\[ \widetilde{\lambda} = \widetilde{\lambda}(\theta) = \lambda(\theta) + \cpxi n^{-1/2}. \]
For the above technical reason, we instead consider replacing $\what{\lambda}^{(W,a)}$ by $\wtilde{\lambda}$.

It is convenient to state this reduction in terms of the following auxiliary function:
\[ \wtilde{\phi}(t_1, \ldots, t_L) \colonequals \phi\left(\frac{-1}{G^{\prime}_{\mu_{\sc}}(\wtilde{\lambda})} t_1,\dots, \frac{-1}{G^{\prime}_{\mu_{\sc}}(\wtilde{\lambda})} t_L\right). \]
Since $-1 / G^{\prime}_{\mu_{\sc}}(\wtilde\lambda)$ is uniformly bounded for all sufficiently large $n$, $\wtilde{\phi}$ satisfies the same boundedness and smoothness assumptions as $\phi$ up to modifying the constants involved (in a way not depending on $n$).
We then define
\begin{equation}
\label{eq:tilde-F}
\wtilde{F}_{ij}(\bw) \colonequals \wtilde{\phi}\left( \left(n \cdot (e_i^* R_{W^{(a)}}(\wtilde{\lambda}) v^{(a)})(v^{(a)*} R_{W^{(a)}}(\wtilde{\lambda}) e_j) \right)_{a=1}^{L} \right).
\end{equation}
\begin{lem}
    \label{lem:F-bound}
    For all $i,j\in[n]$, the following bounds hold:
    \begin{align*}
        |\EE F_{ij}(\bw) - \EE \wtilde{F}_{ij}(\bw)| &\lesssim n^{-1/2 + 3\epsilon_v}, \\
        |\EE F_{ij}(\bx) - \EE \wtilde{F}_{ij}(\bx)| &\lesssim n^{-1/2 + 3\epsilon_v},
        \intertext{and therefore also}
        |\EE F_{ij}(\bw) - \EE F_{ij}(\bx)| &\leq |\EE \wtilde{F}_{ij}(\bw) - \EE \wtilde{F}_{ij}(\bx)| + O(n^{-1/2 + 3\epsilon_v}).
    \end{align*}
\end{lem}
\begin{proof}
We will show the first bound for $\bw$; the same argument will apply symmetrically for $\bx$, and the two bounds combined will give the third bound.

Throughout the proof, we further intersect $\sE$ with the polynomially high probability event on which
Corollary~\ref{cor:dill} holds uniformly on the fixed domains used below, for the finitely many deterministic pairs
\[(x,y)\in \{(e_i,v^{(a)}), (v^{(a)},e_j): a\in[L]\}.\] Since $L$ is fixed, this new event still holds with
polynomially high probability, and we keep the
notation $\sE$ for it.

We first replace the denominator by its deterministic equivalent. By the conditions included in $\mathcal{E}$, we have that, on this event, for $a\in[L]$,
\begin{align*}
&\left|\frac{1}{v^{(a)*} R_{W^{(a)}}(\what{\lambda}^{(W,a)})^2 v^{(a)}} + \frac{1}{G^{\prime}_{\mu_{\sc}}(\wtilde{\lambda})}\right| \\
&\hspace{2cm}\leq \left|\frac{1}{v^{(a)*} R_{W^{(a)}}(\what{\lambda}^{(W,a)})^2 v^{(a)}} + \frac{1}{G^{\prime}_{\mu_{\sc}}(\what{\lambda}^{(W,a)})}\right| + \left|\frac{1}{G^{\prime}_{\mu_{\sc}}(\wtilde{\lambda})} -  \frac{1}{G^{\prime}_{\mu_{\sc}}(\what{\lambda}^{(W,a)})}\right| \\
&\hspace{2cm} = \frac{|v^{(a)*} R_{W^{(a)}}(\what{\lambda}^{(W,a)})^2 v^{(a)} + G^{\prime}_{\mu_{\sc}}(\what{\lambda}^{(W,a)})|}{|v^{(a)*} R_{W^{(a)}}(\what{\lambda}^{(W,a)})^2 v^{(a)}| \cdot |G^{\prime}_{\mu_{\sc}}(\what{\lambda}^{(W,a)})|} + O(n^{-1/2 + \epsilon}) \\
&\hspace{2cm} \leq \frac{n^{-1/2 + \epsilon}}{(\Lambda_1 - n^{-1/2 + \epsilon})\Lambda_1} + O(n^{-1/2 + \epsilon}) \\
&\hspace{2cm} = O(n^{-1/2 + \epsilon}),
\end{align*}
where we also use that $1 / G^{\prime}_{\mu_{\sc}}(z)$ is Lipschitz away from $\{\pm 2\}$, per Proposition~\ref{prop:Gsc}.

Moreover, on $\sE$, the isotropic local law bounds of Theorem~\ref{thm:ill} give
\[
|e_i^* R_{W^{(a)}}(\what{\lambda}^{(W,a)}) v^{(a)}| \le |G_{\mu_{\sc}}(\what{\lambda}^{(W,a)})||v_i^{(a)}| + n^{-1/2+\epsilon} \lesssim n^{-1/2+\epsilon_v}.
\]
Similarly
\[
|v^{(a)*} R_{W^{(a)}}(\what{\lambda}^{(W,a)}) e_j|\lesssim n^{-1/2+\epsilon_v}.
\]
Therefore, we also have that on $\mathcal{E}$,
\begin{align*}
&\Bigg|n \cdot \frac{\left(e_i^* R_{W^{(a)}}(\what{\lambda}^{(W,a)}) v^{(a)}\right) \left(v^{(a)*} R_{W^{(a)}}(\what{\lambda}^{(W,a)}) e_j\right)}{v^{(a)*} R_{W^{(a)}}(\what{\lambda}^{(W,a)})^2 v^{(a)}} \\
&\hspace{1cm} + n \cdot \frac{\left(e_i^* R_{W^{(a)}}(\what{\lambda}^{(W,a)}) v^{(a)}\right) \left(v^{(a)*} R_{W^{(a)}}(\what{\lambda}^{(W,a)}) e_j\right)}{G_{\mu_{\sc}}^{\prime}(\wtilde{\lambda})}\Bigg|\\ &\hspace{2cm} \lesssim n\cdot n^{-1+2\epsilon_v} \cdot n^{-1/2+\epsilon}
= O(n^{-1/2 + \epsilon + 2\epsilon_v}).
\end{align*}

Since $\phi$ is Lipschitz as a function on $\RR^{2L}$ and $L$ is fixed, it follows that
\begin{align*}
\EE \boldone_{\sE}F_{ij}(\bw) &=
\EE \boldone_{\sE}\phi\left( \left(n \cdot
\frac{-1}{G_{\mu_{\sc}}^{\prime}(\wtilde{\lambda})}(e_i^*R_{W^{(a)}}(\what{\lambda}^{(W,a)})v^{(a)})
(v^{(a)*}R_{W^{(a)}}(\what{\lambda}^{(W,a)})e_j)
\right)_{a=1}^L
\right) \\
&\hspace{2cm} + O(n^{-1/2+\epsilon+2\epsilon_v}).
\end{align*}

We next replace $\what{\lambda}^{(W,a)}$ by $\wtilde{\lambda}$ in the two numerator factors. Decompose
\begin{align*} |e_i^*R_{W^{(a)}}(\wtilde{\lambda})v^{(a)} &- e_i^*R_{W^{(a)}}(\what{\lambda}^{(W,a)})v^{(a)}| \le |G_{\mu_{\sc}}(\wtilde{\lambda}) - G_{\mu_{\sc}}(\what{\lambda}^{(W,a)})||v^{(a)}_i|\\ &+ \left| \left(e_i^*R_{W^{(a)}}(\wtilde{\lambda})v^{(a)} - G_{\mu_{\sc}}(\wtilde{\lambda})v^{(a)}_i\right) - \left(e_i^*R_{W^{(a)}}(\what{\lambda}^{(W,a)})v^{(a)} - G_{\mu_{\sc}}(\what{\lambda}^{(W,a)}) v^{(a)}_i\right) \right|.
\end{align*}
Consider the neighborhood \[\wtilde{\Omega} \colonequals \{z \in \CC: |z-\lambda|\le \frac{1}{4}(\lambda -(2+\delta_1))\}.\]
Choose $c,C>0$ such that \[\wtilde{\Omega} \cap \{\Im z>0\} \subset \Omega(c,C) = \{z\in\CC:\Re z>2+c, \Im z>0, |z|\le C\}\]
On $\sE$, for sufficiently large $n$, the segment joining $\what{\lambda}^{(W,a)}$ and $\wtilde{\lambda}$ except possibly the real endpoint $\what{\lambda}^{(W,a)}$ is contained in $\wtilde{\Omega} \cap \{\Im z>0\}$ and is separated from $\spec(W^{(a)})$ by a deterministic positive distance. The derivative local law, Corollary~\ref{cor:dill}, applies on the part of the segment with positive imaginary part, and the same bound holds at the real endpoint by continuity.

By the fundamental theorem of calculus along the segment
$[\what{\lambda}^{(W,a)},\wtilde{\lambda}]$, and by Corollary~\ref{cor:dill},\begin{align*}
&\left| \left(e_i^*R_{W^{(a)}}(\wtilde{\lambda})v^{(a)} - G_{\mu_{\sc}}(\wtilde{\lambda})v^{(a)}_i\right) - \left(e_i^*R_{W^{(a)}}(\what{\lambda}^{(W,a)})v^{(a)} - G_{\mu_{\sc}}(\what{\lambda}^{(W,a)})v^{(a)}_i\right) \right| \\ &\qquad \le |\wtilde{\lambda} -\what{\lambda}^{(W,a)}| \sup_{z\in [\what{\lambda}^{(W,a)},\wtilde{\lambda}]} \left|e^*_i R_{W^{(a)}}(z)^2 v^{(a)} + G_{\mu_{\sc}}^{\prime}(z)v^{(a)}_i\right| \lesssim |\wtilde{\lambda} -\what{\lambda}^{(W,a)}| \cdot n^{-1/2+\epsilon}.
\end{align*}
Using $G_{\mu_{\sc}}$ is Lipschitz and $|\wtilde{\lambda} -\what{\lambda}^{(W,a)}| \lesssim n^{-1/2+\epsilon}$, we obtain
\begin{align*}
|e_i^*R_{W^{(a)}}(\wtilde{\lambda})v^{(a)} - e_i^*R_{W^{(a)}}(\what{\lambda}^{(W,a)})v^{(a)}| &\le |G_{\mu_{\sc}}(\wtilde{\lambda}) - G_{\mu_{\sc}}(\what{\lambda}^{(W,a)})||v^{(a)}_i| + O(|\wtilde{\lambda} -\what{\lambda}^{(W,a)}| n^{-1/2+\epsilon})\\ &\lesssim n^{-1/2+\epsilon} n^{-1/2+\epsilon_v} + n^{-1/2+\epsilon} n^{-1/2+\epsilon}\\ &\lesssim n^{-1+\epsilon+\epsilon_v}
\end{align*}
Also applying this symmetrically to the $v^{(a)*} R_{W^{(a)}}(\what{\lambda}^{(W,a)}) e_j$ term, we find that,
\begin{align*}
&n\left|(e_i^*R_{W^{(a)}}(\wtilde{\lambda})v^{(a)})(v^{(a)*}R_{W^{(a)}}(\wtilde{\lambda})e_j)-(e_i^*R_{W^{(a)}}(\what{\lambda}^{(W,a)})v^{(a)})
(v^{(a)*}R_{W^{(a)}}(\what{\lambda}^{(W,a)})e_j)
\right| \\
&\qquad\lesssim
n\cdot n^{-1+\epsilon+\epsilon_v}\cdot n^{-1/2+\epsilon_v}
+
n\cdot n^{-1+\epsilon+\epsilon_v}\cdot n^{-1+\epsilon+\epsilon_v} \\
&\qquad\lesssim
n^{-1/2+\epsilon+2\epsilon_v}.
\end{align*}
Since $\phi$ is smooth and $\phi^{\prime}$ is bounded, and combining with \eqref{eq:f-trunc},
\begin{align*}
\EE F_{ij}(\bw)
&= \EE \boldone_{\mathcal{E}} \phi \left( \left(n \cdot \frac{-1}{G^{\prime}_{\mu_{\sc}}(\wtilde{\lambda})} \cdot (e_i^* R_{W^{(a)}}(\wtilde{\lambda}) v^{(a)})(v^{(a)*} R_{W^{(a)}}(\wtilde{\lambda}) e_j)\right)_{a=1}^{L} \right) + O(n^{-1/2 + \epsilon + 2\epsilon_v})
\intertext{and, undoing our first truncation step where we introduced the indicator of $\mathcal{E}$, we have}
&= \EE \phi\left( \left(n \cdot \frac{-1}{G^{\prime}_{\mu_{\sc}}(\wtilde{\lambda})} \cdot (e_i^* R_{W^{(a)}}(\wtilde{\lambda}) v^{(a)})(v^{(a)*} R_{W^{(a)}}(\wtilde{\lambda}) e_j)\right)_{a=1}^{L} \right) + O(n^{-1/2 + \epsilon + 2\epsilon_v}) \\
&= \EE \wtilde{\phi}\left( \left(n \cdot (e_i^* R_{W^{(a)}}(\wtilde{\lambda}) v^{(a)})(v^{(a)*} R_{W^{(a)}}(\wtilde{\lambda}) e_j)\right)_{a=1}^{L} \right) + O(n^{-1/2 + \epsilon + 2\epsilon_v}),
\label{eq:EFW-inter}
\numberthis
\end{align*}
completing the proof, where we note that the stated error bound follows since we have taken $\epsilon < \epsilon_v$ by assumption.
\end{proof}

\subsection{Lindeberg method}

Now, we consider comparing $\EE \wtilde{F}_{ij}(\bw)$ and $\EE \wtilde{F}_{ij}(\bx)$ by the Lindeberg method.
Our final result will be the following:
\begin{lem}
    \label{lem:tilde-F-bound}
    For $\wtilde{F}_{ij}$ as defined above in \eqref{eq:tilde-F}, for any $\epsilon > 0$,
    \[ \max_{i,j\in[n]}|\EE \wtilde{F}_{ij}(\bx) - \EE \wtilde{F}_{ij}(\bw)| \lesssim_{\epsilon} n^{-1/2 + 10\epsilon_v + \epsilon}. \]
\end{lem}
\noindent
The proof of Theorem~\ref{thm:main-univ} is then completed by combining Lemma~\ref{lem:F-bound} with Lemma~\ref{lem:tilde-F-bound}.

Following the prescription of the Lindeberg method (our use of the method is also quite similar to the specific one in \cite{knowles2011eigenvectordistributionwignermatrices}), let us construct an interpolating path from tuple $\bW$ to tuple $\bX$.
Let $N = \frac{n(n + 1)}{2}$ be the number of entries on and above the diagonal in a Hermitian matrix.
Fix some enumeration $(p_1, q_1), \dots, (p_N, q_N)$ of all pairs $(p_\ell, q_\ell)$ with $1 \leq p \leq q \leq n$. Throughout this argument, $i,j\in[n]$ are fixed and refer to the entry statistic $F_{ij}$. By contrast, at the $\ell$-th Lindeberg step we write $p_\ell, q_\ell$ for the matrix entry currently being swapped.

For $0\le \ell \le N$, let
\[\bW^{(\ell)} = \left(W^{(\ell,1)},\ldots, W^{(\ell,L)}\right)\]
be the tuple obtained from $\bW$ by replacing the entries $W^{(a)}_{p_1q_1}, \dots, W^{(a)}_{p_{\ell}q_{\ell}}$ with the corresponding entries of $\bX$ for all $a\in[L]$, and by their conjugates symmetrically below the diagonal.
Thus, $\bW^{(0)} = \bW$, $\bW^{(N)} = \bX$.
We emphasize that, while we work with $L$ different Hermitian matrices with a total of $L\cdot N$ entries (possibly complex-valued), our Lindeberg procedure only involves $N$ steps: for each fixed upper-triangular coordinate $(p_{\ell},q_{\ell})$, in one step the whole vector $w^{(p_{\ell}q_{\ell})}\in\RR^{2L}$ is replaced by $x^{(p_{\ell}q_{\ell})}\in\RR^{2L}$.
We swap the whole block at once because the coordinates across the $L$ matrices may be correlated within the same upper-triangular position, while the blocks are independent across different positions. This also ensures that, after the current block is removed, the tuple $\bQ^{(\ell)}$ introduced below is independent of both $w^{(p_{\ell}q_{\ell})}$ and $x^{(p_{\ell}q_{\ell})}$.

We follow the usual prescription of the Lindeberg method.
First, we bound by a telescoping sum,
\begin{align*}
    |\EE \wtilde{F}_{ij}(\bx) - \EE \wtilde{F}_{ij}(\bw)|
    &= \left|\sum_{\ell = 1}^N \left[\EE \wtilde{F}_{ij}(\bW^{(\ell)}) - \EE \wtilde{F}_{ij}(\bW^{(\ell - 1)})\right]\right| \\
    &\leq \sum_{\ell = 1}^N |\EE \wtilde{F}_{ij}(\bW^{(\ell)}) - \EE \wtilde{F}_{ij}(\bW^{(\ell - 1)})|. \numberthis \label{eq:lindeberg-triangle}
\end{align*}

We now consider bounding the effect of each individual swap, i.e., the size of
\begin{equation}
\label{eq:diff-ell}
|\EE \wtilde{F}_{ij}(\bW^{(\ell)}) - \EE \wtilde{F}_{ij}(\bW^{(\ell - 1)})|
\end{equation}
for each fixed $\ell$.
For the sake of brevity, let us write $p = p_{\ell}$ and $q = q_{\ell}$ while analyzing a single step.
Define
\[ \bE^{(W)}_{pq} \colonequals \left(E^{(W,1)}_{pq},\ldots, E^{(W,L)}_{pq}\right), \qquad
\begin{aligned}
    E^{(W,a)}_{pq} &\colonequals \left\{\begin{array}{ll} W^{(a)}_{pq} e_p e_q^{*} + \overline{W^{(a)}_{pq}} e_q e_p^* & \text{if } p \neq q, \\
    W^{(a)}_{pp} e_p e_p^* & \text{if } p = q.\end{array}\right\},
\end{aligned}\]
and similarly $\bE^{(X)}_{pq}$, $E^{(X,a)}_{pq}$.
Then the original tuples decompose as
\[W^{(a)} = \sum_{1\leq \mu \leq \nu \leq n}E^{(W,a)}_{\mu\nu}, \qquad X^{(a)} = \sum_{1\leq \mu \leq \nu \leq n}E^{(X,a)}_{\mu\nu},\qquad a\in[L].\]
Let
\[ Q^{(\ell,a)} = W^{(\ell,a)} - E^{(X,a)}_{pq} = W^{(\ell - 1, a)} - E^{(W,a)}_{pq}, \qquad a\in[L],\]
and write \[\bQ^{(\ell)}\colonequals \left(Q^{(\ell,1)},\ldots, Q^{(\ell,L)}\right).\]
For notational convenience, we suppress the superscript $\ell$ and simply write \[\bQ = \left(Q^{(1)},\ldots,Q^{(L)}\right).\]
Thus, \[\bW^{(\ell)} = \bQ + \bE^{(X)}_{pq},\qquad \bW^{(\ell-1)}=\bQ+\bE^{(W)}_{pq}.\]
In words, $\bQ$ is the tuple after the first $\ell-1$ swaps have been performed but where the entries about to be swapped at step $\ell$ (in positions $(p = p_{\ell}, q = q_{\ell})$ and $(q, p)$) have been set to zero.

It is useful to pause to understand the distribution of $\bQ$. For each fixed $a\in[L]$, the matrix $Q^{(a)}$ is obtained by taking the entries indexed by $(p_r, q_r)$ with $r<\ell$ from $X^{(a)}$, the entries indexed by $(p_r,q_r)$ with $r>\ell$ from $W^{(a)}$, and setting the current entries in positions $(p = p_{\ell}, q = q_{\ell})$ and $(q,p)$ to zero.  Adding back either $E^{(W,a)}_{pq}$ or $E^{(X,a)}_{pq}$ gives respectively $W^{(\ell-1,a)}$ or $W^{(\ell,a)}$, both of which are generalized Wigner matrices with the same variance normalization, since the second moments of $w^{(pq)}$ and $x^{(pq)}$ match. Hence Corollary~\ref{cor:ill-Q} applies to each $Q^{(a)}$. Since $L$ is fixed, we may use these local law bounds simultaneously for all $a\in[L]$.

Moreover, by independence over upper-triangular index pairs, $\bQ$ is independent of the current entry vectors $w^{(pq)},x^{(pq)}\in\RR^{2L}$.

Rewriting our expression in \eqref{eq:diff-ell} in terms of this $\bQ$, we have
\[ \left|\EE \wtilde{F}_{ij}(\bW^{(\ell)}) - \EE \wtilde{F}_{ij}(\bW^{(\ell - 1)})\right| = \left|\EE \wtilde{F}_{ij}\left(\bQ + \bE^{(X)}_{pq}\right) - \EE \wtilde{F}_{ij}\left(\bQ + \bE^{(W)}_{pq}\right)\right|, \]
where the addition of tuples is understood coordinatewise.

We would like to understand the leading order effect of this difference, viewed as a small perturbation.
We do this by first taking a perturbative expansion of each term separately.
Recall that $\wtilde{F}_{ij}$ is a function of $\bW$ only through the resolvent value $R_{W^{(a)}}(\wtilde{\lambda})$ for $a\in[L]$, where $\wtilde{\lambda} = \theta + \theta^{-1} + \cpxi n^{-1/2}$ is a complex constant.
Let us note in passing that this constant falls in the regions treated in the local laws (Theorem~\ref{thm:ill} and Corollary~\ref{cor:ill-Q}) for suitable choices of the parameters there.
All resolvents are evaluated at this $\wtilde{\lambda}$ in the following calculations, so for $W^{(a)}$ and all other matrices involved, let us abbreviate
\[ R_{W^{(a)}}\colonequals R_{W^{(a)}}(\wtilde{\lambda}). \]

We first compute the effect on the resolvent of the above perturbations of $\bQ$.
By the resolvent expansion of Proposition~\ref{prop:res-id}, for each $a\in[L]$, we have
\begin{equation}
R_{W^{(\ell - 1,a)}} = R_{Q^{(a)} + E^{(W,a)}_{pq}} = \sum_{s = 0}^4 R_{Q^{(a)}}\left(E^{(W,a)}_{pq} R_{Q^{(a)}}\right)^s + R_{W^{(\ell - 1,a)}}\left(E^{(W,a)}_{pq} R_{Q^{(a)}}\right)^{5}.
\label{eq:res-exp}
\end{equation}
We will see momentarily why taking the expansion to order 4 is the correct choice here.

To work with these expressions, let us set up some notation for what happens when we expand the $E^{(W,a)}_{pq}$.
\begin{definition}[Conjugation words]\label{def:conjugation-words}
Write $\id: \CC \to \CC$ for the identity map $\id(z) = z$, and $\conj: \CC \to \CC$ for the conjugation map $\conj(z) = \overline{z}$.
For a given $1 \leq p \leq q \leq n$, define
\[ \mathcal{S}_{p,q,s} = \left\{\begin{array}{ll} \{\id, \conj\}^s & \text{if } p \neq q, \\ \{\id\}^s & \text{if } p = q\end{array}\right\}. \]
For $S = (f_1, \dots, f_s) \in \mathcal{S}_{p,q,s}$, viewed as a string of functions of length $s$, let
\[ z^S \colonequals f_1(z) \cdots f_s(z). \]
Also, associate to such $S$ indices $\alpha_{S, 1}, \beta_{S, 1}, \dots, \alpha_{S, s}, \beta_{S, s} \in \{p, q\}$, so that $\alpha_{S, i} = p$ and $\beta_{S, i} = q$ if $f_i = \id$, while $\alpha_{S, i} = q$ and $\beta_{S, i} = p$ if $f_i = \conj$.
\end{definition}

Then, for $a\in[L]$ and $s\ge 1$, we may expand
\[ \left(E^{(W,a)}_{pq} R_{Q^{(a)}}\right)^s = \sum_{S \in \mathcal{S}_{p,q,s}} (W^{(a)}_{pq})^S e_{\alpha_{S,1}} \left(e_{\beta_{S,1}}^* R_{Q^{(a)}} e_{\alpha_{S,2}}\right) \cdots \left(e_{\beta_{S,{s-1}}}^* R_{Q^{(a)}} e_{\alpha_{S,s}}\right) e_{\beta_{S,s}}^*R_{Q^{(a)}}, \]
where we note that the quantities in parentheses are just scalars giving certain entries of $R_{Q^{(a)}}$.
Since such expressions will often come up, let us write $e_p^* R_Y e_q = (R_Y)_{pq} \equalscolon R_Y(p,q)$ for matrices $Y \in \{W^{(\ell - 1,a)}, W^{(\ell,a)}, Q^{(a)}\}$.
Plugging this into \eqref{eq:res-exp}, we have
\begin{align*}
    R_{W^{(\ell - 1,a)}}
    &= \sum_{s = 0}^4 \sum_{S \in \mathcal{S}_{p,q,s}} (W^{(a)}_{pq})^S \cdot R_{Q^{(a)}}(\beta_{S,1}, \alpha_{S,2}) \cdots R_{Q^{(a)}}(\beta_{S,{s-1}}, \alpha_{S,s}) \cdot R_{Q^{(a)}} e_{\alpha_{S,1}}e_{\beta_{S,s}}^*R_{Q^{(a)}} \\
    &\hspace{0.5cm} + \sum_{S \in \mathcal{S}_{p,q,5}} (W^{(a)}_{pq})^S \cdot R_{Q^{(a)}}(\beta_{S,1}, \alpha_{S,2}) \cdots R_{Q^{(a)}}(\beta_{S,t}, \alpha_{S,t+1}) \cdot R_{W^{(\ell - 1,a)}} e_{\alpha_{S,1}}e_{\beta_{S,t+1}}^*R_{Q^{(a)}}.
\end{align*}
The expressions we will finally be interested in are $e_i^* R_Y v^{(a)}$ and $v^{(a)*} R_Y e_j$.
Let us extend the previous notation and write $R_Y(i, v^{(a)})$ and $R_Y(v^{(a)}, j)$ for these, respectively, and extending in the same way to other mixed quadratic forms of this kind, for all $Y$ mentioned above.

Then, firstly we may rewrite
\begin{equation}
\label{eq:F-W-ell}
\wtilde{F}_{ij}(\bW^{(\ell - 1)}) = \wtilde{\phi} \left( \left(n \cdot R_{W^{(\ell - 1,a)}}(i, v^{(a)}) R_{W^{(\ell - 1,a)}}(v^{(a)}, j)\right)_{a=1}^{L} \right),
\end{equation}
and similarly for $\ell - 1$ replaced by $\ell$, and for the above expansion of the inner terms, we have for instance
\begin{align*}
    R_{W^{(\ell - 1,a)}}(i, v^{(a)})
    &= M^{(\ell,a)}(i, v^{(a)}, \bW) + \Delta^{(\ell,a)}(i, v^{(a)}, \bW), \\
    M^{(\ell,a)}(i,v^{(a)}, \bY)
    &= \sum_{s = 0}^4 \sum_{S \in \mathcal{S}_{p,q,s}} (Y^{(a)}_{pq})^S  R_{Q^{(a)}}(i, \alpha_{S,1}) R_{Q^{(a)}}(\beta_{S,1}, \alpha_{S,2}) \cdots \\
    &\hspace{4cm} R_{Q^{(a)}}(\beta_{S,{s-1}}, \alpha_{S,s}) R_{Q^{(a)}}(\beta_{S, s}, v^{(a)}) \\
    &= R_{Q^{(a)}}(i, v^{(a)}) + \sum_{s = 1}^4 \sum_{S \in \mathcal{S}_{p,q,s}} (Y^{(a)}_{pq})^S  R_{Q^{(a)}}(i, \alpha_{S,1}) R_{Q^{(a)}}(\beta_{S,1}, \alpha_{S,2}) \cdots \\
    &\hspace{6cm} R_{Q^{(a)}}(\beta_{S,{s-1}}, \alpha_{S,s}) R_{Q^{(a)}}(\beta_{S, s}, v^{(a)}), \\
    \Delta^{(\ell,a)}(i, v^{(a)}, \bY) &= \sum_{S \in \mathcal{S}_{p,q,5}} (Y^{(a)}_{pq})^S  R_{Q^{(a)} + E^{(Y,a)}_{pq}}(i, \alpha_{S,1}) R_{Q^{(a)}}(\beta_{S,1}, \alpha_{S,2}) \cdots \\
    &\hspace{4cm} R_{Q^{(a)}}(\beta_{S,4}, \alpha_{S,5}) R_{Q^{(a)}}(\beta_{S, 5}, v^{(a)}).
\end{align*}
We define these expressions taking $p = p_{\ell}$ and $q = q_{\ell}$ under our previous ordering, and allow for a general matrix input $Y$, so that we also have
\[ R_{W^{(\ell,a)}}(i, v^{(a)}) = M^{(\ell,a)}(i, v^{(a)}, \bX) + \Delta^{(\ell,a)}(i, v^{(a)}, \bX), \]
reusing the same definitions.
Also, as a convention, we take the term $s = 0$ in this expression to contribute $R_{Q^{(a)}}(i, v^{(a)})$, compatible with our previous calculation.

We will now show the following estimates, which imply that the various $\Delta$ terms in the input to $\wtilde{\phi}$ in \eqref{eq:F-W-ell} may be ignored.
\begin{lem}
    \label{lem:tildeF-poly}
    In the above setting, for each $\ell \in [N]$, we have
    \begin{align}
    &\left|\EE \wtilde{F}_{ij}(\bW^{(\ell - 1)}) - \EE \wtilde{\phi} \left(\left(n \cdot M^{(\ell,a)}(i, v^{(a)}, \bW) M^{(\ell,a)}(v^{(a)}, j, \bW)\right)_{a=1}^{L} \right)\right| \lesssim n^{-5/2 + 3\epsilon_v}, \label{eq:tildeF-poly-Wellm1} \\
    &\left|\EE \wtilde{F}_{ij}(\bW^{(\ell)}) - \EE \wtilde{\phi} \left(\left(n \cdot M^{(\ell,a)}(i, v^{(a)}, \bX) M^{(\ell,a)}(v^{(a)}, j, \bX)\right)_{a=1}^{L}\right)\right| \lesssim n^{-5/2 + 3\epsilon_v}, \label{eq:tildeF-poly-Well}
    \end{align}
    and thus by \eqref{eq:lindeberg-triangle}
    \begin{align*}
    &|\EE \wtilde{F}_{ij}(\bw) - \EE \wtilde{F}_{ij}(\bx)| \\
    &\hspace{0.5cm} \lesssim \sum_{\ell = 1}^{N} \bigg|\EE \wtilde{\phi}\left( \left(n \cdot M^{(\ell,a)}(i, v^{(a)}, \bW) M^{(\ell,a)}(v^{(a)}, j,\bW)\right)_{a=1}^{L}\right) \notag\\ &\hspace{3.4cm} - \EE \wtilde{\phi}\left(\left(n \cdot M^{(\ell,a)}(i, v^{(a)}, \bX) M^{(\ell,a)}(v^{(a)}, j, \bX)\right)_{a=1}^{L}\right)\bigg| + O(n^{-1/2 + 3\epsilon_v}). \numberthis \label{eq:tildeF-poly-diff-bd}
\end{align*}
\end{lem}
\noindent
Having shown this,  we will have again made useful progress. For the $\ell$-th term, the inputs into $\wtilde{\phi}$ are functions of the current entry vector $w^{(pq)}$, $x^{(pq)}$. More precisely, for each coordinate $a\in[L]$, the quantity $M^{(\ell,a)}(i,v^{(a)},\bW)M^{(\ell,a)}(v^{(a)},j,\bW)$ is a polynomial in $W^{(a)}_{pq}$, $\overline{W^{(a)}_{pq}}$, and the entries of $R_{Q^{(a)}}$. Crucially, $\bQ$ is independent of those quantities by definition.

To do this, it will be useful to establish some estimates on the sizes of resolvent quadratic forms and associated expectations.
\begin{prop}[Resolvent quadratic form bounds]
    \label{prop:resolvent-entry-bounds}
    In the above setting, for any $a\in[L]$, any distinct $\alpha, \beta \in [n]$, and any $Y \in \{W^{(\ell,a)}, W^{(\ell - 1,a)}, Q^{(a)}\}$, we have
    \begin{align*}
    |Y(\alpha, \alpha)| &\prec n^{-1/2}, \\
    |Y(\alpha, \beta)| &\prec n^{-1/2}, \\
    |R_Y(\alpha, \beta)|
    &\prec n^{-1/2}, \\
    |R_Y(\alpha, \alpha)|
    &\prec 1, \\
    |R_Y(\alpha, v^{(a)})|
    &\prec n^{-1/2 + \epsilon_v}, \\
    |R_Y(v^{(a)}, v^{(a)})|
    &\prec 1.
    \end{align*}
    Further, fix nonnegative integers $s,t$, and consider a product of terms
    \[ \Pi = R_{Y_1}(\zeta_{11}, \zeta_{12}) \cdots R_{Y_s}(\zeta_{s1}, \zeta_{s2}) \cdot Y_{s + 1}(\alpha_{11}, \alpha_{12}) \cdots Y_{s + t}(\alpha_{t1}, \alpha_{t2}). \]
    For each $m \in [s+t]$, suppose that there is an index $a_m\in [L]$ such that
    \[ Y_m \in \{W^{(\ell, a_m)}, W^{(\ell - 1, a_m)}, Q^{(a_m)}\}. \]
    For $m\in[s]$, each argument $\zeta_{m1}$, $\zeta_{m2}$ is either an index in $[n]$, which represents the corresponding basis vector, or the formal symbol $v^{(a_m)}$.
    Thus, for example
    \[R_{Y_{m}}(\alpha, v^{(a_m)}) = e^{*}_{\alpha}R_{Y_m}v^{(a_m)},\qquad R_{Y_m}(v^{(a_m)},\beta)=v^{(a_m)*} R_{Y_m}e_{\beta}.\] For the entry factors, we have $\alpha_{m1},\alpha_{m2} \in[n]$, $m\in[t]$, so that $Y_{s+m}(\alpha_{m1},\alpha_{m2})$ denotes the $(\alpha_{m1},\alpha_{m2})$-th entry of the matrix $Y_{s+m}$.

    Let $s_v$ be the number of resolvent factors whose two arguments contain the corresponding formal symbol $v^{(a_m)}$ exactly once, and let $s_{\offdiag}$ be the number of resolvent factors whose two arguments are distinct as formal symbols.
    Then, we have
    \begin{equation}
        \label{eq:resolvent-prod-entries-dom}
        |\Pi| \prec n^{-t/2 - s_{\offdiag}/2 + s_v \epsilon_v}
    \end{equation}
    and, for any $\epsilon > 0$,
    \begin{equation}
    \label{eq:resolvent-prod-entries-exp}
    \EE |\Pi| \lesssim_{\epsilon} n^{-t/2 - s_{\offdiag}/2 + s_v \epsilon_v + \epsilon}.
    \end{equation}
\end{prop}
\begin{proof}
    We first prove the individual stochastic domination bounds. The entry bounds follow immediately from the $\xi$-sub-Weibull tail bound in \eqref{eq:entries-subexp}, since every entry of any $Y$ is either zero, an entry of $W^{(a)}$, or an entry of $X^{(a)}$.

    Next, for fixed $a\in[L]$, the matrices $W^{(\ell,a)}$ and $W^{(\ell - 1,a)}$ are generalized Wigner matrices with uniformly controlled parameters. Indeed, the second-moment matching of the vectors $w^{(pq)}$ and $x^{(pq)}$ implies in particular that $W^{(a)}$ and $X^{(a)}$ have the same variance profile. Moreover, $Q^{(a)}$ is obtained from such a matrix by setting one upper-triangular entry and its Hermitian conjugate to zero. Hence the local laws of Theorem~\ref{thm:ill} apply to $W^{(\ell,a)}$ and $W^{(\ell - 1,a)}$, while Corollary~\ref{cor:ill-Q} applies to $Q^{(a)}$, since applying the triangle inequality to those results it suffices to observe that $|G_{\mu_{\sc}}(\wtilde{\lambda})|$ is just a constant, while $|v_{\alpha}^{(a)}| \leq \|v^{(a)}\|_{\infty} \leq n^{-1/2 + \epsilon_v}$ by assumption. Since $L$ is fixed, all estimates are uniform in $a\in[L]$.

    The bound of \eqref{eq:resolvent-prod-entries-dom} then follows by Proposition~\ref{prop:stoch-dom-sum-prod}, since this is just a product of several polynomial stochastic dominations from the first set of results. It remains to pass to expectations. If $s+t=0$, then $\Pi=1$, and there is nothing to show. Otherwise, \eqref{eq:resolvent-prod-entries-exp} says that we may take expectations on either side of this polynomial stochastic domination, provided we insert a factor of $n^{\epsilon}$ on the right-hand side.
    By Proposition~\ref{prop:stochastic-domination-to-moment-bound}, to show this it suffices to check that $\EE|\Pi|^2$ is bounded independently of $n$.
    But, we have
    \[ |\Pi| \leq \|R_{Y_1}\| \cdots \|R_{Y_s}\| \cdot |Y_{s + 1}(\alpha_{11}, \alpha_{12})| \cdots |Y_{s + t}(\alpha_{t1}, \alpha_{t2})|, \]
    and so by H\"{o}lder's inequality,
    \begin{align*}
    \EE|\Pi|^2
    &\leq \left(\EE \|R_{Y_1}\|^{2s + 2t}\right)^{\frac{1}{s + t}} \cdots \left(\EE \|R_{Y_s}\|^{2s + 2t}\right)^{\frac{1}{s + t}} \\
    &\hspace{1cm} \left(\EE |Y_{s + 1}(\alpha_{11}, \alpha_{12})|^{2s + 2t}\right)^{\frac{1}{s + t}} \cdots \left(\EE |Y_{s + t}(\alpha_{t1}, \alpha_{t2})|^{2s + 2t}\right)^{\frac{1}{s + t}}
    \end{align*}
    and we then indeed find that $\EE|\Pi|^2 = O(1)$ by Proposition~\ref{prop:res-moments} and Proposition~\ref{prop:subexp-moments} which show that each factor here is $O(1)$.

    We note that a condition of that Proposition is that $\Im(z) > n^{-\Omega(1)}$, so we are using that we are evaluating our resolvents at $\wtilde{\lambda}$ rather than $\lambda \in \RR$, and this cannot be avoided for working with such expectations (without introducing other technical devices like restricting to particular events) for the reason discussed in Remark~\ref{rem:res-cpx}.
\end{proof}

\begin{proof}[Proof of Lemma~\ref{lem:tildeF-poly}]
    It suffices to show \eqref{eq:tildeF-poly-Wellm1}, then \eqref{eq:tildeF-poly-Well} follows by a symmetric argument, and \eqref{eq:tildeF-poly-diff-bd} follows from the two taken together and the triangle inequality since $N = O(n^2)$.

    By Proposition~\ref{prop:resolvent-entry-bounds}, uniformly in $a\in[L]$, we have that
    \begin{align*}
M^{(\ell,a)}(i, v^{(a)}, \bW)
&\prec n^{-1/2 + \epsilon_v}, \\
\intertext{with the dominant contribution coming from the $s = 0$ term, and}
\Delta^{(\ell,a)}(i, v^{(a)}, \bW)
&\prec n^{-3 + \epsilon_v},
\end{align*}
where in the $\Delta^{(\ell,a)}$ bounds a factor of $n^{-5/2}$ comes from 5 factors of entries of $Y$, and another factor of $n^{-1/2 + \epsilon_v}$ comes from one factor of the form $R_Y(\alpha, v)$ for some $\alpha \in [n]$.
    (We use here that the $M^{(\ell,a)}$ and $\Delta^{(\ell,a)}$ are finite sums of the form treated by Proposition~\ref{prop:resolvent-entry-bounds}, which may be combined over finite sums by Proposition~\ref{prop:stoch-dom-sum-prod}.)
    The same bounds hold for $(i, v^{(a)})$ replaced by $(v^{(a)}, j)$ as well.
    Then, we may expand
    \begin{align*}
    \wtilde{F}_{ij}(\bW^{(\ell - 1)})
    &= \wtilde{\phi}\left( \left(n \cdot R_{W^{(\ell - 1,a)}}(i, v^{(a)}) R_{W^{(\ell - 1,a)}}(v^{(a)}, j)\right)_{a=1}^{L} \right) \\
    &= \wtilde{\phi}\bigg(\Big(n \cdot M^{(\ell,a)}(i, v^{(a)}, \bW)M^{(\ell,a)}(v^{(a)}, j, \bW) \\
    &\hspace{1cm} + n \cdot M^{(\ell,a)}(i, v^{(a)}, \bW)\Delta^{(\ell,a)}(v^{(a)}, j, \bW) \\
    &\hspace{1cm} + n \cdot \Delta^{(\ell,a)}(i, v^{(a)}, \bW)M^{(\ell,a)}(v^{(a)}, j, \bW) \\
    &\hspace{1cm} + n \cdot \Delta^{(\ell,a)}(i, v^{(a)}, \bW)\Delta^{(\ell,a)}(v^{(a)}, j, \bW)\Big)_{a=1}^{L} \bigg) \\
    &\equalscolon \wtilde{\phi}\left(\left(n \cdot M^{(\ell,a)}(i, v^{(a)}, \bW)M^{(\ell,a)}(v^{(a)}, j, \bW) + \Xi_{a}\right)_{a=1}^{L}\right),
    \end{align*}
    where the above bounds imply that
    \[ |\Xi_{a}| \prec n^{-5/2 + 2\epsilon_v}. \]
    Using that $\wtilde{\phi}$ is Lipschitz and $L$ is fixed, we then have
    \[ \left|\wtilde{F}_{ij}(\bW^{(\ell - 1)}) - \wtilde{\phi}\left( \left(n \cdot M^{(\ell,a)}(i, v^{(a)}, \bW)M^{(\ell,a)}(v^{(a)}, j, \bW) \right)_{a=1}^{L} \right) \right| \leq \sum_{a=1}^{L}|\Xi_a|, \]
    and the result follows upon taking expectations and using the triangle inequality, since $\EE |\Xi_a|$ may be bounded by another triangle inequality as a sum of terms each controlled by Proposition~\ref{prop:resolvent-entry-bounds}.
    Applying the Proposition costs another factor of $n^{\epsilon}$ for an arbitrarily small $\epsilon > 0$, and we take $\epsilon = \epsilon_v$ to obtain the result as stated.
\end{proof}

\begin{rmk}[Order of resolvent expansion]
    \label{rmk:res-exp-order}
    We see in the above bounds why our choice of fourth order in the resolvent expansion of \eqref{eq:res-exp} was convenient: in order for $N=O(n^2)$ terms with the above error bound to not contribute in total to our bound on $|\EE \wtilde{F}_{ij}(\bw) - \EE \wtilde{F}_{ij}(\bx)|$, we must have each term to be bounded by $n^{-2-\delta}$ for some $\delta > 0$, which will no longer hold by the above argument if we take a third order (or shorter) expansion.\footnote{It seems that with slightly more care in the combinatorics of how many terms of what ``types'' in terms of the intersection among the $\{p_{\ell}, q_{\ell}\}$ and $\{i, j\}$ index pairs one could take an expansion to order 3, but we take the longer expansion that is more transparently correct for the sake of exposition.}
\end{rmk}

We now continue to finish the proof of the main result of this section, modulo one technical result we will encounter at the end of our calculation that we defer to the following section.

To keep track of conjugated scalar resolvent forms arising from the multivariate Taylor expansion below, let $\sigma\in\{\id,\conj\}$ with the maps defined in Definition~\ref{def:conjugation-words} and write
\[R^{\sigma}_{Y}(\zeta_1, \zeta_2)\colonequals \sigma\big(R_{Y}(\zeta_1,\zeta_2)\big).\]
Since $|R^{\sigma}_{Y}(\zeta_1, \zeta_2)| = |R_{Y}(\zeta_1, \zeta_2)|$, the same bounds in Proposition~\ref{prop:resolvent-entry-bounds} hold if any scalar resolvent form $R_{Y}(\zeta_1, \zeta_2)$ in the statement is replaced by either $R^{\sigma}_{Y}(\zeta_1, \zeta_2)$.

\begin{proof}[Proof of Lemma~\ref{lem:tilde-F-bound}]
Starting from the result of Lemma~\ref{lem:tildeF-poly}, we must control summands of the form
\begin{align*}
\bigg| &\EE \wtilde{\phi}\left( \left(n\cdot M^{(\ell,a)}(i,v^{(a)},\bW)
M^{(\ell,a)}(v^{(a)},j,\bW)
\right)_{a=1}^L \right) \\
&\qquad - \EE \wtilde{\phi}\left(
\left(n\cdot M^{(\ell,a)}(i,v^{(a)},\bX)
M^{(\ell,a)}(v^{(a)},j,\bX)
\right)_{a=1}^L\right) \bigg|.
\end{align*}
For a word $\bkappa\colonequals(\kappa_{1},\ldots,\kappa_{u})\in[2L]^{u}$, write
\[\partial_{\bkappa}\colonequals \partial_{\kappa_1}\ldots\partial_{\kappa_u},\]and define
\[\Phi_{\bkappa, i, j}(\bQ) \colonequals \frac{1}{u!} \partial_{\bkappa} \wtilde{\phi} \left( \left(n \cdot R_{Q^{(a)}}(i, v^{(a)})R_{Q^{(a)}}(v^{(a)}, j)\right)_{a=1}^{L}\right).\]
We now take a multivariate Taylor expansion of $\wtilde{\phi}$ in each term in such a difference.
We define, for each $a\in[L]$,
\[
\Gamma_{a}(\bY) \colonequals M^{(\ell,a)}(i, v^{(a)}, \bY) M^{(\ell,a)}(v^{(a)}, j, \bY) - R_{Q^{(a)}}(i, v^{(a)})R_{Q^{(a)}}(v^{(a)}, j).\]
We view $\wtilde{\phi}$ as a function on $\RR^{2L}$ and define
\[\bGamma(\bY)\colonequals (\Re\Gamma_{1}(\bY),\Im\Gamma_{1}(\bY),\ldots, \Re\Gamma_{L}(\bY),\Im\Gamma_{L}(\bY)),\]
so \[\bGamma_{2a-1}(\bY) = \Re\Gamma_{a}(\bY) = \frac{\id(\Gamma_{a}(\bY)) + \conj(\Gamma_{a}(\bY))}{2}, \]
\[\bGamma_{2a}(\bY) = \Im\Gamma_{a}(\bY) = \frac{\id(\Gamma_{a}(\bY)) - \conj(\Gamma_{a}(\bY))}{2\cpxi}.\]
By the same application of Proposition~\ref{prop:resolvent-entry-bounds} as above in the proof of Lemma~\ref{lem:tildeF-poly}, we have, uniformly in $a\in[L]$,
\[ |\Gamma_{a}(\bY)| \prec n^{-3/2 + 2\epsilon_v} \]
for each $\bY \in \{\bW, \bX\}$.
Taking a Taylor expansion to fourth order around
\[\left( n \cdot R_{Q^{(a)}}(i, v^{(a)}) R_{Q^{(a)}}(v^{(a)}, j)\right)_{a=1}^{L}\] in each expectation, we find
\begin{align*}
&\hspace{-1cm}\wtilde{\phi} \left( \left(n \cdot M^{(\ell,a)}(i, v^{(a)}, \bY) M^{(\ell,a)}(v^{(a)}, j, \bY) \right)_{a=1}^{L} \right)\\
&= \sum_{u = 0}^4  \sum_{\bkappa\in[2L]^{u}}  \Phi_{\bkappa, i, j}(\bQ) \cdot n^u \cdot \prod_{r=1}^{u}\bGamma_{\kappa_r}(\bY) + \Xi^{\prime}\\
&= \sum_{u = 0}^4  \sum_{\bkappa\in[2L]^{u}} \Phi_{\bkappa, i, j}(\bQ) \cdot n^u \cdot \prod_{r=1}^u \left( \sum_{\sigma_r\in\{\id,\conj\}} c_{\sigma_r,\kappa_r}\cdot \sigma_r\big(\Gamma_{\lceil\frac{\kappa_r}{2}\rceil}(\bY)\big) \right) + \Xi^{\prime}\\
&\equalscolon \sum_{u = 0}^4  \sum_{\bkappa\in[2L]^{u}} \Phi_{\bkappa, i, j}(\bQ) \cdot n^u \cdot \sum_{\bsigma\in\{\id,\conj\}^u} c_{\bkappa,\bsigma} \prod_{r=1}^u \sigma_r\big(\Gamma_{\lceil\frac{\kappa_r}{2}\rceil}(\bY)\big) + \Xi^{\prime}
\end{align*}
where $|c_{\bkappa,\bsigma}| = |\prod_{r=1}^u c_{\sigma_r,\kappa_r}| = 2^{-u} \leq 1$ and $\Xi^{\prime}$ is a random error term (the prime mark distinguishing it from the one in the proof of Lemma~\ref{lem:tildeF-poly} above) that is bounded by
\[ |\Xi^{\prime}| \lesssim n^5 \cdot \left(\max_{a\in[L]}|\Gamma_a(\bY)| \right)^5 \prec n^{-5/2 + 10\epsilon_v}. \]
Thus, again using the expectation bounds in Proposition~\ref{prop:resolvent-entry-bounds}, we have
\begin{align*}
&\left| \EE \wtilde{\phi} \left( \left(n \cdot M^{(\ell,a)}(i, v^{(a)}, \bY) M^{(\ell,a)}(v^{(a)}, j, \bY) \right)_{a=1}^{L} \right) - \EE \sum_{u = 0}^4  \sum_{\bkappa\in[2L]^{u}}\Phi_{\bkappa, i, j}(\bQ) \cdot n^u \cdot \prod_{r=1}^{u}\bGamma_{\kappa_r}(\bY)\right| \\
&\hspace{1cm} \leq \EE |\Xi^{\prime}| \\ &\hspace{1cm}\lesssim_{\epsilon} n^{-5/2 + 10\epsilon_v + \epsilon}
\end{align*}
for $\epsilon > 0$ arbitrarily small.

Combining this with Lemma~\ref{lem:tildeF-poly}, we find
\begin{align*}
    &|\EE \wtilde{F}_{ij}(\bw) - \EE \wtilde{F}_{ij}(\bx)| \\
    &\hspace{0.5cm} \lesssim \sum_{\ell = 1}^{N} \left|\EE \sum_{u = 0}^4  \sum_{\bkappa\in[2L]^u}\Phi_{\bkappa, i, j}(\bQ) \cdot n^u \cdot \left( \prod_{r=1}^{u} \bGamma_{\kappa_r}(\bW) - \prod_{r=1}^{u}\bGamma_{\kappa_r}(\bX) \right)\right| + O(n^{-1/2 + 10\epsilon_v + \epsilon}).
\end{align*}
(In the same way as discussed in Remark~\ref{rmk:res-exp-order}, we see from this calculation that fourth order was the correct order of Taylor expansion for this argument to work.)
Let us write
\[ \epsilon^{\prime} \colonequals 10\epsilon_v + \epsilon \]
for this remaining error exponent, which we will carry through to the end of the proof and which gives rise to the dependence on $\epsilon_v$ in Theorem~\ref{thm:main-univ}.

The only properties of $\Phi_{\bkappa, i, j}(\bQ)$ we will need are that it is $O(1)$ almost surely (by the boundedness of $\wtilde{\phi}$) and its derivatives and that it is independent of $w^{(pq)}$ and $x^{(pq)}$ (since it is a function only of $\bQ$, where this entry is set to zero).
Also, let us develop notation for the quantities appearing in powers of $\bGamma(\bY)$.
For $a\in[L]$ and $\sigma\in\{\id,\conj\}$, expanding $\sigma(\Gamma_a(\bY))$ from the definition of $\Gamma_a$, we have
\begin{align*}
    \sigma(\Gamma_{a}(\bY)) &= \sum_{\substack{(s_1, s_2) \in \{0, 1, 2, 3, 4\}^2 \\ (s_1, s_2) \neq (0, 0)}} \sum_{\substack{S_1 \in \mathcal{S}_{p,q,s_1} \\ S_2 \in \mathcal{S}_{p,q,s_2}}} \sigma\left((Y^{(a)}_{pq})^{S_1}(Y^{(a)}_{pq})^{S_2}\right) \\ &\hspace{2cm} \cdot R^{\sigma}_{Q^{(a)}}(i, \alpha_{S_1,1}) R^{\sigma}_{Q^{(a)}}(\beta_{S_1,1}, \alpha_{S_1,2}) \cdots R^{\sigma}_{Q^{(a)}}(\beta_{S_1,{s_1-1}}, \alpha_{S_1,s_1}) R^{\sigma}_{Q^{(a)}}(\beta_{S_1, s_1}, v^{(a)}) \\ &\hspace{2cm} \cdot R^{\sigma}_{Q^{(a)}}(v^{(a)}, \alpha_{S_2,1}) R^{\sigma}_{Q^{(a)}}(\beta_{S_2,1}, \alpha_{S_2,2}) \cdots R^{\sigma}_{Q^{(a)}}(\beta_{S_2,{s_2-1}}, \alpha_{S_2,s_2}) R^{\sigma}_{Q^{(a)}}(\beta_{S_2, s_2}, j) \\  &\equalscolon \sum_{\substack{(s_1, s_2) \in \{0, 1, 2, 3, 4\}^2 \\ (s_1, s_2) \neq (0, 0)}} \sum_{\substack{S_1 \in \mathcal{S}_{p,q,s_1} \\ S_2 \in \mathcal{S}_{p,q,s_2}}} (Y^{(a)}_{pq})^{S_1}(Y^{(a)}_{pq})^{S_2} J^{(a,\sigma)}_{i,j,S_1}(Q^{(a)}) J^{(a,\sigma)}_{i,j,S_2}(Q^{(a)}),
    \intertext{for $\sigma \in\{\id,\conj\}$. It will be useful to rearrange a little bit more, by defining the union $\mathcal{S}_{p,q} \colonequals \mathcal{S}_{p,q,0} \sqcup \cdots \sqcup \mathcal{S}_{p,q,4}$. For $S \in \mathcal{S}_{p,q}$, is a string of indeterminate length now, write $s(S)$ for its length, so that $S \in \mathcal{S}_{p,q,s(S)}$. Then, writing $\emptyset$ for the empty string that is the only element of $\mathcal{S}_{p,q,0}$, we may write the above as a single sum,}
    &= \sum_{(S_1, S_2) \in \mathcal{S}_{p,q}^2 \setminus \{(\emptyset, \emptyset)\}} \left((Y^{(a)}_{pq})^{S_1}(Y^{(a)}_{pq})^{S_2}\right)^\sigma J^{(a,\sigma)}_{i,j,S_1}(Q^{(a)}) J^{(a,\sigma)}_{i,j,S_2}(Q^{(a)}).
\end{align*}
We may then write the result of taking products of the real-coordinate components of $\bGamma(\bY)$ in terms of matrices with entries taking values in $\mathcal S_{p,q}$. For each $r\in[u]$, let $a_r \in [L]$ be the unique index such that $\kappa_r \in \{2a_r-1, 2a_r\}$.
Expanding the real and imaginary parts of $\Gamma_{a_r}(\bY)$ gives a finite linear combination of terms of the following form:
\[ \prod_{r=1}^{u}\bGamma_{\kappa_r}(\bY) = \sum_{\bsigma\in\{\id,\conj\}^u}c_{\bkappa,\bsigma} \sum_{\substack{S \in \mathcal{S}_{p,q}^{u \times 2} \\ \text{no row of } S \text{ equals } (\emptyset, \emptyset)}} \bY^{\bkappa,S,\bsigma}_{pq} J^{\bsigma}_{\bkappa,i,j,S}(\bQ),\]

\noindent
Lastly, for such a matrix $S$, define
\begin{align*}
    \bY^{\bkappa,S,\bsigma}_{pq} &\colonequals \prod_{r=1}^u \sigma_r\left((Y^{(a_r)}_{pq})^{S_{r1}} (Y^{(a_r)}_{pq})^{S_{r2}}\right), \\
    J^{\bsigma}_{\bkappa,i,j,S}(\bQ) &\colonequals \prod_{r=1}^u J^{(a_r,\sigma_r)}_{i,j,S_{r1}}(Q^{(a_r)}) J^{(a_r,\sigma_r)}_{i,j,S_{r2}}(Q^{(a_r)}), \\
    |S| &\colonequals \sum_{r, f} s(S_{rf}), \\
    |S|_0 &\colonequals \#\{(r, f): S_{rf} \neq \emptyset\}.
\end{align*}
Using this, we may rewrite our bound concisely as
\begin{align*}
    &|\EE \wtilde{F}_{ij}(\bw) - \EE \wtilde{F}_{ij}(\bx)| \\
    & \lesssim \sum_{\ell = 1}^{N} \bigg|\EE \sum_{u = 0}^4 \sum_{\bkappa\in[2L]^u} \sum_{\bsigma\in\{\id,\conj\}^u} \\
    &\hspace{1cm} \sum_{\substack{S \in \mathcal{S}_{p,q}^{u \times 2} \\ \text{no row of } S \text{ equals } (\emptyset, \emptyset)}} \Phi_{\bkappa, i, j}(\bQ)J^{\bsigma}_{\bkappa,i,j,S}(\bQ) \cdot n^u \cdot \left(\bW^{\bkappa,S,\bsigma}_{pq} - \bX^{\bkappa,S,\bsigma}_{pq}\right)\bigg| + O(n^{-1/2 + \epsilon^{\prime}})
    \intertext{where we may use independence of $\bQ$ and $(W^{(a)}_{pq}, X^{(a)}_{pq})$, giving when combined with a triangle inequality}
    &\leq \sum_{\ell = 1}^N \sum_{u = 0}^4 \sum_{\bkappa\in[2L]^u} \sum_{\bsigma\in\{\id,\conj\}^u} \\
    &\hspace{1cm} \sum_{\substack{S \in \mathcal{S}_{p,q}^{u \times 2} \\ \text{no row of } S \text{ equals } (\emptyset, \emptyset)}} \left|\EE[\Phi_{\bkappa, i, j}(\bQ)J^{\bsigma}_{\bkappa,i,j,S}(\bQ)]\right| \cdot n^u \cdot \left|\EE\big[\bW^{\bkappa,S,\bsigma}_{pq}\big] - \EE\big[\bX^{\bkappa,S,\bsigma}_{pq}\big]\right| + O(n^{-1/2 + \epsilon^{\prime}})
    \intertext{Here, for each fixed $\bsigma$, the expressions $\bW_{pq}^{\bkappa,S,\bsigma}$ and $\bX_{pq}^{\bkappa,S,\bsigma}$ are fixed complex linear combinations of monomials of total degree $|S|$ in the real and imaginary parts of the block variables $w^{(pq)},x^{(pq)}\in\RR^{2L}$. Indeed, applying either $\id$ or $\conj$ only decides whether we use an entry monomial or its complex conjugate, and this does not change the total degree. Therefore the difference
    $\EE[\bW_{pq}^{\bkappa,S,\bsigma}] - \EE[\bX_{pq}^{\bkappa,S,\bsigma}]$
    vanishes whenever $|S|\leq 2$, by the assumed matching of joint moments up to order two.
    On the other hand, we have $|J^{\bsigma}_{\bkappa,i,j,S}(\bQ)| \prec n^{-u + 2u\epsilon_v} \prec n^{-u + 8\epsilon_v}$ for all $S \in \mathcal{S}^{u\times 2}_{p,q}$ by Proposition~\ref{prop:resolvent-entry-bounds}.
    Thus, also applying the expectation bounds from Proposition~\ref{prop:resolvent-entry-bounds} and using the boundedness of $\Phi_{\bkappa,i,j}(\bQ)$, each term above has a simple \emph{a priori} bound of $O(n^u\cdot n^{-u+8\epsilon_v+\epsilon}\cdot n^{-|S|/2}) = O(n^{-|S| / 2 + 9\epsilon_v})$. (Here and below we let $\epsilon \in (0, \epsilon_v)$ be a temporary parameter as needed.) In particular, if $|S| \geq 5$ then this is at most $O(n^{-5/2 + 9\epsilon_v})$, and since $9\epsilon_v < \epsilon^{\prime}$ this may be subsumed into the error term we already have, giving}
    &\leq \sum_{\ell = 1}^N \sum_{u = 0}^4 \sum_{\bkappa\in[2L]^u} \sum_{\bsigma\in\{\id,\conj\}^u} \\ &\hspace{1cm} \sum_{\substack{S \in \mathcal{S}_{p,q}^{u \times 2} \\ \text{no row of } S \text{ equals } (\emptyset, \emptyset) \\ |S| \in \{3, 4\}}} n^u \cdot |\EE[\Phi_{\bkappa,i,j}(\bQ)
    J^{\bsigma}_{\bkappa,i,j,S}(\bQ)]| \cdot \left|\EE[\bW^{\bkappa,S,\bsigma}_{pq}]-\EE[\bX^{\bkappa,S,\bsigma}_{pq}]\right| + O(n^{-1/2 + \epsilon^{\prime}}),
    \intertext{Now, for the remaining terms we will not have any particular control of the difference of expectations involving $\bW$ and $\bX$, so a direct bound on these using Proposition~\ref{prop:subexp-moments} reduces this to}
    &\lesssim \sum_{\ell = 1}^N \sum_{u = 0}^4 \sum_{\bkappa\in[2L]^u} \sum_{\bsigma\in\{\id,\conj\}^u} \sum_{\substack{S \in \mathcal{S}_{p,q}^{u \times 2} \\ \text{no row of } S \text{ equals } (\emptyset, \emptyset) \\ |S| \in \{3, 4\}}} n^{u - |S| / 2} \cdot \left|\EE\left[ \Phi_{\bkappa,i,j}(\bQ) J^{\bsigma}_{\bkappa,i,j,S}(\bQ) \right]\right| + O(n^{-1/2 + \epsilon^{\prime}}).
\end{align*}

Now, consider grouping the terms in the outer sum over $\ell$ according to, firstly, whether $p_{\ell} = q_{\ell}$ or not, and second according to the size of $\{p_{\ell}, q_{\ell}\} \cap \{i, j\}$.
The total number of terms where \emph{either} $p_{\ell} = q_{\ell}$ or $|\{p_{\ell}, q_{\ell}\} \cap \{i, j\}| \geq 1$ is $O(n)$.
On the other hand, note first that, since $|J^{\bsigma}_{\bkappa,i,j,S}(\bQ)| \prec n^{-u + 8\epsilon_v}$ as we derived earlier, by Proposition~\ref{prop:resolvent-entry-bounds} (and the Cauchy-Schwarz inequality to extract the $\Phi_{\bkappa,i,j}(\bQ)$ factor) every term in the sum above is $O(n^{-|S| / 2 + 8\epsilon_v + \epsilon})$.
Further, whenever $S \in \sS^{u \times 2}_{p,q}$ and no row of $S$ equals $(\emptyset, \emptyset)$, then in particular every row contributes at least 1 to $|S|$, so $|S| \geq 3$.
Therefore, we have the simpler bound that every term in the sum above is $O(n^{-|S| / 2 + 8\epsilon_v + \epsilon)} \leq O(n^{-3/2 + 9\epsilon_v})$.
So, the sum of $O(n)$ such terms is always $O(n^{-1/2 + 9\epsilon_v})$, and up to such error, which is again subsumed in our current error term, we may restrict our attention to the case $p_{\ell} \neq q_{\ell}$ and $\{p_{\ell}, q_{\ell}\} \cap \{i, j\} = \emptyset$.
Since there are $O(n^2)$ such terms, we have:
\begin{align*}
    |\EE \wtilde{F}_{ij}(\bw) - \EE \wtilde{F}_{ij}(\bx)|
    &\lesssim \max_{\substack{p,q \in [n] \text{ distinct} \\ \{p,q\} \cap \{i, j \} = \emptyset \\ 0\le u\leq 4,\ \bkappa\in[2L]^u, \bsigma\in\{\id,\conj\}^u \\ S \in \sS^{u \times 2}_{p,q} \\ \text{no row of $S$ equals $(\emptyset, \emptyset)$} \\ |S| \in \{3, 4\}}} n^{2 + u - |S| / 2} \left| \EE\left[ \Phi_{\bkappa,i,j}(\bQ)J^{\bsigma}_{\bkappa,i,j,S}(\bQ) \right] \right| \\
    &\hspace{3cm} + O(n^{-1/2 + \epsilon^{\prime}}).
\end{align*}

Now, we note that when $\{p,q\} \cap \{i, j\} = \emptyset$, then we can improve our above bound to $|J^{\bsigma}_{\bkappa,i,j,S}(\bQ)| \prec n^{-u - |S|_0/2 + 2|S|\epsilon_v}$.
Thus, every expression in the maximum above, for a given value of $u$, is bounded by $O(n^{2 - |S|/2 - |S|_0/2 + 8\epsilon_v + \epsilon})$.
All $S$ in the maximum have $|S|_0 \geq 1$ (since otherwise they would have a row identically zero), so if $|S| = 4$ then any such term in the maximum has value $O(n^{-1/2 + 9\epsilon_v})$, again smaller than our current error term, and so such terms can effectively be ignored.

So, we may restrict our attention to terms with $|S| = 3$.
In this case, if $|S|_0 \geq 2$, then again such a term has value $O(n^{-1/2 + 9\epsilon_v})$ since $|S| / 2 + |S|_0 / 2 \geq 5/2$, and such terms can likewise be ignored.
So, we may finally restrict our attention to the case $|S| = 3$ and $|S|_0 = 1$.
In this case, we must have $u = 1$, and $S = [T \,\,\, \emptyset]$ or $S = [\emptyset \,\,\, T]$ for some $T \in \sS_{p,q,3}$.
Expanding the definitions back out, we find
\begin{align*}
&|\EE \wtilde{F}_{ij}(\bw) - \EE \wtilde{F}_{ij}(\bx)| \\
&\lesssim \max_{\substack{i, j \in [n] \\ a\in [L], \kappa\in\{2a-1,2a\}\\ p,q\in [n]\text{ distinct } \\ \{p,q\} \cap \{i, j \} = \emptyset \\ T \in \sS_{p,q,3} \\ \sigma\in\{\id,\conj\}  } } n^{3/2} \bigg| \EE\bigg[\Phi_{\kappa,i,j}(\bQ) R^{\sigma}_{Q^{(a)}}(i,\alpha_{T,1})
R^{\sigma}_{Q^{(a)}}(\beta_{T,1},\alpha_{T,2})
R^{\sigma}_{Q^{(a)}}(\beta_{T,2},\alpha_{T,3})
R^{\sigma}_{Q^{(a)}}(\beta_{T,3},v^{(a)}) \\
&\hspace{3.2cm} \cdot
R^{\sigma}_{Q^{(a)}}(v^{(a)},j) \bigg] \bigg|  + O(n^{-1/2 + \epsilon^{\prime}}),
\intertext{where the case of $S = [\emptyset \,\,\, T]$ is treated analogously with the decoupling performed from the right. Now, recall that given $T \in \sS_{p,q,3}$, we have $\alpha_{T,d}, \beta_{T,d} \in \{p,q\}$, and $\alpha_{T,d}\neq \beta_{T,d}$ for $d=1,2,3$ (so one is $p$ and the other is $q$). Thus, whenever $\beta_{T,1} \neq \alpha_{T,2}$ or $\beta_{T,2} \neq \alpha_{T,3}$, by Proposition~\ref{prop:resolvent-entry-bounds} the expression in the maximum is $O(n^{3/2} \cdot n^{-2 + 2\epsilon_v + \epsilon}) \leq O(n^{-1/2 + 2\epsilon_v + \epsilon})$, and the corresponding term is absorbed into the current error term. So, we may further reduce to a specific combination of resolvent quadratic forms, removing the dependence on $T$ and arriving at an entirely concrete expression:}
&\lesssim \max_{\substack{i, j \in [n] \\  a\in[L], \kappa \in \{2a-1,2a\}\\p,q \in [n] \text{ distinct } \\ \{p,q\} \cap \{i, j \} = \emptyset \\ \sigma\in\{\id,\conj\} }} n^{3/2} \bigg| \EE \bigg[
\Phi_{\kappa,i,j}(\bQ)
R^{\sigma}_{Q^{(a)}}(i,p)
R^{\sigma}_{Q^{(a)}}(p,p)
R^{\sigma}_{Q^{(a)}}(q,q)
R^{\sigma}_{Q^{(a)}}(q,v^{(a)})
R^{\sigma}_{Q^{(a)}}(v^{(a)},j)
\bigg] \bigg| \\&\hspace{3cm} + O(n^{-1/2 + \epsilon^{\prime}}).
\end{align*}
Expanding out the definition of $\Phi_{\kappa,i,j}(\bQ)$, we see that, for $\kappa\in\{2a-1,2a\}$, this factor is one of the first real partial derivatives of $\wtilde{\phi}$ evaluated at $(n\cdot R_{Q^{(c)}}(i,v^{(c)})R_{Q^{(c)}}(v^{(c)},j)
)_{c=1}^{L}$. Thus, by the boundedness of the derivatives of $\wtilde{\phi}$, it is a bounded Lipschitz function of this $L$-tuple. The proof is then completed upon using the result of Lemma~\ref{lem:bad-res-term} below.
\end{proof}

The last technical ingredient we will need is the following.
We remove some of the specific details of the form of the $\Phi_{\kappa,i,j}(\bQ)$ factor for the sake of brevity in the proof.
\begin{lem}
    \label{lem:bad-res-term}
In the above setting, suppose $i, j \in [n]$ and that $p,q \in [n] \setminus \{i, j\}$ are distinct.
Fix $a_0\in[L]$ and $\sigma\in\{\id,\conj\}$. Let $f: \CC^{L} \to \CC$ be a bounded Lipschitz function.
Then, for any $\epsilon > 0$,
\begin{align*}
    &\bigg| \EE\bigg[ f \left( \left( n\cdot R_{Q^{(a)}}(i,v^{(a)})R_{Q^{(a)}}(v^{(a)},j) \right)_{a=1}^{L} \right) \\
    &\hspace{2cm}\cdot R^{\sigma}_{Q^{(a_0)}}(i,p) R^{\sigma}_{Q^{(a_0)}}(p,p) R^{\sigma}_{Q^{(a_0)}}(q,q) R^{\sigma}_{Q^{(a_0)}}(q,v^{(a_0)}) R^{\sigma}_{Q^{(a_0)}}(v^{(a_0)},j) \bigg] \bigg| \lesssim_{\epsilon,f} n^{-2+4\epsilon_v+\epsilon}.
\end{align*}
\end{lem}

We note that this will conclude the proof of Theorem~\ref{thm:main-univ}: Lemma~\ref{lem:bad-res-term} completes the proof of Lemma~\ref{lem:tilde-F-bound}, and combining this with Lemma~\ref{lem:F-bound} in turn completes the proof of Theorem~\ref{thm:main-univ}, as described earlier.

\subsection{Resolvent monomial expectation via decoupling: Proof of Lemma~\ref{lem:bad-res-term}}

It suffices to prove the case $\sigma=\id$. Fix $i,j,p,q$ and $a_0$ as in Lemma~\ref{lem:bad-res-term}. For each $a\in[L]$, let $Q^{(a,p)}\in\CC^{(n-1)\times (n-1)}_{\herm}$ be obtained from $Q^{(a)}$ by removing the $p$-th row and column.
Also, denote by \[\sF^{(p)} \colonequals \sigma\{Q^{(a)}_{\mu\nu}:a\in [L], p \notin \{\mu,\nu\}\}\] the $\sigma$-algebra generated by all entries of the tuple $\bQ$ whose matrix indices do not
touch $p$. We write $\EE_{p}[\ \cdot\ ] \colonequals \EE[\ \cdot\ |\sF^{(p)}]$ for the operation of taking conditional expectation with respect to this $\sigma$-algebra. In simpler language, $\EE_p$ is the operation of averaging over the $p$-th row and column of all matrices $Q^{(a)}$, $a\in[L]$, jointly.

Since in this section we will only work with resolvents of matrices related to $\bQ$, and the proof is the same for both signs, let us slightly change notation and define, for $a\in[L]$,
\begin{align*}
    R_a &\colonequals R_{Q^{(a)}}(\wtilde{\lambda}),
\end{align*}
and write $R^{(p)}_a$ for $R^{\sigma}_{Q^{(a,p)}}(\wtilde{\lambda})$ extended to have dimension $n \times n$ by adding a $p$-th row and column equal to zero.
Here $\wtilde{\lambda} = \theta + \theta^{-1} + \cpxi n^{-1/2}$ is the same complex value at which we have been evaluating all resolvents in the previous section as well.

For the distinguished tuple coordinate $a_0$ appearing in Lemma~\ref{lem:bad-res-term}, we abbreviate \[ R \colonequals R_{a_0},\qquad R^{(p)} \colonequals R^{(p)}_{a_0}, \qquad v \colonequals v^{(a_0)}\] while the test function $f$ may still depend on the whole tuple $\left( n \cdot R_a(i,v^{(a)}) R_a(v^{(a)},j) \right)_{a=1}^{L}$. The identities below are stated for a general tuple coordinate $a\in[L]$; in the proof of Lemma~\ref{lem:bad-res-term} we will mostly use them with $a=a_0$, in which case they are written using the abbreviated notation $R,R^{(p)},v$.

For the remaining analysis, we use the method of \emph{fluctuation averaging} (see Section~7 of \cite{BGK-2016-LecturesLocalSemicircleLaw}).
This method is based on the observation that $\EE_p R(i, p)$ is much smaller than the typical value of $R(i, p)$ with no averaging.
The following are the main technical devices used to make this argument, identities relating resolvents to their minors.

\begin{prop}\label{prop:resolvent-minor}
For any $a\in[L]$ and any $\mu,\nu,k\in [n]$ with $k \notin \{\mu,\nu\}$, we have
\[
R_{a}(\mu,\nu) = R^{(k)}_{a}(\mu,\nu) + \frac{R_{a}(\mu,k) R_{a}(k,\nu)}{R_{a}(k,k)}.
\]
\end{prop}
\begin{proof}
    This is a special case of the formula for the inverse of the entries of a block matrix, where we view the $(k, k)$ entry as a $1 \times 1$ block.
\end{proof}

\begin{prop}\label{prop:resolvent-minor-entries}
For any $a\in[L]$ and any distinct $\mu, \nu \in [n]$, we have
\begin{align*}
R_{a}(\mu,\nu) &= R_{a}(\mu,\mu) \sum_{k\neq \mu} Q^{(a)}_{\mu k} R^{(\mu)}_{a}(k,\nu), \\ R_{a}(\mu,\nu) &= R_{a}(\nu,\nu) \sum_{k\neq \nu} R^{(\nu)}_{a}(\mu,k) Q^{(a)}_{k\nu}.
\end{align*}
\end{prop}
\begin{proof}
This is a special case of the same formula referenced above, but now where we view the $(\mu, \mu)$ or $(\nu, \nu)$ entry as the $1 \times 1$ block with respect to which we expand the inverse.
\end{proof}

The following is then the main estimate that we will use.
Note that Proposition~\ref{prop:resolvent-entry-bounds} implies only the a priori bound $\EE |R_{a}(m, r)|^2 \lesssim_{\epsilon} n^{-1 + \epsilon}$, on which this result improves considerably.
\begin{lem}
    \label{lem:Ei-estimate}
    For any $a\in[L]$ and any distinct $i, p \in [n]$ , for any $\epsilon > 0$,
    \[ \EE |\EE_p R_a(i, p)|^2  \lesssim_{\epsilon} n^{-2 + \epsilon}. \]
\end{lem}
\begin{proof}
Fix $a\in[L]$ and distinct $i, p \in [n]$. Recall the abbreviation
\[R=R_a,\qquad R^{(p)}=R^{(p)}_{a}, \qquad Q=Q^{(a)}.\]
By Proposition~\ref{prop:resolvent-minor-entries}, for indices $i\neq p$,
\[
R(i,p) = R(p,p) \sum_{k\neq p} R^{(p)}(i,k)Q_{kp}.
\]
Decompose $R(p,p)=G_{\mu_{\sc}}(\wtilde{\lambda})+\big(R(p,p)-G_{\mu_{\sc}}(\wtilde{\lambda})\big)$, where $G_{\mu_{\sc}}$ is the Cauchy transform of the semicircle law. Then
\[
R(i,p) = G_{\mu_{\sc}}(\wtilde{\lambda}) \cdot \sum_{k\neq p} R^{(p)}(i,k)Q_{kp} + \big(R(p,p)-G_{\mu_{\sc}}(\wtilde{\lambda})\big) \cdot \sum_{k\neq p} R^{(p)}(i,k)Q_{kp}.
\]
Now we consider taking the conditional expectation $\EE_{p}$ (conditional on $\sF^{(p)}$).
Since $R^{(p)}$ is $\sF^{(p)}$-measurable (i.e., independent of row and column $p$ of $Q=Q^{(a)}$) and by the definition of generalized Wigner matrices we have $\EE_{p}[Q_{kp}]=0$, we have
\[
\EE_{p}\Big[G_{\mu_{\sc}}(\wtilde{\lambda}) \cdot \sum_{k\neq p} R^{(p)}(i,k)Q_{kp}\Big] = G_{\mu_{\sc}}(\wtilde{\lambda}) \cdot \sum_{k\neq p} R^{(p)}(i,k) \EE_{p}[Q_{kp}]=0.
\]
So, we are left with just
\[
\EE_{p}[R(i,p)]=\EE_{p}\left[\big(R(p,p)-G_{\mu_{\sc}}(\wtilde{\lambda})\big) \cdot \sum_{k\neq p} R^{(p)}(i,k)Q_{kp}\right].
\]
Applying the Cauchy-Schwarz inequality on this conditional expectation, we have
\begin{equation}\label{eq:bad-term-1}
\Big|\EE_{p}[R(i,p)]\Big|^{2} \le \Big(\EE_{p}|R(p,p)-G_{\mu_{\sc}}(\wtilde{\lambda})|^2\Big)\left(\EE_{p}\left|\sum_{k\neq p} R^{(p)}(i,k)Q_{kp}\right|^2\right),
\end{equation}
and taking the unconditional expectation,
\begin{equation}\label{eq:bad-term-2}
\EE\Big|\EE_{p}[R(i,p)]\Big|^{2} \le \EE\left[\EE_{p}|R(p,p)-G_{\mu_{\sc}}(\wtilde{\lambda})|^2\ \cdot \EE_{p}\left|\sum_{k\neq p} R^{(p)}(i,k)Q_{kp}\right|^2\right].
\end{equation}
Expanding the second factor, we see that we may simplify it to
\begin{align*}
\EE_{p}\Big|\sum_{k\neq p} R^{(p)}(i,k)Q_{kp}\Big|^2 &=
\EE_{p}\bigg[\sum_{k,l\neq p} R^{(p)}(i,k)\overline{R^{(p)}(i,l)} \cdot Q_{kp}\overline{Q_{lp}}\bigg]\\
&=\sum_{k,l\neq p}
R^{(p)}(i,k)\overline{R^{(p)}(i,l)}
\EE_p\left[Q_{kp}\overline{Q_{lp}}\right] \\
&=
\sum_{k\neq p}
\EE_p|Q_{kp}|^2\, |R^{(p)}(i,k)|^2
\end{align*}
We have $\EE_{p}\left[|Q_{kp}|^2\right]=\EE\left[|Q_{kp}|^2\right] \leq Cn^{-1}$ for some $C > 0$ using our definition of generalized Wigner matrices and that $Q$ is formed from a matrix of this kind by setting some entries to zero.
So, we have
\[
\EE_{p}\Big|\sum_{k\neq p} R^{(p)}(i,k)Q_{kp}\Big|^2\lesssim n^{-1} \cdot \sum_{k}\big|R^{(p)}(i,k)\big|^2,
\]
which is an $\sF^{(p)}$-measurable random variable (not depending on row and column $p$ of $Q$).

Substituting back, we may then rewrite as a full expectation,
\begin{align*}
    \EE\Big|\EE_{p}[R(i,p)]\Big|^{2}
    &\lesssim n^{-1} \cdot \EE \left[ \EE_{p} [|R(p,p) - G_{\mu_{\sc}}(\wtilde{\lambda})|^2] \cdot \sum_{k}\big|R^{(p)}(i,k)\big|^2\right] \\
    &= n^{-1} \cdot \EE \left[ |R(p,p) - G_{\mu_{\sc}}(\wtilde{\lambda})|^2 \cdot \sum_{k}\big|R^{(p)}(i,k)\big|^2\right]
    \intertext{using that $\sum_k |R^{(p)}(i,k)|^2 = \|(R^{(p)})^{*}e_{i}\|^2 \leq \|R^{(p)}\|^2$, we have}
    &\leq n^{-1} \cdot \EE \left[ |R(p,p) - G_{\mu_{\sc}}(\wtilde{\lambda})|^2 \cdot \|R^{(p)}\|^2\right] \\
    &\leq n^{-1} \cdot \left(\EE |R(p,p) - G_{\mu_{\sc}}(\wtilde{\lambda})|^4\right)^{1/2} \cdot \left(\EE \|R^{(p)}\|^4\right)^{1/2}.
\end{align*}
By Proposition~\ref{prop:res-moments}, the resolvent norm factor is $O(1)$.
By Corollary~\ref{cor:ill-Q} (the isotropic local law for $Q$), Proposition~\ref{prop:stochastic-domination-to-moment-bound}, and Proposition~\ref{prop:res-moments}, the second factor is $O_{\epsilon}(n^{-1 + \epsilon})$ for any $\epsilon > 0$.
The result then follows from combining these bounds.
\end{proof}

Now we use this to give the proof of the main result of this section.
\begin{proof}[Proof of Lemma~\ref{lem:bad-res-term}]
We note before continuing that the bounds on resolvent entries and quadratic forms with $e_{\alpha}$ and $v^{(a)}$ proved in Proposition~\ref{prop:resolvent-entry-bounds} all apply equally well to $R^{(p)}_{a}$, for every $a\in[L]$. Indeed, $R^{(p)}_{a}$ is zero outside of a principal $(n - 1) \times (n - 1)$ submatrix that equals the resolvent of $Q^{(a,p)}$, where $Q^{(a,p)}$ is either a Wigner matrix of dimension $n - 1$ or such a matrix with one or two entries set to zero, to which Proposition~\ref{prop:resolvent-entry-bounds} applies directly.

We would like to move towards applying Lemma~\ref{lem:Ei-estimate}.
To that end, define
\begin{align*}
\sA &\colonequals f\bigg( \left(
    n\cdot R_a(i,v^{(a)})R_a(v^{(a)},j)
    \right)_{a=1}^{L} \bigg) \cdot R(p,p) R(q,q) R(q, v) R(v, j), \\
\sA^{(p)} &\colonequals f\bigg( \left(
    n\cdot R_a^{(p)}(i,v^{(a)})R_a^{(p)}(v^{(a)},j)
    \right)_{a=1}^{L}\bigg) \cdot G_{\mu_{\sc}}(\wtilde{\lambda}) R^{(p)}(q,q) R^{(p)}(q, v) R^{(p)}(v, j).
\end{align*}
Recall that we are trying to prove a bound on $\EE [\sA \cdot R(i, p)]$.
Since $\sA^{(p)}$ is $\sF^{(p)}$-measurable, we will try to replace $\sA$ by $\sA^{(p)}$, noting that we expect these to be close by since we make only small perturbations by replacing various $R_a$ by $R_a^{(p)}$, and another small perturbation in replacing $R(p, p)$ by $G_{\mu_{\sc}}(\widetilde{\lambda})$ according to the local law in Corollary~\ref{cor:ill-Q}.
We have:
\begin{align*}
\big|\EE[\sA \cdot R(i,p)]\big|
&= |\EE[(\sA-\sA^{(p)}) \cdot R(i,p)] + \EE[\sA^{(p)} \cdot R(i,p)]| \\
&\le \big|\EE[(\sA-\sA^{(p)}) \cdot R(i,p)]\big| + \big|\EE[\sA^{(p)} \cdot \EE_{p}[R(i,p)]]\big|, \numberthis \label{eq:bad-term-5}
\end{align*}
where in the second term we have used that $\sA^{(p)}$ is $\sF^{(p)}$-measurable.

We now bound each term individually. Consider the second term of \eqref{eq:bad-term-5} first. Using Cauchy–Schwarz,
\[
\Big|\EE\big[\sA^{(p)} \cdot \EE_{p}[R(i,p)]\big]\Big|\le \Big(\EE|\sA^{(p)}|^2\Big)^{1/2}\Big(\EE\big|\EE_{p}[R(i,p)]\big|^2\Big)^{1/2}
\]
As shown in Lemma~\ref{lem:Ei-estimate}, we already have
\[
\EE\big|\EE_{p}[R(i,p)]\big|^2 \lesssim n^{-2 + \epsilon}
\]
for any $\epsilon > 0$.
It remains to bound $\EE|\sA^{(p)}|^2$.
Since $f$ is bounded, we have by Proposition~\ref{prop:resolvent-entry-bounds}
\[
|\sA^{(p)}| \le \|f\|_{\infty} \cdot |G_{\mu_{\sc}}(\wtilde{\lambda})| \cdot |R^{(p)}(q,q)| \cdot |R^{(p)}(q,v)| \cdot |R^{(p)}(v,j)| \prec n^{-1+2\epsilon_{v}},
\]
and by the expectation bounds in Proposition~\ref{prop:resolvent-entry-bounds} we have
\[
\EE|\sA^{(p)}|^2 \lesssim n^{-2+4\epsilon_{v}+ \epsilon}.
\]
for any $\epsilon > 0$.
Thus, the entire second term of \eqref{eq:bad-term-5} is bounded by
\begin{equation}\label{eq:bad-term-3}
\big|\EE[\sA^{(p)} \cdot \EE_{p}[R(i,p)]]\big| \lesssim n^{-1+2\epsilon_{v}+\epsilon / 2}\cdot n^{-1+ \epsilon/2} = n^{-2+2\epsilon_{v}+\epsilon}.
\end{equation}

Now we consider the first term in \eqref{eq:bad-term-5}.
We first reorganize the expression $\sA - \sA^{(p)}$.
Define
\begin{align*}
\sB &\colonequals R(p,p) R(q,q) R(q, v) R(v, j), \\
\sB^{(p)} &\colonequals G_{\mu_{\sc}}(\wtilde{\lambda}) R^{(p)}(q,q) R^{(p)}(q, v) R^{(p)}(v, j).
\end{align*}
Then, we have
\begin{equation}
\begin{aligned}
&\sA-\sA^{(p)} = \left( f\Big(\big( n\cdot R_a(i,v^{(a)})R_a(v^{(a)},j) \big)_{a=1}^{L}\Big) - f\Big(\big( n\cdot R_a^{(p)}(i,v^{(a)})R_a^{(p)}(v^{(a)},j) \big)_{a=1}^{L}\Big) \right) \cdot \sB\\ &\hspace{1cm}+ f\left( \left( n\cdot R_a^{(p)}(i,v^{(a)})R_a^{(p)}(v^{(a)},j) \right)_{a=1}^{L}\right) \cdot (\sB - \sB^{(p)}).
\end{aligned}
\label{eq:AB}
\end{equation}
We first bound the first term on the right-hand side of \eqref{eq:AB}. For each $a\in[L]$, Proposition~\ref{prop:resolvent-minor} gives
\[
R_a(i,v^{(a)})-R_a^{(p)}(i,v^{(a)}) = \frac{R_a(i,p)R_a(p,v^{(a)})}{R_a(p,p)},
\]
and
\[
R_a(v^{(a)},j)-R_a^{(p)}(v^{(a)},j) = \frac{R_a(v^{(a)},p)R_a(p,j)}{R_a(p,p)}.
\]
Therefore, by Proposition~\ref{prop:resolvent-entry-bounds},
\[
|R_a(i,v^{(a)})-R_a^{(p)}(i,v^{(a)})| \prec n^{-1+\epsilon_v}, \qquad |R_a(v^{(a)},j)-R_a^{(p)}(v^{(a)},j)| \prec n^{-1+\epsilon_v}.
\]
It follows that
\begin{align*}
& \left| n\cdot R_a(i,v^{(a)})R_a(v^{(a)},j) - n\cdot R_a^{(p)}(i,v^{(a)})R_a^{(p)}(v^{(a)},j) \right| \le n|R_a(i,v^{(a)})-R_a^{(p)}(i,v^{(a)})||R_a(v^{(a)},j)| \\ &\hspace{3.2cm} + n|R_a^{(p)}(i,v^{(a)})|\,|R_a(v^{(a)},j)-R_a^{(p)}(v^{(a)},j)| \prec n^{-1/2+2\epsilon_v}.
\end{align*}

Since $f$ is Lipschitz and $L$ is fixed
\[\bigg|f\Big(\big( n\cdot R_a(i,v^{(a)})R_a(v^{(a)},j) \big)_{a=1}^{L}\Big) - f\Big(\big( n\cdot R_a^{(p)}(i,v^{(a)})R_a^{(p)}(v^{(a)},j) \big)_{a=1}^{L}\Big)\bigg| \prec n^{-1/2+2\epsilon_v}.\]
Again by Proposition~\ref{prop:resolvent-entry-bounds}, we have
\[ |\sB| \prec n^{-1 + 2\epsilon_v}, \]
so by Proposition~\ref{prop:stoch-dom-sum-prod}, the first term of \eqref{eq:AB} is polynomially stochastically dominated as
\begin{equation}
\bigg|\Big(f\Big(\big( n\cdot R_a(i,v^{(a)})R_a(v^{(a)},j) \big)_{a=1}^{L}\Big) - f\Big(\big( n\cdot R_a^{(p)}(i,v^{(a)})R_a^{(p)}(v^{(a)},j) \big)_{a=1}^{L}\Big) \Big)\cdot \sB \bigg| \prec n^{-3/2 + 4\epsilon_v}. \label{eq:AB1}
\end{equation}

Meanwhile, for the second term of \eqref{eq:AB}, we may apply a telescoping expansion to $\sB - \sB^{(p)}$, obtaining
\begin{align*}
\sB-\sB^{(p)} &= \big(R(p,p)-G_{\mu_{\sc}}(\wtilde\lambda)\big)R(q,q)R(q,v)R(v,j)\\ &\hspace{1cm} + G_{\mu_{\sc}}(\wtilde\lambda) \big(R(q,q)-R^{(p)}(q,q)\big)R(q,v)R(v,j)\\ &\hspace{1cm} + G_{\mu_{\sc}}(\wtilde\lambda) R^{(p)}(q,q) \big(R(q,v)-R^{(p)}(q,v)\big)R(v,j)\\ &\hspace{1cm} + G_{\mu_{\sc}}(\wtilde\lambda) R^{(p)}(q,q)R^{(p)}(q,v) \big(R(v,j)-R^{(p)}(v,j)\big).
\end{align*}
We now control the individual differences that appear here. The first term above is
\[\big(R(p,p)-G_{\mu_{\sc}}(\wtilde\lambda)\big)R(q,q)R(q,v)R(v,j) \prec n^{-1/2+\epsilon} \cdot 1\cdot n^{-1/2+\epsilon_v}\cdot n^{-1/2+\epsilon_v} \prec n^{-3/2+2\epsilon_v+\epsilon} \]
For the remaining differences, Proposition~\ref{prop:resolvent-minor} gives
\[R(q,q)-R^{(p)}(q,q) = \frac{R(q,p)R(p,q)}{R(p,p)} \prec n^{-1}, \]
and
\[R(q,v)-R^{(p)}(q,v) = \frac{R(q,p)R(p,v)}{R(p,p)} \prec n^{-1+\epsilon_v},\]
with the analogous bound
\[R(v,j) - R^{(p)}(v,j) \prec n^{-1+\epsilon_v}.\]
Putting everything together, we find
\[ |\sB - \sB^{(p)}| \prec n^{-3/2 + 2\epsilon_v+\epsilon}, \]
and therefore the whole second term in \eqref{eq:AB} is bounded as
\[ \bigg|f\Big(\big( n\cdot R_a^{(p)}(i,v^{(a)})R_a^{(p)}(v^{(a)},j) \big)_{a=1}^{L}\Big) \cdot (\sB - \sB^{(p)})\bigg| \prec n^{-3/2 + 2\epsilon_v + \epsilon} \]
as well, since $f$ is bounded.
Combining this with \eqref{eq:AB1} and plugging both bounds into \eqref{eq:AB}, we have
\[ |\sA - \sA^{(p)}| \prec n^{-3/2 + 4\epsilon_v + \epsilon}.\]
Since $R(i,p)\prec n^{-1/2}$, we obtain
\[ |(\sA - \sA^{(p)}) \cdot R(i, p)| \prec n^{-2 + 4\epsilon_v}, \]
and thus
\[ |\EE[(\sA - \sA^{(p)}) \cdot R(i, p)]| \lesssim_{\epsilon,f} n^{-2 + 4\epsilon_v + \epsilon}. \]
Finally, applying our bounds to either term of \eqref{eq:bad-term-5}, we get
\[ |\EE[\sA \cdot R(i, p)]| \lesssim_{\epsilon} n^{-2 + 4\epsilon_v + \epsilon}, \]
completing the proof.
\end{proof}

\section{Gaussian formula for entrywise averages: Proof of Theorem~\ref{thm:main-gauss}}

Recall the setting of the Theorem: unlike the above proof of Theorem~\ref{thm:main-univ}, the statement of this result depends on the fields to which the entries of the matrices in the tuple $\bW$ belong.
So, as in the Introduction, for each $a\in[L]$, let us write
\[ \FF_a = \left\{\begin{array}{ll} \RR & \text{if } W^{(a)} \text{ is real symmetric}, \\ \CC & \text{if } W^{(a)} \text{ is complex Hermitian}\end{array}\right\}. \]
Recall that we write $\sN_{\FF_a}(0, 1)$ for the standard Gaussian scalar distribution associated to the field $\FF_a$ (normalized so that $\EE |Z|^2 = 1$ for $Z \in \sN_{\FF_a}(0, 1)$ in either case), and we write $\sN_{\FF_a}(0, I_n)$ for the law of a vector of $n$ i.i.d.\ such Gaussians.

In Theorem~\ref{thm:main-gauss}, we have $\psi: \CC^L \times \CC^L \to \RR$ a $\sC^5$ function with bounded values and first five derivatives, and an associated function $\Psi: (\CC^{n \times n}_{\herm})^L \times (\CC^{n \times n}_{\herm})^L \to \RR$ defined by
    \[ \Psi(\bA, \bB) \colonequals \frac{1}{n^2} \sum_{i, j = 1}^n \psi\left(\left(A^{(a)}_{ij}\right)_{a=1}^{L}, \left(B^{(a)}_{ij}\right)_{a=1}^{L}\right). \]
Theorem~\ref{thm:main-gauss} then bounds, for $\what{v}^{(W,a)} = v_1\left(\theta v^{(a)}v^{(a)*} + W^{(a)}\right)$ and $g^{(a)} \sim \sN_{\FF_a}(0, I_n)$, with the Gaussian vectors $g^{(a)}$ independent for $a\in[L]$, the quantity
\begin{align*}
&\Bigg|\EE\Psi\left(\left(nv^{(a)}v^{(a)*}\right)_{a=1}^L,\left(n\what{v}^{(W,a)}\what{v}^{(W,a)*}\right)_{a=1}^L \right) \\ &\hspace{2cm} -\EE\Psi\left(\left(nv^{(a)}v^{(a)*}\right)_{a=1}^L, \left((\rho(\theta)\sqrt{n} v^{(a)}+\tau(\theta)g^{(a)})(\rho(\theta)\sqrt{n}v^{(a)}+\tau(\theta)g^{(a)})^* \right)_{a=1}^L \right) \Bigg|.
\end{align*}

Before proceeding, let us establish a few tools for working with such functions $\Psi$.
First, it follows immediately from our assumptions that $\Psi$ is $O(n^{-1})$-Lipschitz:
\begin{prop}
    \label{prop:Psi-lip}
    For any $\bA, \bB, \bA^{\prime}, \bB^{\prime} \in (\CC^{n \times n}_{\herm})^L$, writing $\bA = (A^{(a)})_{a=1}^L$, $\bB = (B^{(a)})_{a=1}^L$, $\bA^{\prime} = (A^{\prime(a)})_{a=1}^L$ and $\bB^{\prime} = (B^{\prime(a)})_{a=1}^L$, we have
    \[ |\Psi(\bA, \bB) - \Psi(\bA^{\prime}, \bB^{\prime})| \leq \frac{C_{\psi}} {n} \sum_{a=1}^L \left(\|A^{(a)} - A^{\prime(a)}\|_{\F} + \|B^{(a)} - B^{\prime(a)}\|_{\F}\right), \]
    for a constant $C_{\psi}$ depending only on $L$ and the bounds on $\psi$ and its first derivatives.
\end{prop}
\begin{proof}
For each $i,j\in[n]$, write $\bA_{ij} \colonequals (A^{(a)}_{ij})_{a=1}^L$, $\bB_{ij} \colonequals (B^{(a)}_{ij})_{a=1}^L$, and define $\bA^{\prime}_{ij}$, $\bB^{\prime}_{ij}$ similarly.
Then, the result follows since we have entrywise
\begin{align*}
\left|\psi(\bA_{ij}, \bB_{ij}) - \psi(\bA^{\prime}_{ij}, \bB^{\prime}_{ij})\right|
&\leq C_{\psi} \sum_{a=1}^L \left(|A^{(a)}_{ij} - A^{\prime(a)}_{ij}| + |B^{(a)}_{ij} - B^{\prime(a)}_{ij}|\right)
\end{align*}
for a suitable $C_{\psi}$.
Thus, we have
\begin{align*}
|\Psi(\bA, \bB) - \Psi(\bA^{\prime}, \bB^{\prime})| &\leq \frac{C_{\psi}}{n^2}\sum_{i, j = 1}^n \sum_{a=1}^L \left(\big|A^{(a)}_{ij} - A^{\prime(a)}_{ij}\big| + \big|B^{(a)}_{ij} - B^{\prime(a)}_{ij}\big|\right) \\ &\leq \frac{C_{\psi}}{n} \sum_{a=1}^L\left(\|A^{(a)} - A^{\prime(a)}\|_{\F} + \|B^{(a)} - B^{\prime(a)}\|_{\F}\right)
\end{align*}
using the Cauchy-Schwarz inequality in the sum over $(i,j)$ for each fixed $a$.
\end{proof}
\noindent
Next, we give a standard bound on the change in eigenvector projections under small perturbations of a matrix; the result we give is a simple reformulation of the Davis-Kahan inequality; see \cite{davis1970rotation,yu2014usefulvariantdaviskahantheorem}.
\begin{prop}
    \label{prop:dk}
    Let $H, \Delta \in \CC^{n \times n}_{\herm}$.
    Let $\lambda_1(H)$ be a simple eigenvalue and suppose $\|\Delta\| < \frac{1}{2}\big(\lambda_1(H) - \lambda_2(H)\big)$.
    Then,
    \[ \|v_1(H)v_1(H)^* - v_1(H + \Delta) v_1(H + \Delta)^*\|_{\F} \leq 2\sqrt{2} \cdot \frac{\|\Delta\|}{\lambda_1(H) - \lambda_2(H)}. \]
\end{prop}
\noindent
Combining the two results and rescaling per our setting gives the following useful corollary.
\begin{cor}
\label{cor:dk}
Let $\bW=(W^{(a)})_{a=1}^L$ and $\bW^{\prime}=(W^{\prime(a)})_{a=1}^L$ be two tuples of Hermitian matrices and for each $a\in[L]$, consider $H^{(a)} = \theta v^{(a)}v^{(a)*} + W^{(a)}$ and $H^{\prime(a)} = \theta v^{(a)}v^{(a)*} + W^{\prime(a)}$.
Write $\what{v}^{(W,a)} = v_1(H^{(a)})$ and $\what{v}^{(W^{\prime},a)} = v_1(H^{\prime(a)})$ for the associated top eigenvectors.
If, for every $a\in[L]$, $\|W^{(a)} - W^{\prime(a)}\| < \frac{1}{2} \big(\lambda_1(H^{(a)}) - \lambda_2(H^{(a)})\big)$, then
\begin{align*}
&\Bigg|\Psi\left(\left(nv^{(a)}v^{(a)*}\right)_{a=1}^L, \left(n\what{v}^{(W,a)}\what{v}^{(W,a)*}\right)_{a=1}^L \right) - \Psi\left(\left(nv^{(a)}v^{(a)*}\right)_{a=1}^L, \left(n\what{v}^{(W',a)}\what{v}^{(W',a)*}\right)_{a=1}^L \right) \Bigg|  \\ &\hspace{2cm}\leq 2\sqrt{2} C_{\psi} \sum_{a=1}^L\frac{\|W^{(a)}-W^{\prime(a)}\|}{\lambda_1(H^{(a)})-\lambda_2(H^{(a)})},
\end{align*}
where $C_{\psi}$ is the constant from Proposition~\ref{prop:Psi-lip}.
\end{cor}
\noindent
By Theorem~\ref{thm:spiked-non-asymp}, applied separately to each coordinate $a\in[L]$, and since $L$ is fixed, we have that $\min_{a\in[L]}\left(\lambda_1(H^{(a)}) - \lambda_2(H^{(a)})\right) = \Omega(1)$ with polynomially high probability in all settings below where the matrices $W^{(a)}$ are generalized Wigner matrices with uniformly controlled parameters.
Therefore, in effect, coordinatewise perturbations of the tuple $\bW$ satisfying $\max_{a\in[L]}\|W^{(a)} - W^{\prime(a)}\| = o(1)$ have a small effect on the value of $\Psi$ we are interested in.
We will use this observation several times below.

Returning to the main proof, for each $a\in[L]$, let us define
\[ \GFE_a(n) \colonequals \left\{\begin{array}{ll} \GOE(n) & \text{if } \FF_a = \RR, \\ \GUE(n) & \text{if } \FF_a = \CC\end{array}\right\}. \]
We write $\bX\sim \prod_{a=1}^L\GFE_a(n)$ to mean that $X^{(1)},\ldots,X^{(L)}$ are independent and $X^{(a)}\sim\GFE_a(n)$ for each $a\in[L]$.

We prove the Theorem in several steps. First, we show using the above argument that we may, at the cost of a small error, replace the weakly Wigner tuple $\bW$ with a tuple $\bW^{(\reg)}$ such that each coordinate is a generalized Wigner matrix satisfying the exact variance normalization \eqref{eq:gen-wigner-var-sums}.

\begin{thm}[Special case of symmetric Sinkhorn scaling theorem, Theorem 5.4 of \cite{idel2016reviewmatrixscalingsinkhorns}]\label{thm:symmetric-sinkhorn-scaling}
Let $n \geq 2$ and $S\in\RR^{n\times n}_{\sym}$ be a non-negative matrix with $S_{ij} > 0$ for all $i, j \in [n]$ distinct. Then, there exists a diagonal matrix $X=\Diag(x_1,\ldots,x_n)$ with positive entries such that $X S X$ has row sums equal $1$.
\end{thm}

\begin{lem}\label{lem:sinkhorn-near-constant}
Let $S\colonequals (s_{ij})_{i,j=1}^n\in\RR^{n\times n}_{\sym}$ be a nonnegative matrix. Suppose that,
\[\max_{i\in[n]}\sum_{j=1}^n\bigg|s_{ij} - \frac{1}{n}\bigg| \lesssim n^{-1} + n^{-\epsilon_W},\]
and $s_{ij}>0$ for all $i\neq j$. Then, for sufficiently large $n$, there exist positive $d_1,\ldots,d_n$ such that
\begin{equation}\label{eq:constant}
\sum_{j=1}^n\frac{s_{ij}}{d_i d_j}=1,\qquad i\in[n]
\end{equation}
with \begin{equation}\label{eq:d_i-1}
\max_{i\in[n]}|1-d_i|\lesssim n^{-1} + n^{-\epsilon_W}.
\end{equation}
\end{lem}
\begin{proof}
By Theorem~\ref{thm:symmetric-sinkhorn-scaling}, there exists positive $X=\Diag(x_1,\ldots,x_n)$ such that
$x_i\sum_{j=1}^n s_{ij} x_j = 1$ for all $i\in[n]$.
Setting $d_i=x_i^{-1}$ gives \eqref{eq:constant}.
Thus, for all $i\in[n]$, $d_i=\sum_{j=1}^n\frac{s_{ij}}{d_j}$.
To show \eqref{eq:d_i-1}, first define $d_{i_M} = \max_{i\in[n]}d_i$ and $d_{i_m} = \min_{i\in[n]}d_i$. Then
$|d_{i_M}-d_{i_m}| = \left|\sum_{j=1}^n \frac{s_{i_M j}-s_{i_m j}}{d_j}\right| \lesssim \frac{n^{-1} + n^{-\epsilon_W}}{d_{i_m}}.$
Also, $1- O(n^{-1} + n^{-\epsilon_W}) \leq \sum_{j=1}^n s_{ij} \leq 1+ O(n^{-1} + n^{-\epsilon_W})$. Hence, $d_{i_m} = \sum_{j=1}^n\frac{s_{i_m j}}{d_j}\geq \frac{1-O(n^{-1} + n^{-\epsilon_W})}{d_{i_M}}$ and $ d_{i_M} = \sum_{j=1}^n\frac{s_{i_M j}}{d_j}\leq \frac{1+O(n^{-1} + n^{-\epsilon_W})}{d_{i_m}}$,
which implies $d_{i_M} d_{i_m} = 1 + O(n^{-1} + n^{-\epsilon_W})$. Also, $d_{i_m}(d_{i_M}-d_{i_m}) \lesssim n^{-1} + n^{-\epsilon_W}$.
Then, $d_{i_m} = (d_{i_M} d_{i_m} - d_{i_m}(d_{i_M}-d_{i_m}))^{1/2} = 1+ O(n^{-1} + n^{-\epsilon_W})$ and hence $d_{i_M} = 1+ O(n^{-1} + n^{-\epsilon_W})$. The claim follows.
\end{proof}

\begin{lem}
    \label{lem:W-reg}
    Let $L\geq 1$ be fixed, $\theta > 1$, and let $\bW = (W^{(a)})_{a=1}^L$ be a weakly Wigner tuple with parameters $(\xi_W,\epsilon_W,C_W)$. Let $v^{(a)} \in \SS^{n - 1}(\FF_a)$ be deterministic for each $a\in[L]$.

    For each $a\in[L]$, define
    \begin{align*}
        s^{(a)}_{ij} &\colonequals \EE |W^{(a)}_{ij}|^2, \\
        S^{(a)} &\colonequals (s^{(a)}_{ij})_{i,j=1}^n,
    \end{align*}
    Let $D^{(a)}\colonequals \Diag(d^{(a)}_1,\ldots, d^{(a)}_n)$ and
    write
    \[ W^{(\reg,a)} \colonequals (D^{(a)})^{-1/2} W^{(a)} (D^{(a)})^{-1/2}, \qquad \bW^{\reg} \colonequals \left(W^{(\reg,a)}\right)_{a=1}^L \]
    Then, for each $a\in[L]$, $W^{(\reg,a)}$ is a generalized Wigner matrix satisfying \eqref{eq:gen-wigner-var-sums} with uniformly controlled parameters. Moreover, $\bW^{(\reg)}$ is again a weakly Wigner tuple with parameters depending only on $(\xi_W,\epsilon_W,C_W)$ and $L$.
    Also, we have, for any $\epsilon > 0$,
    \begin{align*}
    &\bigg| \EE \Psi \left( \left(nv^{(a)}v^{(a)*}\right)_{a=1}^L, \left(n\what{v}^{(W,a)}\what{v}^{(W,a)*}\right)_{a=1}^L \right) \\ &\hspace{2cm} - \EE \Psi \left( \left(nv^{(a)}v^{(a)*}\right)_{a=1}^L, \left(n\what{v}^{(W^{(\reg)},a)}\what{v}^{(W^{(\reg)},a)*}\right)_{a=1}^L \right) \bigg| \lesssim n^{-1+\epsilon} + n^{-\epsilon_W + \epsilon} \numberthis \label{eq:W-reg-tuple},
    \end{align*}
    where
    \[\what{v}^{(W^{(\reg)},a)}\colonequals v_1\left(\theta v^{(a)} v^{(a)*} + W^{(\reg,a)}\right).\]
\end{lem}
\begin{proof}
    For each $a\in[L]$, by the weakly Wigner tuple assumption, we have $s^{(a)}_{ij} = \EE|W^{(a)}_{ij}|^2 = \frac{1}{n} + O(n^{-1-\epsilon_W}) > 0$ for $i\neq j$ and $s^{(a)}_{ii} = O(n^{-1})$ uniformly in $a$. Then \[\max_{i\in[n]}\sum_{j=1}^n\left|s^{(a)}_{ij} - \frac{1}{n}\right| \leq \max_{i\in[n]}\sum_{j\neq i} \left|s^{(a)}_{ij} -\frac{1}{n}\right| + \max_{i\in[n]} \left|s^{(a)}_{ii} -\frac{1}{n}\right| \lesssim n^{-1} + n^{-\epsilon_W}.\]
    Apply Lemma~\ref{lem:sinkhorn-near-constant} to $S^{(a)}$ for each $a\in[L]$, we then have \[\sum_{j=1}^n\frac{s^{(a)}_{ij}}{d^{(a)}_i d^{(a)}_j} = 1 \qquad\text{for all }i\in [n].\]
    Choose $D^{(a)}\colonequals \Diag(d^{(a)}_1,\ldots, d^{(a)}_n)$ so that $W^{(\reg,a)}$ satisfies \eqref{eq:gen-wigner-var-sums}.

    By Lemma~\ref{lem:sinkhorn-near-constant}, we have $\max_{a\in[L]}\max_{i\in[n]}\big|1 - d^{(a)}_i\big| \lesssim n^{-1}+n^{-\epsilon_W}$, and it then follows from elementary bounds that $\bW^{(\reg)}$ remains a weakly Wigner tuple as well.

    For the final bound, note that we have by the above observation $\|I-D^{(a)}\|+\|I-(D^{(a)})^{-1/2}\|\lesssim n^{-1}+n^{-\epsilon_W}$, and Proposition~\ref{prop:wig-norm} gives $\|(D^{(a)})^{-1/2} W^{(a)} (D^{(a)})^{-1/2}\|\prec 1$.
    So, we may bound
    \begin{align*}
        \|W^{(a)} - W^{(\reg,a)}\|
        &= \|W^{(a)} - (D^{(a)})^{-1/2}W^{(a)} (D^{(a)})^{-1/2}\| \\
        &\leq \|W^{(a)}(I - (D^{(a)})^{-1/2}) + (I - (D^{(a)})^{-1/2})W^{(a)}(D^{(a)})^{-1/2}\| \\
        &\leq \|W^{(a)}\| \cdot \|I - (D^{(a)})^{-1/2}\| \cdot \big(1 + \|(D^{(a)})^{-1/2}\|\big) \\
        &\prec n^{-1} + n^{-\epsilon_W}.
    \end{align*}
    The result then follows from Corollary~\ref{cor:dk}, with Theorem~\ref{thm:spiked-non-asymp} used as mentioned after the Corollary.
\end{proof}

Next, we show that we may, at the cost of a small error, replace $\bW^{(\reg)}$ in the above quantity with a Gaussian tuple.
\begin{lem}
    \label{lem:psi-univ}
    Let $L\geq 1$ be fixed, $\theta>1$, and let $v^{(a)} \in \SS^{n - 1}(\FF_a)$ deterministic satisfying Assumption~\ref{asm:v} uniformly for all $a\in[L]$.
    Let $\bW^{(\reg)} = \left(W^{(\reg,a)}\right)_{a=1}^L$ be a weakly Wigner tuple such that each coordinate
    $W^{(\reg,a)}$ is a generalized Wigner matrix with uniformly controlled parameters and satisfies the exact variance normalization \eqref{eq:gen-wigner-var-sums}.
    Let $\bX^{(0)} = \left(X^{(0,a)}\right)_{a=1}^L$ be a centered Gaussian tuple with the same first two moments as $\bW^{(\reg)}$ such that for each $1\leq p \leq q \leq n$, the real-coordinate block
    \[x^{(0,pq)} = \left(\Re X^{(0,1)}_{pq}, \Im X^{(0,1)}_{pq},\ldots, \Re X^{(0,L)}_{pq}, \Im X^{(0,L)}_{pq} \right)\] is centered Gaussian, the blocks are independent over $1\leq p \leq q \leq n$, and \[\EE x^{(0,pq)} x^{(0,pq)\top} = \EE w^{(\reg,pq)} w^{(\reg,pq)\top}, \] where $w^{(\reg,pq)}$ is the corresponding real-coordinate block of $\bW^{(\reg)}$.
    Write $H^{(Y,a)}\colonequals \theta v^{(a)} v^{(a)*} + Y^{(a)}$ for $Y^{(a)}\in\{W^{(\reg,a)}, X^{(0,a)}\}$ and $\what{v}^{(Y,a)} \colonequals v_1(H^{(Y,a)})$. Let $\psi, \Psi$ be as above. Then, for any $\epsilon > 0$,
    \begin{align*}
    &\bigg| \EE \Psi\left( \left(n v^{(a)}v^{(a)*}\right)_{a=1}^L, \left(n \what{v}^{(W^{(\reg)},a)}\what{v}^{(W^{(\reg)},a)*}\right)_{a=1}^L\right)\\
    &\hspace{2cm} - \EE \Psi\left( \left(n v^{(a)}v^{(a)*}\right)_{a=1}^L, \left(n \what{v}^{(X^{(0)},a)}\what{v}^{(X^{(0)},a)*}\right)_{a=1}^L\right)\bigg| \lesssim n^{-1/2 + 10\epsilon_v + \epsilon}.
    \end{align*}
\end{lem}
\begin{proof}
    This is a simple consequence of Theorem~\ref{thm:main-univ}: note that $\bW^{(\reg)}$ and $\bX^{(0)}$ satisfy the assumptions of that result, and we may bound
    \begin{align*}
    &\bigg| \EE \Psi\left( \left(n v^{(a)}v^{(a)*}\right)_{a=1}^L, \left(n \what{v}^{(W^{(\reg)},a)}\what{v}^{(W^{(\reg)},a)*}\right)_{a=1}^L\right)\\
    &\hspace{3.5cm} - \EE \Psi\left( \left(n v^{(a)}v^{(a)*}\right)_{a=1}^L, \left(n \what{v}^{(X^{(0)},a)}\what{v}^{(X^{(0)},a)*}\right)_{a=1}^L\right)\bigg|\\
    &\leq \frac{1}{n^2}\sum_{i,j=1}^n \Bigg|\EE\psi\left(\left(nv^{(a)}_i \overline{v^{(a)}_j}\right)_{a=1}^L, \left(n \what{v}^{(W^{(\reg)},a)}_i \overline{\what{v}^{(W^{(\reg)},a)}_j}\right)_{a=1}^L\right) \\&\hspace{3.5cm}- \EE \psi \left(\left(n v^{(a)}_i\overline{v^{(a)}_j}\right)_{a=1}^L, \left(n\what{v}^{(X^{(0)},a)}_i \overline{\what{v}^{(X^{(0)},a)}_j}\right)_{a=1}^L \right) \Bigg|\\
    &\lesssim \frac{1}{n^2} \sum_{i,j=1}^n n^{-1/2+10\epsilon_v+\epsilon} \lesssim n^{-1/2+10\epsilon_v+\epsilon}
    \end{align*}
    as claimed.
\end{proof}

Note that above $\bX^{(0)}$ is itself a weakly Wigner tuple, and each $X^{(0,a)}$ is a generalized Wigner matrix, since these conditions depend only on the first two moments and Gaussian entries have the required tail bounds.
We next show that the top eigenvectors of any such $\bX^{(0)}$ are close to those of a Gaussian tuple
$\bX = (X^{(1)},\ldots, X^{(L)}) \sim \bigotimes_{a=1}^L \GFE_a(n)$. The simple idea is that the weakly Wigner tuple assumption implies that we may couple $\bX^{(0)}$ and $\bX$ so that $\max_{a\in[L]}\|X^{(0,a)}-X^{(a)}\|\prec n^{-1/2} + n^{-\epsilon_W/2}$, whereby the top eigenvectors of outlier eigenvalues of rank-one perturbations of these matrices are close for all $a\in[L]$ by standard eigenvector perturbation inequalities.

\begin{lem}
    \label{lem:gauss-to-gfe}
    In the setting of Lemma~\ref{lem:psi-univ}, let $\bX = (X^{(1)},\ldots, X^{(L)}) \sim \bigotimes_{a=1}^L \GFE_a(n)$.
    Then, there exists a coupling between $\bX^{(0)}$ and $\bX$ such that
    \begin{equation}
    \sum_{a=1}^L\left\| \what{v}^{(X^{(0)},a)}\what{v}^{(X^{(0)},a)*} - \what{v}^{(X,a)}\what{v}^{(X,a)*} \right\|_{\F} \prec n^{-1/2} + n^{-\epsilon_W / 2}
    \end{equation}
\end{lem}
\begin{proof}
    From the condition of being weakly Wigner tuple, we may couple $\bX^{(0)}$ with a tuple $\bX\sim \bigotimes_{a=1}^L\GFE_a(n)$ so that, for each $a\in[L]$, \[X^{(0,a)} - X^{(a)} = -\delta X^{(a)}_{\offdiag} + \Delta^{(a)}_{\offdiag} + \Delta^{(a)}_{\diag}, \] where $\delta = O(n^{-\epsilon_W})$, $\|X^{(a)}_{\offdiag}\| \prec 1$, $\Delta^{(a)}_{\diag}$ is diagonal with independent Gaussian diagonal entries having variance $O(n^{-1})$, and $\Delta^{(a)}_{\offdiag}$ is the Hermitian matrix formed from the $a$-th coordinate of the off-diagonal residual blocks.
    To see that this decomposition exists, with $\delta = O(n^{-\epsilon_W})$, for each off-diagonal block, $\Sigma^{(pq)}_{0}-(1-\delta)^2\Sigma^{(pq)}_{\GFE}=(\Sigma^{(pq)}_{0}-\Sigma^{(pq)}_{\GFE}) + (1-(1-\delta)^2)\Sigma^{(pq)}_{\GFE}\succeq 0$, and this residual covariance is $O(n^{-1-\epsilon_W})$.

    We next bound $\Delta^{(a)}_{\offdiag}$. From standard bounds such as those of Theorem 1.1 and Corollary 3.9 of \cite{10.1214/15-AOP1025}, it follows that $\|\Delta^{(a)}_{\offdiag}\| \prec n^{-\epsilon_W/2}$. Also, $\|\Delta^{(a)}_{\diag}\|\prec n^{-1/2}$.
    Thus, using a union bound over fixed $L$, we have $\max_{a\in[L]}\|X^{(0,a)} - X^{(a)}\|\prec n^{-1/2} + n^{-\epsilon_W/2}$.

    By Theorem~\ref{thm:spiked-non-asymp}, again with a union bound over fixed $L$, with polynomially high probability, the gap between the top two eigenvalues of $\theta v^{(a)}v^{(a)*} + X^{(a)}$ is of size $\Omega(1)$.
    The result then follows by Proposition~\ref{prop:dk}.
\end{proof}

The third and more substantial part of the proof is a direct analysis of such expressions for $\GFE$ tuples.
For the sake of brevity below, for each $a\in[L]$, let us define
\[ y^{(a)} = y(\theta, v^{(a)}, g^{(a)}) \colonequals \rho(\theta) v^{(a)} +\frac{\tau(\theta)}{\sqrt{n}} g^{(a)},\] where $g^{(a)}\sim\sN_{\FF_a}(0,I_n)$ are independent over $a\in[L]$ and $\by\colonequals (y^{(a)})_{a=1}^L$.
\begin{lem}
    \label{lem:gauss-coupling}
    In the setting of Lemma~\ref{lem:gauss-to-gfe}, let $\bg =(g^{(1)},\ldots, g^{(L)})$ where $g^{(a)} \sim \sN_{\FF_a}(0, I_n)$ are independent over $a\in[L]$ and have i.i.d.\ entries.
    Then there exists a coupling between $\bX$ and $\bg$ such that
    \begin{equation}
    \sum_{a=1}^L \left\|\what{v}^{(X,a)}\what{v}^{(X,a)*} - y^{(a)}y^{(a)*} \right\|_{\F} \prec n^{-1/2}.
    \end{equation}
\end{lem}
\noindent
This is our main technical tool, which we prove in the following section.
For now, let us see how this together with Lemma~\ref{lem:psi-univ} completes the proof of Theorem~\ref{thm:main-gauss}.

\begin{cor}
    \label{cor:gauss-avg}
    In the setting of Lemma~\ref{lem:gauss-coupling}, we have that, for any $\epsilon > 0$,
    \begin{align*}
    &\Bigg|\EE \Psi\left(\left(n v^{(a)} v^{(a)*}\right)_{a=1}^L, \left(n \what{v}^{(X^{(0)},a)} \what{v}^{(X^{(0)},a)*}\right)_{a=1}^L\right) -\EE \Psi\left(\left(n v^{(a)}v^{(a)*}\right)_{a=1}^L, \left(ny^{(a)}y^{(a)*}\right)_{a=1}^L\right) \Bigg|\\ &\hspace{3cm}\lesssim n^{-1/2+\epsilon} + n^{-\epsilon_W/2+\epsilon}.
    \end{align*}
\end{cor}
\begin{proof}
    Combining Lemma~\ref{lem:gauss-to-gfe} and Lemma~\ref{lem:gauss-coupling}, we may realize $\bX^{(0)}$, $\bX$ and $\bg$ on the same probability space so that
    \[ \sum_{a=1}^L \left\|\what{v}^{(X^{(0)},a)}\what{v}^{(X^{(0)},a)*} - y^{(a)}y^{(a)*}\right\|_{\F} \prec n^{-1/2} + n^{-\epsilon_W/2}. \]
    Fix $\epsilon > 0$, and define the event
    \[ \sE = \left\{\sum_{a=1}^L \left\|\what{v}^{(X^{(0)},a)}\what{v}^{(X^{(0)},a)*} - y^{(a)}y^{(a)*}\right\|_{\F} \leq n^{-1/2 + \epsilon} + n^{-\epsilon_W/2 + \epsilon}\right\}. \]
    Then, by polynomial stochastic domination, for every fixed $D > 0$ we have $\PP[\sE^c] \leq n^{-D}$ for all sufficiently large $n$.
    Set
    \begin{align*}
        A &\colonequals \Psi\left(\left(n v^{(a)} v^{(a)*}\right)_{a=1}^L, \left(n \what{v}^{(X^{(0)},a)} \what{v}^{(X^{(0)},a)*}\right)_{a=1}^L\right), \\
        B &\colonequals \Psi\left(\left(n v^{(a)}v^{(a)*}\right)_{a=1}^L, \left(ny^{(a)}y^{(a)*}\right)_{a=1}^L\right).
    \end{align*}
    Since $\psi$ is bounded, $\Psi$ is bounded. Then by Proposition~\ref{prop:Psi-lip},
    on $\sE$ we have
    \[
    |A-B|\leq C_{\psi} \sum_{a=1}^L\left\| \what{v}^{(X^{(0)},a)}\what{v}^{(X^{(0)},a)*} - y^{(a)}y^{(a)*}\right\|_{\F} \leq C_{\psi} \left(n^{-1/2 + \epsilon} + n^{-\epsilon_W/2 + \epsilon}\right),\]
    where the factor $n$ in the inputs cancels by the $1/n$ Lipschitz factor
    from Proposition~\ref{prop:Psi-lip}.
    Using the boundedness of $\Psi$ on $\sE^c$,
    \[|\EE A - \EE B| \leq \EE\boldone_{\sE}|A-B| + \EE\boldone_{\sE^c}|A-B| \lesssim n^{-1/2+\epsilon} + n^{-\epsilon_W/2 + \epsilon} +n^{-D}. \]
    Taking $D$ large enough gives the claim.
\end{proof}

Theorem~\ref{thm:main-gauss} then follows from combining Corollary~\ref{cor:gauss-avg} with Lemma~\ref{lem:W-reg} and Lemma~\ref{lem:psi-univ}.
It remains to prove Lemma~\ref{lem:gauss-coupling}, which we do below.

\subsection{Analysis of Gaussian noise models: Proof of Lemma~\ref{lem:gauss-coupling}}
\label{sec:gauss-models}

It will also be useful to define
\[ \sU_{\FF}(n) = \left\{\begin{array}{ll} \sO(n), \text{ the $n \times n$ orthogonal group} & \text{if } \FF = \RR, \\ \sU(n), \text{ the $n \times n$ unitary group} & \text{if } \FF = \CC\end{array}\right\}. \]
The following is the main property of $\GFE$ matrices that we will use.
\begin{prop}
    \label{prop:invar}
    For each $\FF \in \{\RR, \CC\}$, the following hold:
    \begin{itemize}
        \item $g \sim \sN_{\FF}(0, I_n)$ has a $\sU_{\FF}(n)$ invariant law as a vector; that is, for all $U \in \sU_{\FF}(n)$, $Ug \eqlaw g$.
        \item $X \sim \GFE(n)$ has a $\sU_{\FF}(n)$-invariant law as Hermitian matrix; that is, for all $U \in \sU_{\FF}(n)$, $UXU^* \eqlaw X$.
    \end{itemize}
\end{prop}

\begin{proof}[Proof of Lemma~\ref{lem:gauss-coupling}]
    It suffices to prove the desired estimate for each fixed $a\in[L]$, uniformly in $a$, since $L$ is fixed.
    Fix $a\in[L]$ and suppress the superscript $a$ throughout this proof, writing
    \[
    \FF=\FF_a,\qquad v=v^{(a)},\qquad X=X^{(a)},\qquad g=g^{(a)}.
    \]
    Let us write $\what{v}(X) \colonequals v_1(\theta vv^* + X)$.
    For $U \in \sU_{\FF}(n)$, we have by Proposition~\ref{prop:invar} the equality of distributions
    \[ \what{v}(X) \eqlaw \what{v}(UXU^*) = v_1(\theta vv^* + UXU^*) = Uv_1(\theta (U^*v)(U^*v)^* + X). \]
    We may choose such a deterministic $U$ with $U^*v = e_1$ (note here the delicate point that we are using that $v \in \FF^n$, avoiding the case where $\FF = \RR$ while $v$ has complex values), so we find that there is some such (deterministic) $U$ for which
    \[ \what{v}(X) \eqlaw Uv_1(\theta e_1e_1^* + X). \]
    Define
    \[ \wtilde{y} \colonequals \rho(\theta) e_1 + \frac{\tau(\theta)}{\sqrt{n}} g. \]
    By Proposition~\ref{prop:invar} again, we have
    \[ \wtilde{y} \eqlaw U^* y. \]
    Let us couple the two sides of this equality so that, on the same probability space, we have $\wtilde{y} = U^* y$.
    Then, we have
    \begin{align*}
        \left\| v_1(\theta e_1e_1^* + X) v_1(\theta e_1e_1^* + X)^* - \wtilde{y}\wtilde{y}^*\right\|_{\F}
        &= \left\| Uv_1(\theta e_1e_1^* + X) v_1(\theta e_1e_1^* + X)^*U^* - U\wtilde{y}\wtilde{y}^*U^*\right\|_{\F} \\
        &\eqlaw \left\|\what{v}(X) \what{v}(X)^* - yy^*\right\|_{\F}
    \end{align*}
    for a corresponding coupling of $X$ and $g$.
    Thus, in effect we may take $v = e_1$ without loss of generality, so let us make this assumption going forward, in which case we may also take $y = \wtilde{y}$.

    We will only consider a single distribution of $X$ now, so let us take $X \sim \GFE(n)$ and $\what{v} \colonequals \what{v}(X) = v_1(\theta e_1e_1^* + X)$.
    We remove the sign (for $\FF = \RR$) or phase (for $\FF = \CC$) ambiguity by choosing $\what{v}$ to have its first coordinate real and non-negative.
    Write $\what{v}^{\perp}$ for $\what{v}$ with the first entry set to zero, so that we have
    \begin{equation}
    \label{eq:gauss-whatv-decomp}
    \what{v} = |\langle \what{v}, e_1 \rangle| e_1 + \what{v}^{\perp}.
    \end{equation}

    Let $\sU_{\FF}^{(1)}(n)$ be the subgroup of matrices in $\sU_{\FF}(n)$ that fix $e_1$, equivalently matrices of the form $\Diag(1,U_0)$ with $U_0\in\sU_{\FF}(n-1)$.
    Then, $\theta e_1e_1^* + X$ is $\sU_{\FF}^{(1)}(n)$-invariant, in the same sense as in the second part of Proposition~\ref{prop:invar}, since the first summand is fixed by the action of these matrices and the distribution of the second summand is unchanged.
    Thus, its top eigenvector $\what{v}$ is $\sU_{\FF}^{(1)}(n)$-invariant, in the same sense as in the first part of Proposition~\ref{prop:invar}.
    Since the two summands in \eqref{eq:gauss-whatv-decomp} each belong to fixed subspaces of this group action, they are each individually invariant as well.
    In particular, $\what{v}^{\perp} / \|\what{v}^{\perp}\|$ has the law of an $(n - 1)$-dimensional random vector drawn uniformly at random from $\SS^{n - 2}(\FF)$, with a zero entry prepended to it.

    Let us introduce $g \sim \sN_{\FF}(0, I_n)$ and write $g^{\perp}$ for $g$ with the first entry set to zero.
    By the above observations, we have the equality of laws
    \[ \frac{\what{v}^{\perp}}{\|\what{v}^{\perp}\|} \eqlaw \frac{g^{\perp}}{\|g^{\perp}\|}. \]
    Let us couple $X$ and $g$ such that this equality holds.
    From this $g$, we set
    \[ y = \rho(\theta)e_1 + \frac{\tau(\theta)}{\sqrt{n}}g. \]
    Then, we have
    \begin{align*}
        &\hspace{-0.125cm}\left\|\what{v}\what{v}^* - yy^*\right\|_{\F} \\
        &= \left\|(\what{v} - y)\what{v}^* + y(\what{v} - y)^*\right\|_{\F} \\
        &\leq (1 + \|y\|)\|\what{v} - y\| \\
        &\leq \left(1 + \rho(\theta) + \tau(\theta) \frac{\|g\|}{\sqrt{n}}\right)\left(\big||\langle \what{v}, e_1 \rangle| - \rho(\theta)\big| + \left\|\frac{\tau(\theta)}{\sqrt{n}}g - \frac{\|\what{v}^{\perp}\|}{\|g^{\perp}\|} g^{\perp}\right\|\right) \\
        &\leq \left(1 + \rho(\theta) + \tau(\theta) \frac{\|g\|}{\sqrt{n}}\right)\left(\big||\langle \what{v}, e_1 \rangle| - \rho(\theta)\big| + \tau(\theta)\frac{|\langle g,e_1\rangle|}{\sqrt{n}} + \left|\frac{\tau(\theta)}{\sqrt{n}} - \frac{\|\what{v}^{\perp}\|}{\|g^{\perp}\|}\right|\|g^{\perp}\|\right) \\
        &= \left(1 + \rho(\theta) + \tau(\theta) \frac{\|g\|}{\sqrt{n}}\right)\left(\big||\langle \what{v}, e_1 \rangle| - \rho(\theta)\big| + \tau(\theta)\frac{|\langle g,e_1\rangle|}{\sqrt{n}} + \left|\tau(\theta) \frac{\|g^{\perp}\|}{\sqrt{n}} - \|\what{v}^{\perp}\|\right|\right) \\
        &\leq \left(1 + \rho(\theta) + \tau(\theta) \frac{\|g\|}{\sqrt{n}}\right)\left(\big||\langle \what{v}, e_1 \rangle| - \rho(\theta)\big| + \tau(\theta)\frac{|\langle g, e_1 \rangle|}{\sqrt{n}} + \left| \tau(\theta) - \|\what{v}^{\perp}\|\right|\frac{\|g^{\perp}\|}{\sqrt{n}} + \left|\frac{\|g^{\perp}\|}{\sqrt{n}} - 1\right|\right).
    \end{align*}
    We recall here that $\tau(\theta) = \sqrt{1 - \rho(\theta)^2}$ and by definition $\|\what{v}^{\perp}\| = \sqrt{1 - |\langle \what{v}, e_1\rangle|^2}$, so quantities in the third term of the second factor are determined by those in the first term.

    To complete the proof, it suffices to show that the above expression is small with high probability.
    To that end, fix $\epsilon > 0$ and define the events
    \begin{align*}
        \sE_1 &\colonequals \left\{\left|\frac{\|g\|}{\sqrt{n}}-1\right| \leq n^{-1/2+\epsilon},\quad \left|\frac{\|g^{\perp}\|}{\sqrt{n}}-1\right| \leq n^{-1/2+\epsilon}\right\}, \\
        \sE_2 &\colonequals \left\{|\langle g, e_1\rangle| \leq \log n\right\}, \\
        \sE_3 &\colonequals \left\{\big||\langle \what{v}, e_1 \rangle| - \rho(\theta)\big| \leq n^{-1/2 + \epsilon}\right\}, \\
        \sE &\colonequals \sE_1 \cap \sE_2 \cap \sE_3.
    \end{align*}
    We see from elementary manipulations that, on the event $\sE$, we have
    \[ \|\what{v}\what{v}^* - yy^*\|_{\F} \lesssim n^{-1/2 + \epsilon}. \]
    Fix some $K > 0$.
    We have by the tail bounds proved in \cite{10.1214/aos/1015957395} that
    \[ \PP[\sE_1^c] \lesssim n^{-K}, \]
    and by well-known Gaussian tail bounds we also have
    \[ \PP[\sE_2^c] \lesssim n^{-K}. \]
    Finally, by Theorem~\ref{thm:spiked-non-asymp}, we also have $\left||\langle \what{v}, e_1 \rangle|^2 - \rho(\theta)^2\right| \prec n^{-1/2}$, and thus
    \[ \PP[\sE_3^c] \lesssim n^{-K}. \]
    Thus $\PP[\sE^c] \lesssim n^{-K}$. Hence, for this fixed $a$, \[ \left\|\what{v}^{(X,a)}\what{v}^{(X,a)*}-y^{(a)}y^{(a)*}\right\|_{\F}\prec n^{-1/2}.\]
    After constructing the coupling for each $a$, we take these couplings independently over $a\in[L]$. Since $L$ is fixed, the sum of the resulting $L$ errors is still $\prec n^{-1/2}$.
\end{proof}

\subsection{Numerical evaluation of Gaussian fluctuations}

\begin{figure}
    \centering
    \includegraphics[width=1\linewidth]{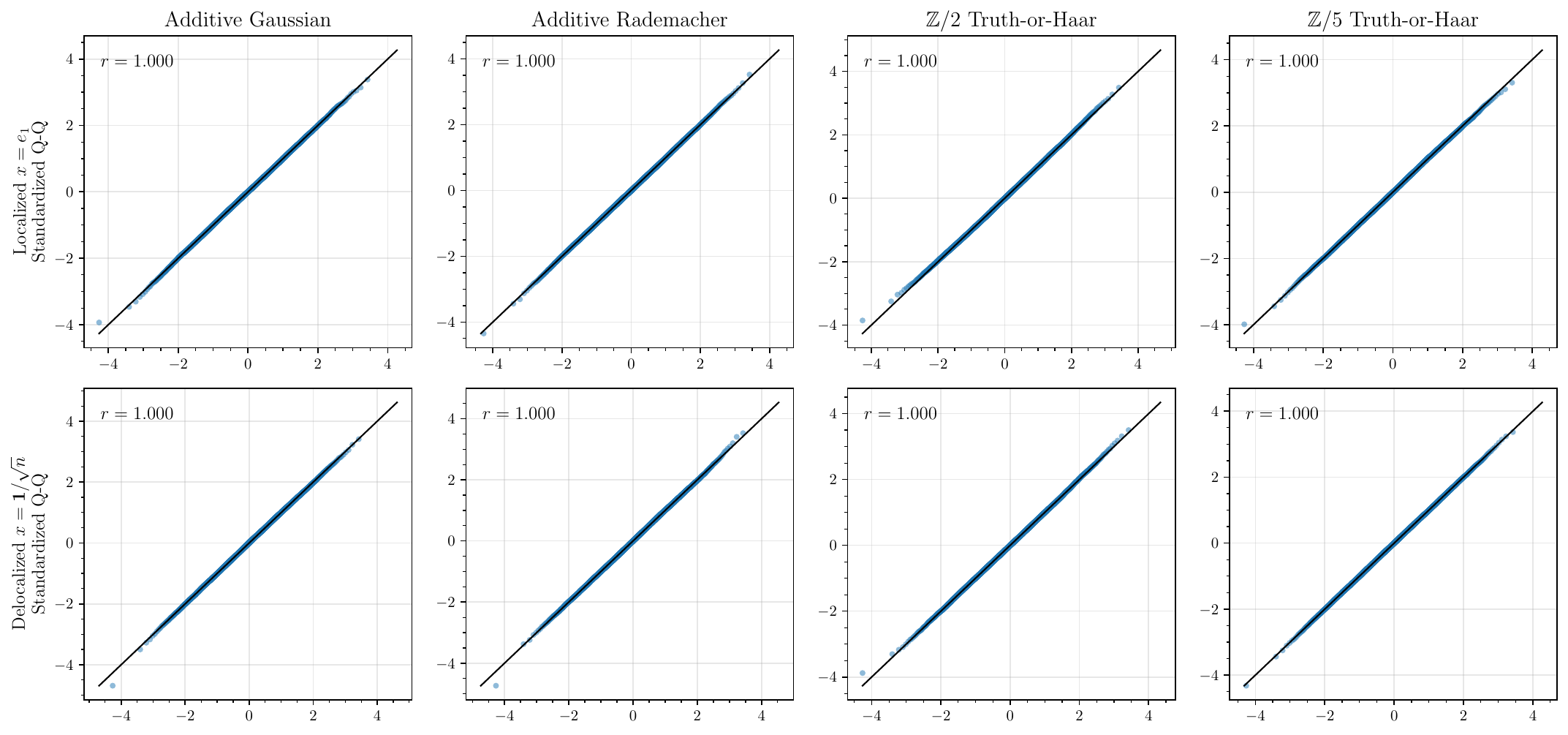}
    \caption{\textbf{Q--Q plots of one-dimensional projections of eigenvector fluctuations.} We present Q--Q plots justifying the approximate Gaussianity of $\what{v} - \rho(\theta)v$ proposed in \eqref{eq:gauss-approx}. In particular, we present such plots for the inner product of this vector with two test directions $x$, one localized and one delocalized, and four different noise models. We find strong evidence of Gaussianity in all cases.}
    \label{fig:qq}
\end{figure}

As a simple test of the Gaussianity given as an informal approximation in \eqref{eq:gauss-approx} and made precise in Theorem~\ref{thm:main-gauss}, in Figure~\ref{fig:qq} we present Gaussian Q--Q plots of eigenvector fluctuations projected in a delocalized and in a localized direction.
Both have a distribution very close to Gaussian, and this persists under several different noise models: we consider additive Gaussian and Rademacher noise (as in Figure~\ref{fig:non-univ}) as well as the truth-or-Haar model of group synchronization as detailed in Section~\ref{sec:intro-appl} for two cyclic groups.

\section{Applications}
\label{sec:appl}

\subsection{Angular synchronization: Proof of Theorem~\ref{thm:sync}}

We recall the construction of the random matrix tuple $H^{(1)},\ldots,H^{(L)}$ that this result concerns: we have an underlying family of group elements $x_1, \dots, x_n \stackrel{\mathrm{i.i.d.}}{\sim}\Haar(G)$, and define the group-valued matrix $M_{ij} = x_ix_j^{-1}$.
From this, we define $Y$ to be a noisy version of $M$ by \eqref{eq:sync-Y}, and $H^{(a)}$ the entrywise image of $Y$ under the character $\chi^{(a)}$ given in \eqref{eq:sync-chi}.
Let us write \[v^{(a)}_i \colonequals \frac{\chi^{(a)}(x_i)}{\sqrt{n}}, \qquad a\in[L].\]
We have $|v^{(a)}_i| = n^{-1/2}$, so $v^{(a)} \in \SS^{n - 1}(\FF_a)$, where $\FF_a = \RR$ if $\chi^{(a)}(G) \subseteq \RR$ and $\FF_a = \CC$ otherwise, $v^{(a)}$ has $\|v^{(a)}\|_{\infty} = n^{-1/2}$, thus satisfying Assumption~\ref{asm:v} for any $\epsilon_v\in(0,1/20)$.

Conditional on $x$, then, the upper-triangular off-diagonal blocks of $H^{(1)},\ldots,H^{(L)}$ are distributed as, independently for $1 \leq i < j \leq n$,
\begin{equation}
\label{eq:sync-H}
\Big(H^{(a)}_{ij}\Big)_{a=1}^L \sim \left\{\begin{array}{ll} \Big(\sqrt{n} \cdot v^{(a)}_i\overline{v^{(a)}_j}\Big)_{a=1}^L & \text{with probability } p, \\ \Big(\chi^{(a)}(y) / \sqrt{n}\Big)_{a=1}^L,  y\sim\Haar(G)
& \text{with probability } 1 - p\end{array}\right\},
\end{equation}
where $p = \theta / \sqrt{n}$.
Since we construct $Y$ to be $G$-Hermitian, we also have $H^{(a)} = H^{(a)*}$.
Under our specific definition from the Introduction we will always have $H^{(a)}_{ii} = 1 / \sqrt{n} = \chi^{(a)}(e) / \sqrt{n}$ for $e$ the identity of $G$, but this detail will be inconsequential as we have shown above.

We then have, conditional on the underlying randomness of the $x_i$, for the off-diagonal entries,
\[ \EE\big[H^{(a)}_{ij} \mid x\big] = p\sqrt{n} \cdot v^{(a)}_i \overline{v^{(a)}_j} = \theta \cdot v^{(a)}_i\overline{v^{(a)}_j}, \qquad i\neq j,\  a\in[L]. \]
We then define
\[ W^{(a)} \colonequals H^{(a)} - \theta v^{(a)}v^{(a)*}, \qquad \bW \colonequals \Big(W^{(a)}\Big)_{a=1}^L, \]
which will be centered conditional on $x$, away from the inconsequential diagonal entries.
Specifically, its entries conditional on $x$ are distributed as
\begin{equation}
\label{eq:sync-W}
\Big(W^{(a)}_{ij}\Big)_{a=1}^L \sim \left\{\begin{array}{ll} \Big((\sqrt{n} - \theta) \cdot v^{(a)}_i\overline{v^{(a)}_j}\Big)_{a=1}^L & \text{with probability } p, \\ \Big(\chi^{(a)}(y) / \sqrt{n} - \theta v^{(a)}_i\overline{v^{(a)}_j}\Big)_{a=1}^L, y\sim\Haar(G) & \text{with probability } 1 - p\end{array}\right\}.
\end{equation}
This is a centered Hermitian tuple whose upper-triangular blocks are independent conditional on $x$.
We have $\big|v^{(a)}_i\overline{v^{(a)}_j}\big| = 1 / n$, and thus the covariance matrix of the vector of real coordinates
\[w^{(ij)} \colonequals \left(\Re W^{(1)}_{ij}, \Im W^{(1)}_{ij},\dots,\Re W^{(L)}_{ij}, \Im W^{(L)}_{ij}\right)\]
is $n^{-1}$ times that of
\[\sX(y) \colonequals \left(\Re \chi^{(1)}(y), \Im \chi^{(1)}(y),\dots,\Re \chi^{(L)}(y), \Im \chi^{(L)}(y)\right), \qquad y\sim \Haar(G),\]
up to an additive error $O(n^{-3/2})$:
\[\EE \left[w^{(ij)} w^{(ij)\top} \mid x \right] = \frac{1}{n}\Ex_{y\sim\Haar(G)} \sX(y) \sX(y)^{\top} + O(n^{-3/2}).\]
Here and below, the $O(n^{-3/2})$ term denotes a $2L\times 2L$ matrix with entries $O(n^{-3/2})$.
Since $L$ is fixed, the same bound also holds in operator norm, up to changing the implicit constant.

We use the following group-theoretic observation, a slight variation on the classical orthogonality of characters, to verify the covariance assumption of the definition of weakly Wigner tuples.
We note that this applies equally well to one-dimensional irreducible representations of any compact Lie group, but we focus on the cases that will appear in our applications for the sake of simplicity.
\begin{prop}[Real-coordinate character orthogonality]\label{prop:character-orthogonality}
Let $G=\ZZ/K$ or $G=U(1)$, let $y\sim\Haar(G)$ and let $\chi^{(1)},\ldots,\chi^{(L)}: G \to U(1)$ be nontrivial characters such that $\chi^{(a)}\notin \{\chi^{(b)},\overline{\chi^{(b)}}\}$ for $a\neq b$. For $a\in[L]$, let $\FF_a=\RR$ if $\chi^{(a)}(G) \subseteq \RR$ and $\FF_a=\CC$ otherwise. Then
\[\EE\sX(y)=0,\qquad \EE\big[\sX(y) \sX(y)^{\top}\big] = \Diag(\Sigma_{\FF_1},\ldots,\Sigma_{\FF_L}),\] where
\[\Sigma_{\RR} = \begin{pmatrix} 1&0\\ 0&0 \end{pmatrix},\qquad \Sigma_{\CC} = \begin{pmatrix} 1/2&0\\ 0&1/2 \end{pmatrix}.\]
\end{prop}
\begin{proof}
Using translation invariance of normalized Haar measure, we have
\[ \Ex_{y\sim\Haar(G)}\chi(y) = \boldone\{\chi \text{ is trivial}\}. \]
By our assumption, each $\chi^{(a)}$ is nontrivial, so $\EE\sX(y)=0$.
For $a,b\in[L]$, set
\[A_{ab} \colonequals \EE\big[\chi^{(a)}(y)\overline{\chi^{(b)}(y)}\big],\qquad B_{ab} \colonequals \EE\big[\chi^{(a)}(y)\chi^{(b)}(y)\big],\]
then \[A_{ab} = \boldone\{\chi^{(a)} = \chi^{(b)}\}, \qquad B_{ab} = \boldone\{\chi^{(a)} = \overline{\chi^{(b)}}\}.\]
Then the $(a,b)$-block of $\EE\big[\sX(y) \sX(y)^{\top}\big]$ is
\begin{equation}\label{eq:covariance}
\EE \begin{pmatrix} \Re \chi^{(a)}(y)\\ \Im \chi^{(a)}(y) \end{pmatrix} \begin{pmatrix} \Re \chi^{(b)}(y)& \Im \chi^{(b)}(y)\end{pmatrix} = \frac{1}{2} \begin{pmatrix} A_{ab}+B_{ab} & 0 \\ 0 & A_{ab}-B_{ab} \end{pmatrix}.
\end{equation}
If $a\neq b$, the assumptions give $A_{ab}=B_{ab}=0$, so the corresponding off-diagonal block vanishes.
If $a=b$, then $A_{aa}=1$ while $B_{aa} = \boldone\{\chi^{(a)}=\overline{\chi^{(a)}}\} = \boldone\{\chi^{(a)}(G)\subseteq\RR\}$. Thus, the $(a,a)$-th block is $\Sigma_{\RR}$ when $\FF_a=\RR$, and $\Sigma_{\CC}$ when $\FF_a=\CC$.
\end{proof}
By Proposition \ref{prop:character-orthogonality}, uniformly for $1\leq i < j \leq n$,
\[\EE \left[w^{(ij)} w^{(ij)\top} \mid x \right] = \frac{1}{n}\Diag(\Sigma_{\FF_1},\ldots,\Sigma_{\FF_L}) + O(n^{-3/2}).\]
Replacing the diagonal of $W^{(a)}$ by zero only subtracts the scalar matrix $(n^{-1/2}-\theta/n)I$ from $H^{(a)}$, and hence does not change its eigenvectors. After this diagonal modification, conditional on $x$, $\bW$ is a weakly Wigner tuple, with $C_W$ an absolute constant and $\epsilon_W=1/2$.

Now, let $\what{v}^{(a)}$ be the top eigenvector of this $H^{(a)}$ for each $a\in[L]$.
Choosing the parameters $\epsilon_v>0$ in Assumption~\ref{asm:v} and $\epsilon>0$ in Theorem~\ref{thm:main-gauss} sufficiently small, we then obtain that, provided that $\psi$ is $\sC^5$ with bounded values and first five derivatives, conditionally on $x$ (whereby $v^{(a)}$ is not random), we have
\begin{align*}
&\Bigg| \EE \left[\Psi\left( \left( n\cdot v^{(a)}v^{(a)*}\right)_{a=1}^L, \left(n\cdot \what{v}^{(a)} \what{v}^{(a)*}\right)_{a=1}^L \right)\middle|x\right] \\ &- \Ex_{\substack{g^{(a)}\sim\sN_{\FF_a}(0,I_n)\\ a\in[L]}} \Psi \left( \left(n \cdot v^{(a)}v^{(a)*}\right)_{a=1}^L,
\left((\rho(\theta)\sqrt{n}\cdot v^{(a)}+\tau(\theta)g^{(a)}) (\rho(\theta)\sqrt{n}\cdot v^{(a)}+\tau(\theta)g^{(a)})^*
\right)_{a=1}^L
\right)\Bigg|\\ &\hspace{1cm}\lesssim n^{-1/5}.
\end{align*}
Next, we show that the expectation over $x$ (and thus over the randomness in $v^{(a)}$) of the second expression converges to a single-letter formula.
Indeed, expanding the definition of $\Psi$, we have
\begin{align*}
    &\Ex_{x_1, \dots, x_n \stackrel{\mathrm{i.i.d.}}{\sim}\Haar(G)}
    \Ex_{\substack{g^{(a)}\sim\sN_{\FF_a}(0,I_n)\\ a\in[L]}}
    \Psi\Bigg(
    \left(n \cdot v^{(a)}v^{(a)*}\right)_{a=1}^L, \\
    &\hspace{3cm}
    \left(
    (\rho(\theta)\sqrt{n}\cdot v^{(a)}+\tau(\theta)g^{(a)})
    (\rho(\theta)\sqrt{n}\cdot v^{(a)}+\tau(\theta)g^{(a)})^*
    \right)_{a=1}^L
    \Bigg) \intertext{ }
    &=
    \Ex_{x_1, \dots, x_n \stackrel{\mathrm{i.i.d.}}{\sim} \Haar(G)}
    \Ex_{\substack{g^{(a)}\sim\sN_{\FF_a}(0,I_n)\\ a\in[L]}}
    \frac{1}{n^2}\sum_{i,j=1}^n
    \psi\Bigg(
    \left(\chi^{(a)}(x_i)\overline{\chi^{(a)}(x_j)}\right)_{a=1}^L,\\
    &\hspace{3cm}
    \left(
    (\rho(\theta)\chi^{(a)}(x_i)+\tau(\theta)g_i^{(a)})
    \overline{(\rho(\theta)\chi^{(a)}(x_j)+\tau(\theta)g_j^{(a)})}
    \right)_{a=1}^L
    \Bigg) \\
    &=
    \frac{1}{n^2}\sum_{i\neq j}
    \Ex_{\substack{x_i,x_j\sim\Haar(G)\\
    g_i^{(a)},g_j^{(a)}\sim\sN_{\FF_a}(0,1),\ a\in[L]}}
    \psi\Bigg(
    \left(\chi^{(a)}(x_i)\overline{\chi^{(a)}(x_j)}\right)_{a=1}^L,\\
    &\hspace{3cm}
    \left(
    (\rho(\theta)\chi^{(a)}(x_i)+\tau(\theta)g_i^{(a)})
    \overline{(\rho(\theta)\chi^{(a)}(x_j)+\tau(\theta)g_j^{(a)})}
    \right)_{a=1}^L
    \Bigg)
    +O(n^{-1})\\
    &=
    \Ex_{\substack{x,y\sim\Haar(G)\\
    g_a,h_a\sim\sN_{\FF_a}(0,1),\ a\in[L]}}
    \psi\Bigg(
    \left(\chi^{(a)}(x)\overline{\chi^{(a)}(y)}\right)_{a=1}^L,\\
    &\hspace{3cm}
    \left(
    (\rho(\theta)\chi^{(a)}(x)+\tau(\theta)g_a)
    \overline{(\rho(\theta)\chi^{(a)}(y)+\tau(\theta)h_a)}
    \right)_{a=1}^L
    \Bigg)
    +O(n^{-1}),
\end{align*}
since the diagonal terms contribute negligibly and all other expectations are equal.

To derive the actual statement of Theorem~\ref{thm:sync}, we choose $\psi$ such that
\[\psi\left(\left(\chi^{(a)}(x) \overline{\chi^{(a)}(y)}\right)_{a=1}^L, \bz\right) = \ell(xy^{-1}, \Round(\bz)),\]
for $x, y \in G$ and $\bz \in \CC^L$.
Since $\left(\chi^{(a)}(x)\overline{\chi^{(a)}(y)}\right)_{a=1}^L=\bm{\chi}(xy^{-1})$ and $\bm{\chi}$ is injective, this is always possible.
With this choice,
\[ \Psi\left(\left(n \cdot v^{(a)}v^{(a)*}\right)_{a=1}^L,\left(n \cdot \what{v}^{(a)}\what{v}^{(a)*}\right)_{a=1}^L\right) = \frac{1}{n^2}\sum_{i,j=1}^n \ell(M_{ij},\what{M}_{ij}),
\]
and the single-letter expression above becomes
\[
\Ex_{\substack{x,y\sim\Haar(G)\\ g_a,h_a\sim\sN_{\FF_a}(0,1),\ a\in[L]}}\ell\bigg(xy^{-1},\Round\left(\left((\rho(\theta)\chi^{(a)}(x)+\tau(\theta)g_a)\overline{(\rho(\theta)\chi^{(a)}(y)+\tau(\theta)h_a)}\right)_{a=1}^L\right)\bigg).
\]
The resulting $\psi$ may not be smooth.
However, conditional on $x$ and $y$, the random vector
\[
\left((\rho(\theta)\chi^{(a)}(x)+\tau(\theta)g_a)\overline{(\rho(\theta)\chi^{(a)}(y)+\tau(\theta)h_a)}\right)_{a=1}^L
\]
has a density with respect to Lebesgue measure on $\FF_1\times\cdots\times\FF_L$.
It is therefore a continuity point of $\Round$ almost surely.
Theorem~\ref{thm:main-gauss} gives weak convergence against bounded smooth test functions, and the standard extension of weak convergence to bounded functions that are continuous almost surely under the limiting measure gives the result.

\subsection{Discussion and numerical experiments}
\label{sec:numerical}

Let us give some additional discussion and present numerical experiments verifying these results.
We first validate that the prediction for the average entrywise error of Theorem~\ref{thm:sync} is sound with some extra results extending those presented for the circle group $G = U(1)$ in Figure~\ref{fig:emp-pred-U} above.
Then we study the implications of this prediction for the choice of the parameters $L$ and $\varrho$ in a multi-frequency spectral algorithm.

\begin{figure}
    \centering
    \includegraphics[width=1\linewidth]{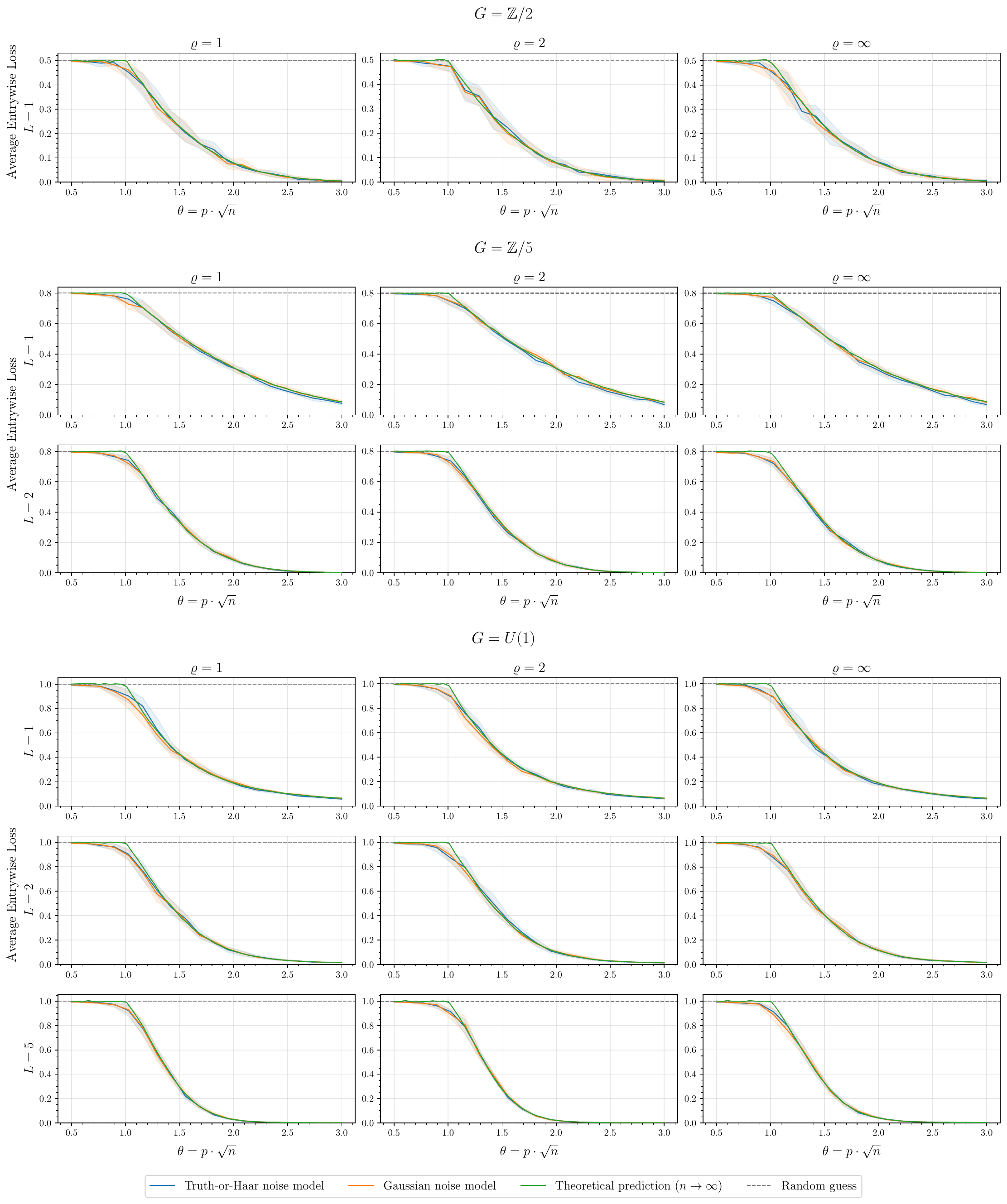}

    \caption{\textbf{Empirical versus predicted synchronization performance.} We plot the empirical average loss of the rounded multi-frequency spectral algorithm for group synchronization for several choices of group $G$, number of frequencies $L$, and rounding function parameter $\varrho$.}
    \label{fig:emp-pred-all}
\end{figure}

First, in Figure~\ref{fig:emp-pred-all}, we present experiments of the same kind as in Figure~\ref{fig:emp-pred-U} comparing the prediction of Theorem~\ref{thm:sync} with empirical distributions of average entrywise loss, now considering several groups $G$, several numbers of frequencies $L$, and several parameters $\varrho$ of the rounding function.
For $G = U(1)$ we use the loss function $\ell(x, y) = 1 - \cos(x - y)$ as before, while for $G = \ZZ/K$ (for $K \in \{2, 5\}$ in the experiments) we use the indicator loss function $\ell(x, y) = \boldone\{x \neq y\}$.
We note that, for the indicator loss function, there is a ``baseline'' fraction of errors $\frac{K - 1}{K}$ achieved by the trivial estimator of $M$ that guesses each entry uniformly at random, which is indeed what our $\what{M}$ achieves for $\theta < 1$ and which we plot with a dotted line labelled ``Random guess'' in the figures below.
Agreement with the theoretical prediction is strong in all cases.

\begin{figure}
    \centering
    \includegraphics[width=1\linewidth]{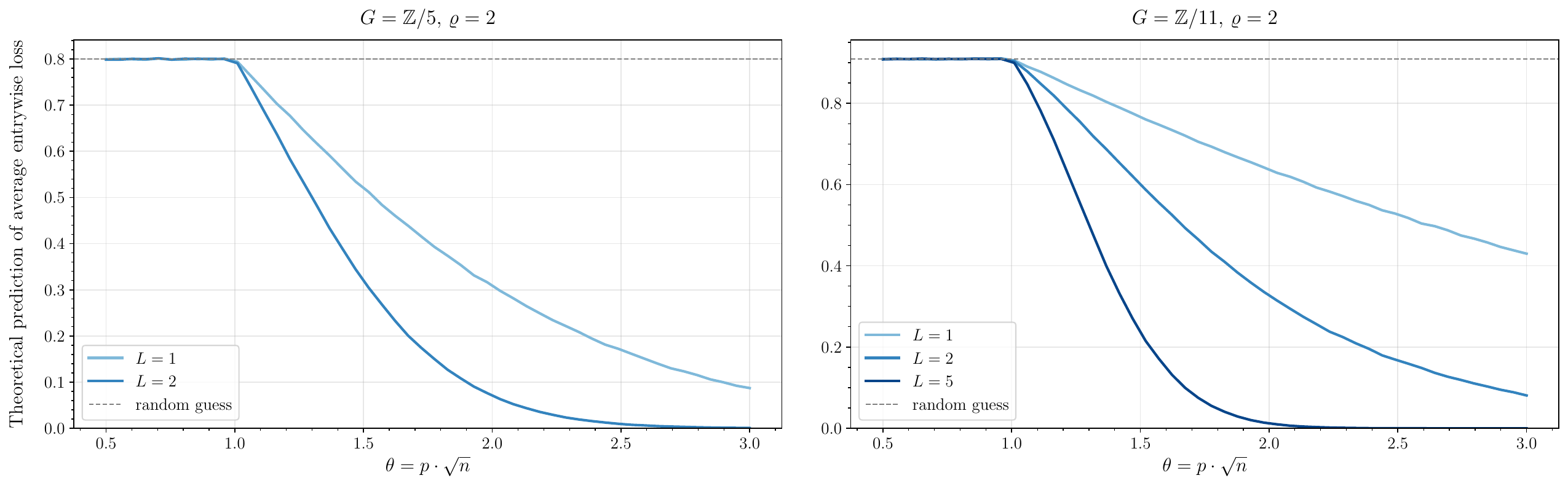}
    \caption{\textbf{Predicted superiority of multiple frequencies for cyclic groups.} We plot the theoretical predictions of average entrywise loss per Theorem~\ref{thm:sync} for group synchronization in two cyclic groups for various numbers of frequencies $L$.}
    \label{fig:multifreq-cyc}
\end{figure}

Next, taking as a given the accuracy of the prediction of Theorem~\ref{thm:sync} tested in these experiments, we consider how to best choose the values of the parameters $L$ and $\varrho$.

\subsubsection{Choice of \texorpdfstring{$L$}{L} parameter: Superiority of multi-frequency algorithms}

We showed the effect of different choices of $L$ on the formula of Theorem~\ref{thm:sync} in the case $G = U(1)$ in Figure~\ref{fig:multifreq-U}.
In Figure~\ref{fig:multifreq-cyc}, we show analogous results for two cyclic groups $\ZZ / K$ (note that for multi-frequency algorithms to be defined, we must take $K \geq 3$, so we do not include the case $K = 2$ here unlike in the previous experiments).
In both cases, at all values of $\theta$, algorithms using more frequencies and $\varrho = 2$ (see below for discussion of this choice) achieve lower loss.

We note that the phenomenon observed earlier for the case of $G = U(1)$ with the cosine loss function $\ell(x, y) = 1 - \cos(x - y)$ does not appear here: for cyclic groups with the indicator loss function, our calculations indicate that using more frequencies gives a uniformly superior algorithm for all values of signal strength $\theta$.
On the other hand, we \emph{do} observe (in experiments omitted here) the same non-monotonicity using the cosine loss function for cyclic groups, suggesting that the phenomenon is associated to the choice of loss function $\ell$ rather than to the underlying group $G$.

\subsubsection{Choice of \texorpdfstring{$\varrho$}{varrho} parameter: Superiority of Euclidean rounding}

\begin{figure}
    \centering
    \includegraphics[width=1\linewidth]{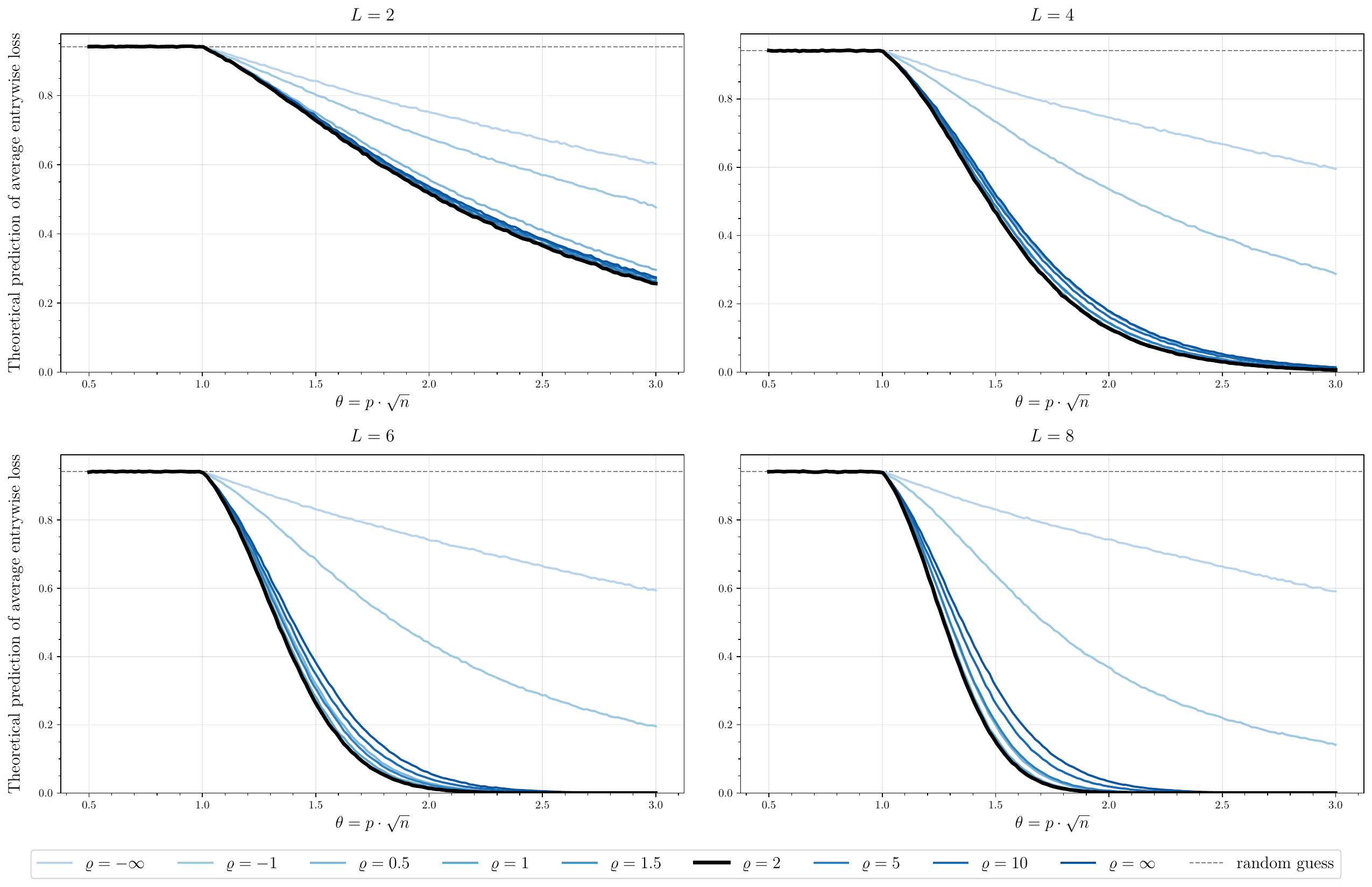}
    \caption{\textbf{Predicted performance of different rounding functions.} We plot the theoretical predictions of average entrywise loss per Theorem~\ref{thm:sync} for group synchronization in $G = \ZZ / 17$ for various numbers of frequencies $L$ and, in each plot, for various choices of the parameter $\varrho$ controlling the choice of rounding function. }
    \label{fig:rho}
\end{figure}

For the parameter $\varrho$, in Figure~\ref{fig:rho} we fix a cyclic group $G = \ZZ / K$, fix various numbers of frequencies $L$, and consider the average loss predicted by Theorem~\ref{thm:sync} over a range of values of $\varrho \in [-\infty, \infty]$.
Perhaps surprisingly, we find that in all cases, among the values of $\varrho$ we include, $\varrho = 2$ achieves the lowest average loss for all values of $\theta$, lower than choices of $\varrho$ both larger (more ``max-like'' loss functions) and smaller (more ``min-like'' loss functions) than $\varrho = 2$.

\begin{rmk}
    Another special property of the choice $\varrho = 2$ is that, if we consider including all frequencies associated to characters $\chi_1, \dots, \chi_K$, in the objective expression $\sum_{i = 1}^K |\chi_i(x) - z_i|^{\varrho}$ for some $x \in G$, when $\varrho = 2$ this is the same as the squared distance between the inverse discrete Fourier transform of the vector $z$ and the indicator vector associated to $x \in G$, by the $G$-valued version of Parseval's theorem.
    Still, it is unclear what this isometry might have to do with the superiority of the associated rounding scheme for a spectral algorithm.
\end{rmk}

\subsubsection{Sketch of extension to non-abelian groups}

We do not pursue it here, but we remark that it is in principle straightforward (though technically tedious) to extend these results to multi-frequency synchronization over non-abelian groups $G$.
Let us sketch how such an extension might look.

In this case, following for instance the ideas of \cite{Perry_2018}, we would like to apply possibly higher-dimensional irreducible representations entrywise to the matrix $Y$ to form a block matrix $H$; if the underlying representation is of dimension $d$, then $H$ will be $dn \times dn$.
By the Peter-Weyl theorem on orthogonality of matrix coefficients, such matrices will obey a covariance property analogous to that in Definition~\ref{def:weakly-Wigner-tuples}, and thus can be compared to spiked matrices with independent Gaussian noise via an analog of Theorem~\ref{thm:main-gauss}.
We expect that a version of Theorem~\ref{thm:sync} should therefore apply, where now if our $\chi_1, \dots, \chi_L$ were replaced with representations of dimensions $d_1, \dots, d_L$, then our rounding function would map from $\CC^{d_1 \times d_1} \times \cdots \times \CC^{d_L \times d_L}$ to $G$, and accordingly the appropriate version of Theorem~\ref{thm:sync} would involve expectations of this dimension.
One subtle point is that the higher-dimensional  matrix-valued replacements for the $g_a$ and $h_a$ in the statement (scalars in our setting) must be matched to the type (real, complex, or quaternionic) of the representation $\chi_a$.
The main contours of our proof technique should also still apply; perhaps the main technical subtlety is that we would need to invoke local laws dealing with block matrices, such as those of \cite{erdHos2025cusp}.

\section*{Acknowledgments}
\addcontentsline{toc}{section}{Acknowledgments}

We thank Afonso Bandeira, Tatiana Brailovskaya, and Ke Wang for helpful suggestions during the course of this project.

\bibliographystyle{alpha}
\addcontentsline{toc}{section}{References}
\bibliography{bib}

\newcommand{\etalchar}[1]{$^{#1}$}
\begin{thebibliography}{BDWW22}

\bibitem[AEK17]{AEK-2017-UniversalityWignerTypeMatrices}
Oskari~H Ajanki, L{\'a}szl{\'o} Erd{\H{o}}s, and Torben Kr{\"u}ger.
\newblock Universality for general {Wigner}-type matrices.
\newblock {\em Probability Theory and Related Fields}, 169(3):667--727, 2017.

\bibitem[AFWZ20]{abbe2019entrywiseeigenvectoranalysisrandom}
Emmanuel Abbe, Jianqing Fan, Kaizheng Wang, and Yiqiao Zhong.
\newblock Entrywise eigenvector analysis of random matrices with low expected
  rank.
\newblock {\em The Annals of Statistics}, 48(3):1452--1474, 2020.

\bibitem[AGZ09]{Anderson_Guionnet_Zeitouni_2009}
Greg~W. Anderson, Alice Guionnet, and Ofer Zeitouni.
\newblock {\em An Introduction to Random Matrices}, volume 118 of {\em
  Cambridge Studies in Advanced Mathematics}.
\newblock Cambridge University Press, 2009.

\bibitem[BBAP05]{10.1214/009117905000000233}
Jinho Baik, G{\'e}rard Ben~Arous, and Sandrine P{\'e}ch{\'e}.
\newblock Phase transition of the largest eigenvalue for nonnull complex sample
  covariance matrices.
\newblock {\em The Annals of Probability}, 33(5):1643--1697, 2005.

\bibitem[BCH26]{bhattacharya2026outliereigenvalueseigenvectors}
Bishakh Bhattacharya, Arijit Chakrabarty, and Rajat~Subhra Hazra.
\newblock Outlier eigenvalues and eigenvectors of generalized {Wigner} matrices
  with finite-rank perturbations, 2026.

\bibitem[BDW21]{bao2020singularvectorsingularsubspace}
Zhigang Bao, Xiucai Ding, and Ke~Wang.
\newblock Singular vector and singular subspace distribution for the matrix
  denoising model.
\newblock {\em The Annals of Statistics}, 49(1):370--392, 2021.

\bibitem[BDWW22]{bao2020statisticalinferenceprincipalcomponents}
Zhigang Bao, Xiucai Ding, Jingming Wang, and Ke~Wang.
\newblock Statistical inference for principal components of spiked covariance
  matrices.
\newblock {\em The Annals of Statistics}, 50(2):1144--1169, 2022.

\bibitem[BEK{\etalchar{+}}14]{bloemendal2015isotropiclocallawssample}
Alex Bloemendal, L{\'a}szl{\'o} Erd{\H{o}}s, Antti Knowles, Horng-Tzer Yau, and
  Jun Yin.
\newblock Isotropic local laws for sample covariance and generalized {Wigner}
  matrices.
\newblock {\em Electronic Journal of Probability}, 19(33):1--53, 2014.

\bibitem[BGGM11]{benaychgeorges2011fluctuationsextremeeigenvaluesfinite}
Florent Benaych-Georges, Alice Guionnet, and Myl{\`e}ne Ma{\"i}da.
\newblock Fluctuations of the extreme eigenvalues of finite rank deformations
  of random matrices.
\newblock {\em Electronic Journal of Probability}, 16(60):1621--1662, 2011.

\bibitem[BGK16]{BGK-2016-LecturesLocalSemicircleLaw}
Florent Benaych-Georges and Antti Knowles.
\newblock Lectures on the local semicircle law for {Wigner} matrices, 2016.

\bibitem[BGN11]{benaychgeorges2010eigenvalueseigenvectorsfinitelow}
Florent Benaych-Georges and Raj~Rao Nadakuditi.
\newblock The eigenvalues and eigenvectors of finite, low rank perturbations of
  large random matrices.
\newblock {\em Advances in Mathematics}, 227(1):494--521, 2011.

\bibitem[BvH16]{10.1214/15-AOP1025}
Afonso~S. Bandeira and Ramon van Handel.
\newblock Sharp nonasymptotic bounds on the norm of random matrices with
  independent entries.
\newblock {\em The Annals of Probability}, 44(4):2479--2506, 2016.

\bibitem[Cap17]{Capitaine-2017-HDRRandomMatrices}
Mireille Capitaine.
\newblock {\em Deformed ensembles, polynomials in random matrices and free
  probability theory}.
\newblock Habilitation {\`a} diriger des recherches, Universit{\'e} Paul
  Sabatier -- Toulouse 3, 2017.
\newblock HAL Id: tel-01978065, version 1.

\bibitem[CDM21]{capitaine2020nonuniversalityfluctuationsoutlier}
Mireille Capitaine and Catherine Donati-Martin.
\newblock Non universality of fluctuations of outlier eigenvectors for block
  diagonal deformations of {Wigner} matrices.
\newblock {\em ALEA. Latin American Journal of Probability and Mathematical
  Statistics}, 18:129--165, 2021.

\bibitem[CDMF09]{Capitaine_2009}
Mireille Capitaine, Catherine Donati-Martin, and Delphine F{\'e}ral.
\newblock The largest eigenvalues of finite rank deformation of large {Wigner}
  matrices: Convergence and nonuniversality of the fluctuations.
\newblock {\em The Annals of Probability}, 37(1):1--47, January 2009.

\bibitem[CDMF12]{capitaine2011centrallimittheoremseigenvalues}
Mireille Capitaine, Catherine Donati-Martin, and Delphine F{\'e}ral.
\newblock Central limit theorems for eigenvalues of deformations of {Wigner}
  matrices.
\newblock {\em Annales de l'Institut Henri Poincar{\'e}, Probabilit{\'e}s et
  Statistiques}, 48(1):107--133, 2012.

\bibitem[CH13]{couillet2012fluctuationsspikedrandommatrix}
Romain Couillet and Walid Hachem.
\newblock Fluctuations of spiked random matrix models and failure diagnosis in
  sensor networks.
\newblock {\em IEEE Transactions on Information Theory}, 59(1):509--525, 2013.

\bibitem[CR24]{cademartori2023nonasymptoticanalysisgeneralizedapproximate}
Collin Cademartori and Cynthia Rush.
\newblock A non-asymptotic analysis of generalized vector approximate message
  passing algorithms with rotationally invariant designs.
\newblock {\em IEEE Transactions on Information Theory}, 70(8):5811--5856,
  August 2024.

\bibitem[CSC12]{Cucuringu_2012}
Mihai Cucuringu, Amit Singer, and David Cowburn.
\newblock Eigenvector synchronization, graph rigidity and the molecule problem.
\newblock {\em Information and Inference: A Journal of the IMA}, 1(1):21--67,
  December 2012.

\bibitem[CT22]{cucuringu2021extensionangularsynchronizationproblem}
Mihai Cucuringu and Hemant Tyagi.
\newblock An extension of the angular synchronization problem to the
  heterogeneous setting.
\newblock {\em Foundations of Data Science}, 4(1):71--122, 2022.

\bibitem[DK70]{davis1970rotation}
Chandler Davis and W.~M. Kahan.
\newblock The rotation of eigenvectors by a perturbation. {III}.
\newblock {\em SIAM Journal on Numerical Analysis}, 7(1):1--46, March 1970.

\bibitem[EHR25]{erdHos2025cusp}
L{\'a}szl{\'o} Erd{\H{o}}s, Sven~Joscha Henheik, and Volodymyr Riabov.
\newblock Cusp universality for correlated random matrices.
\newblock {\em Communications in Mathematical Physics}, 406(10):253, 2025.

\bibitem[EYY12]{erdos2011rigidityeigenvaluesgeneralizedwigner}
L{\'a}szl{\'o} Erd{\H{o}}s, Horng-Tzer Yau, and Jun Yin.
\newblock Rigidity of eigenvalues of generalized {Wigner} matrices.
\newblock {\em Advances in Mathematics}, 229(3):1435--1515, 2012.

\bibitem[FFHL22]{fan2020asymptotictheoryeigenvectorsrandom}
Jianqing Fan, Yingying Fan, Xiao Han, and Jinchi Lv.
\newblock Asymptotic theory of eigenvectors for random matrices with diverging
  spikes.
\newblock {\em Journal of the American Statistical Association},
  117(538):996--1009, 2022.

\bibitem[FK81]{FurediKomlos1981Eigenvalues}
Zolt{\'a}n F{\"u}redi and J{\'a}nos Koml{\'o}s.
\newblock The eigenvalues of random symmetric matrices.
\newblock {\em Combinatorica}, 1(3):233--241, September 1981.

\bibitem[FP07]{F_ral_2007}
Delphine F{\'e}ral and Sandrine P{\'e}ch{\'e}.
\newblock The largest eigenvalue of rank one deformation of large {Wigner}
  matrices.
\newblock {\em Communications in Mathematical Physics}, 272(1):185--228, March
  2007.

\bibitem[FVRS22]{feng2021unifyingtutorialapproximatemessage}
Oliver~Y. Feng, Ramji Venkataramanan, Cynthia Rush, and Richard~J. Samworth.
\newblock A unifying tutorial on approximate message passing.
\newblock {\em Foundations and Trends in Machine Learning}, 15(4):335--536,
  2022.

\bibitem[GZ19]{gao2019multifrequencyphasesynchronization}
Tingran Gao and Zhizhen Zhao.
\newblock Multi-frequency phase synchronization.
\newblock In Kamalika Chaudhuri and Ruslan Salakhutdinov, editors, {\em
  Proceedings of the 36th International Conference on Machine Learning},
  volume~97 of {\em Proceedings of Machine Learning Research}, pages
  2132--2141. PMLR, 09--15 Jun 2019.

\bibitem[Han25]{han2025entrywisedynamicsuniversalitygeneral}
Qiyang Han.
\newblock Entrywise dynamics and universality of general first order methods.
\newblock {\em The Annals of Statistics}, 53(4):1783--1807, 2025.

\bibitem[Ide16]{idel2016reviewmatrixscalingsinkhorns}
Martin Idel.
\newblock A review of matrix scaling and {Sinkhorn}'s normal form for matrices
  and positive maps, 2016.

\bibitem[Joh01]{10.1214/aos/1009210544}
Iain~M. Johnstone.
\newblock On the distribution of the largest eigenvalue in principal components
  analysis.
\newblock {\em The Annals of Statistics}, 29(2):295--327, 2001.

\bibitem[JP18]{JohnstonePaul2018PCAHighDim}
Iain~M. Johnstone and Debashis Paul.
\newblock {PCA} in high dimensions: An orientation.
\newblock {\em Proceedings of the IEEE}, 106(8):1277--1292, 2018.

\bibitem[KBK26]{kireeva2024computationallowerboundsmultifrequency}
Anastasia Kireeva, Afonso~S. Bandeira, and Dmitriy Kunisky.
\newblock Computational lower bounds for multi-frequency group synchronization.
\newblock {\em Applied and Computational Harmonic Analysis}, 84:101880, 2026.

\bibitem[Kun25]{Kunisky-2025-LCDA2}
Dmitriy Kunisky.
\newblock Low coordinate degree algorithms {II}: Categorical signals and
  generalized stochastic block models.
\newblock In Nika Haghtalab and Ankur Moitra, editors, {\em Proceedings of the
  Thirty-Eighth Conference on Learning Theory}, volume 291 of {\em Proceedings
  of Machine Learning Research}, pages 3486--3526. PMLR, 2025.

\bibitem[KY13a]{knowles2011eigenvectordistributionwignermatrices}
Antti Knowles and Jun Yin.
\newblock Eigenvector distribution of {Wigner} matrices.
\newblock {\em Probability Theory and Related Fields}, 155(3--4):543--582,
  2013.

\bibitem[KY13b]{knowles2012isotropicsemicirclelawdeformation}
Antti Knowles and Jun Yin.
\newblock The isotropic semicircle law and deformation of {Wigner} matrices.
\newblock {\em Communications on Pure and Applied Mathematics},
  66(11):1663--1750, 2013.

\bibitem[KY14]{Knowles_2014}
Antti Knowles and Jun Yin.
\newblock The outliers of a deformed {Wigner} matrix.
\newblock {\em The Annals of Probability}, 42(5):1980--2031, September 2014.

\bibitem[KY17]{knowles2016anisotropiclocallawsrandom}
Antti Knowles and Jun Yin.
\newblock Anisotropic local laws for random matrices.
\newblock {\em Probability Theory and Related Fields}, 169(1--2):257--352,
  2017.

\bibitem[LCC24]{lebeau2024asymptoticgaussianfluctuationseigenvectors}
Hugo Lebeau, Florent Chatelain, and Romain Couillet.
\newblock Asymptotic gaussian fluctuations of eigenvectors in spectral
  clustering.
\newblock {\em IEEE Signal Processing Letters}, 31:1920--1924, 2024.

\bibitem[LFW23]{li2023approximatemessagepassingrandom}
Gen Li, Wei Fan, and Yuting Wei.
\newblock Approximate message passing from random initialization with
  applications to $\mathbb{Z}_{2}$ synchronization.
\newblock {\em Proceedings of the National Academy of Sciences},
  120(31):e2302930120, 2023.

\bibitem[Li26]{Li-2026-LowerBoundMultiFrequencySynchronization}
Zhangsong Li.
\newblock Improved computational lower bound of estimation for multi-frequency
  group synchronization.
\newblock {\em arXiv preprint arXiv:2601.20522}, 2026.

\bibitem[LM00]{10.1214/aos/1015957395}
B.~Laurent and P.~Massart.
\newblock Adaptive estimation of a quadratic functional by model selection.
\newblock {\em The Annals of Statistics}, 28(5):1302--1338, 2000.

\bibitem[LW23]{li2023nonasymptoticframeworkapproximatemessage}
Gen Li and Yuting Wei.
\newblock A non-asymptotic framework for approximate message passing in spiked
  models, 2023.

\bibitem[MTV22]{mondelli2021optimalcombinationlinearspectral}
Marco Mondelli, Christos Thrampoulidis, and Ramji Venkataramanan.
\newblock Optimal combination of linear and spectral estimators for generalized
  linear models.
\newblock {\em Foundations of Computational Mathematics}, 22(5):1513--1566,
  2022.

\bibitem[MV21]{montanari2019estimationlowrankmatricesapproximate}
Andrea Montanari and Ramji Venkataramanan.
\newblock Estimation of low-rank matrices via approximate message passing.
\newblock {\em The Annals of Statistics}, 49(1):321--345, 2021.

\bibitem[MV22]{mondelli2021approximatemessagepassingspectral}
Marco Mondelli and Ramji Venkataramanan.
\newblock Approximate message passing with spectral initialization for
  generalized linear models.
\newblock {\em Journal of Statistical Mechanics: Theory and Experiment},
  2022(11):114003, 2022.

\bibitem[MY22]{marcinek2020highdimensionalnormalitynoisy}
Jake Marcinek and Horng-Tzer Yau.
\newblock High dimensional normality of noisy eigenvectors.
\newblock {\em Communications in Mathematical Physics}, 395(3):1007--1096,
  2022.

\bibitem[PA14]{PAUL20141}
Debashis Paul and Alexander Aue.
\newblock Random matrix theory in statistics: A review.
\newblock {\em Journal of Statistical Planning and Inference}, 150:1--29, 2014.

\bibitem[P{\'e}c06]{Peche2006LargestEigenvalue}
Sandrine P{\'e}ch{\'e}.
\newblock The largest eigenvalue of small rank perturbations of {Hermitian}
  random matrices.
\newblock {\em Probability Theory and Related Fields}, 134(1):127--173, January
  2006.

\bibitem[PWBM16]{perry2016optimalitysuboptimalitypcaspiked}
Amelia Perry, Alexander~S. Wein, Afonso~S. Bandeira, and Ankur Moitra.
\newblock Optimality and sub-optimality of {PCA} for spiked random matrices and
  synchronization, 2016.

\bibitem[PWBM18]{Perry_2018}
Amelia Perry, Alexander~S. Wein, Afonso~S. Bandeira, and Ankur Moitra.
\newblock Message-passing algorithms for synchronization problems over compact
  groups.
\newblock {\em Communications on Pure and Applied Mathematics},
  71(11):2275--2322, April 2018.

\bibitem[RG20]{Romanov_2020}
Elad Romanov and Matan Gavish.
\newblock The noise-sensitivity phase transition in spectral group
  synchronization over compact groups.
\newblock {\em Applied and Computational Harmonic Analysis}, 49(3):935--970,
  November 2020.

\bibitem[RS13]{renfrew2012finiterankdeformationswigner}
David Renfrew and Alexander Soshnikov.
\newblock On finite rank deformations of {Wigner} matrices {II}: Delocalized
  perturbations.
\newblock {\em Random Matrices: Theory and Applications}, 2(1):1250015, 2013.

\bibitem[RV18]{Rush_2018}
Cynthia Rush and Ramji Venkataramanan.
\newblock Finite sample analysis of approximate message passing algorithms.
\newblock {\em IEEE Transactions on Information Theory}, 64(11):7264--7286,
  November 2018.

\bibitem[Sin11]{singer2009angularsynchronizationeigenvectorssemidefinite}
Amit Singer.
\newblock Angular synchronization by eigenvectors and semidefinite programming.
\newblock {\em Applied and Computational Harmonic Analysis}, 30(1):20--36,
  2011.

\bibitem[YWF25]{YWF-2024-MutualInformationQuadraticCompactGroups}
Kaylee~Y. Yang, Timothy L.~H. Wee, and Zhou Fan.
\newblock Asymptotic mutual information in quadratic estimation problems over
  compact groups.
\newblock {\em Information and Inference: A Journal of the IMA}, 14(3):iaaf024,
  2025.

\bibitem[YWS15]{yu2014usefulvariantdaviskahantheorem}
Yi~Yu, Tengyao Wang, and Richard~J. Samworth.
\newblock A useful variant of the {Davis}--{Kahan} theorem for statisticians.
\newblock {\em Biometrika}, 102(2):315--323, 2015.

\end{thebibliography}

\end{document}